\documentclass[12pt]{article}

\usepackage{indentfirst}
\usepackage{amsfonts}
\usepackage{amsmath}
\usepackage{amssymb}
\usepackage{amsbsy}
\usepackage[margin=1in]{geometry}
\usepackage{tikz}
\usepackage{circuitikz}
\usepackage{pgflibrarysnakes}
\usetikzlibrary{snakes}

\newtheorem{dfn}{Definition}[section]
\newtheorem{theorem}[dfn]{Theorem}
\newtheorem{lemma}[dfn]{Lemma}

\newtheorem{corollary}[dfn]{Corollary}
\newenvironment{pf}{\noindent{\bf Proof.}}
{\enspace\vrule height5pt depth0pt width5pt} 

\def\F {{\mathcal F}}
\def\Se {{\mathcal S}}
\def\A {{\mathcal A}}
\def\P {{\mathcal P}}

\begin{document}

\title{Well-quasi-ordering digraphs with no long alternating paths by the strong immersion relation}
\author{Chun-Hung Liu\thanks{chliu@math.tamu.edu. Partially supported by NSF under Grant No.~DMS-1929851 and DMS-1954054.} \\
\small Department of Mathematics, \\
\small Texas A\&M University,\\
\small College Station, TX, USA \\
\\
Irene Muzi\thanks{irene.muzi@gmail.com. Supported by the European Research Council (ERC) under the European Union's Horizon 2020 research and innovation programme (ERC Consolidator Grant DISTRUCT, grant agreement No.~648527).} \\
\small Institut fur Softwaretechnik und Theoretische Informatik \\
\small (Logic and Semantics Research Group) \\
\small Technical University Berlin \\
\small Berlin, Germany}

\maketitle

\begin{abstract}
Nash-Williams' Strong Immersion Conjecture states that graphs are well-quasi-ordered by the strong immersion relation.
That is, given infinitely many graphs, one graph contains another graph as a strong immersion.
In this paper we study the analogous problem for directed graphs.
It is known that digraphs are not well-quasi-ordered by the strong immersion relation, but for all known such infinite antichains, paths that change direction arbitrarily many times can be found.
This paper proves that the converse statement is true: for every positive integer $k$, the digraphs that do not contain a path that changes direction $k$ times are well-quasi-ordered by the strong immersion relation, even when vertices are labelled by a well-quasi-order.
This result is optimal for classes of digraphs closed under taking subgraphs since paths that change direction arbitrarily many times with vertex-labels form an infinite antichain with respect to the strong immersion relation.
\end{abstract}

\section{Introduction}

In this paper, graphs and directed graphs are finite, loopless and allowed to have parallel edges, unless otherwise specified.

A {\it quasi-ordering} is a reflexive and transitive binary relation.
A quasi-ordering $\preceq$ on a set $S$ is a {\it well-quasi-ordering} if for every infinite sequence $a_1,a_2,...$ over $S$, there exist $1 \leq i <j$ such that $a_i \preceq a_j$.
We say that $Q=(S,\preceq)$ is a {\it quasi-order} (or a {\it well-quasi-order}, respectively) if $\preceq$ is a quasi-ordering (or a well-quasi-ordering, respectively) on $S$.

The study of well-quasi-ordering on graphs can be traced back to a conjecture of V\'{a}zsonyi proposed in 1940s: Trees are well-quasi-ordered by the topological minor relation.
We say that a graph $G$ contains another graph $H$ as a {\it topological minor} if some subgraph of $G$ is isomorphic to a subdivision of $H$.
This conjecture was proved by Kruskal \cite{k_tree} and independently by Tarkowski \cite{t_wqo_tree}. 
Another conjecture proposed by V\'{a}zsonyi states that subcubic graphs are well-quasi-ordered by the topological minor relation.
This conjecture is significantly more difficult than the previous conjecture on trees.
The only known proof of this conjecture is via the celebrated Graph Minor Theorem of Robertson and Seymour \cite{rs XX}.

A graph $G$ contains another graph $H$ as a {\it minor} if $H$ is isomorphic to a graph that can be obtained from a subgraph of $G$ by repeatedly contracting edges.
The Graph Minor Theorem \cite{rs XX} states that graphs are well-quasi-ordered by the minor relation.
It is one of the deepest theorems in graph theory, and its proof spans over around 20 papers.
As for subcubic graphs, the minor relation is equivalent to the topological minor relation.
The aforementioned conjecture of V\'{a}zsonyi on subcubic graphs is then an immediate corollary of the Graph Minor Theorem.

One strength of well-quasi-ordering is an implication of the existence of a finite characterization of a property that is closed under a well-quasi-ordering.
A property is {\it closed under a quasi-ordering $\preceq$} if an element $a$ satisfies this property implies that every element $b$ with $b \preceq a$ also satisfies this property.
For a property $\P$ that is closed under a quasi-ordering $\preceq$, we denote the set of $\preceq$-minimal elements that do not satisfy $\P$ by $m(\P)$.
If $\preceq$ is a well-quasi-ordering, then $m(\P)$ is finite since $m(\P)$ is an antichain with respect to $\preceq$. 
Note that one can determine whether a given input element $x$ satisfies $\P$ or not by testing whether $y \preceq x$ for every $y \in m(\P)$.
So $\P$ is uniquely determined by $m(\P)$.
And if $\lvert m(\P) \rvert$ is finite, and if for each fixed $y \in m(\P)$, testing whether an input element $x$ satisfies $y \preceq x$ or not can be done in polynomial time, then one can decide whether $x$ satisfies $\P$ or not in polynomial time.

This leads to prominent applications of the Graph Minor Theorem.
It implies that every minor-closed property (such as the embeddability in any fixed surface, linkless embeddability or knotless embeddability etc.) can be characterized by finitely many graphs.
As Robertson and Seymour \cite{rs XIII} also proved that for any fixed graph $H$, deciding whether an input graph $G$ contains $H$ as a minor or not can be done in polynomial time, the discussion in the previous paragraph implies that every minor-closed property can be decided in polynomial time.

Due to the power of well-quasi-ordering and the success of the Graph Minor Theorem, one might consider whether the Graph Minor Theorem can be generalized.
One possible generalization would be to extend the result to infinite graphs.
However, it was disproved by Thomas \cite{t_counter}.
But it remains open whether the Graph Minor Theorem is true for countable graphs.

Another possible generalization is to extend the Graph Minor Theorem to relations that are finer than the minor relation.
The topological minor relation is an example, as if a graph $G$ contains another graph $H$ as a topological minor, then $G$ contains $H$ as a minor.
Hence one might ask whether V\'{a}zsonyi's two conjectures on topological minors can be extended to all graphs.
However, it is not true, and there are many different constructions for infinite antichains with respect to the topological minor relation.
Robertson conjectured a common generalization of the two V\'{a}zsonyi's conjecture in 1980s.
This conjecture was proved by the first author and Thomas \cite{l_thesis}.
See \cite{l_thesis,lp} for more details.

Though the topological minor relation does not well-quasi-order all graphs, it is still possible to extend both V\'{a}zsonyi's conjectures to all graphs.
Nash-Williams conjectured that the weak immersion relation \cite{n_weak} and the strong immersion relation \cite{n_strong} are well-quasi-ordering on graphs.
Both of these conjectures imply both V\'{a}zsonyi's conjectures as for trees and subcubic graphs, the weak and strong immersion relations are equivalent to the topological minor relation.
Nash-Williams' Weak Immersion Conjecture was proved by Robertson and Seymour \cite{rs XXIII} by strengthening the Graph Minor Theorem.
The Strong Immersion conjecture remains open\footnote{Robertson and Seymour believe that they had a proof of the Strong Immersion Conjecture at one time, but even if it was correct, it was very complicated, and it is unlikely that they will write it down (see \cite{rs XXIII}).}, and progress on it seems rare in the literature.
It can be easily shown that the conjecture is true for graphs of bounded maximum degree by using the result on weak immersion (see \cite{l_survey}).
Andreae \cite{a_strong} proved the conjecture for the class of simple graphs satisfying that either they do not contain $K_{2,3}$ as a strong immersion, or all blocks are either complete graphs, cycles, or balanced complete bipartite graphs.

Another possible extension of the Graph Minor Theorem is to consider directed graphs.
There are different notions of minors for directed graphs.
We only consider butterfly minors here. 
(See \cite{l_survey} for a survey about well-quasi-ordering on graphs, including results for different minor containments for digraphs.)
Again, the butterfly minor relation does not well-quasi-order all digraphs.
Every construction of infinite antichains involves paths that change direction arbitrarily many times.
Chudnovsky, the second author, Oum, Seymour and Wollan (see \cite{m}) proved that this obstruction is the only obstruction: for every positive integer $k$, digraphs whose underlying graphs do not contain a path that change direction $k$ times are well-quasi-ordered by the butterfly minor relation.

This paper addresses a combination of two directions mentioned above: we consider the strong immersion relation on digraphs.
We need some notions to formally state our result.

Let $G$ and $H$ be digraphs possibly with loops.
A function $f$ is a {\it strong immersion embedding} from $H$ to $G$ if the following hold.
	\begin{itemize}
		\item $f$ maps $V(H)$ to $V(G)$ injectively.
		\item $f$ maps each non-loop edge of $H$ with tail $x$ and head $y$ to a directed path in $G$ from $f(x)$ to $f(y)$; $f$ maps each loop of $H$ with end $x$ to a directed cycle passing through $f(x)$.
		\item If $e_1,e_2$ are different edges of $H$, then $f(e_1)$ and $f(e_2)$ are edge-disjoint.
		\item For every edge $e$ of $H$ and every vertex $v$ of $H$, if $v$ is not an end of $e$, then $f(v) \not \in V(f(e))$.
	\end{itemize}
We say that {\it $G$ contains $H$ as a strong immersion} if there exists a strong immersion embedding from $H$ to $G$.

The strong immersion relation does not well-quasi-order digraphs.
A {\it thread} is a digraph whose underlying graph is a path.
A {\it pivot} in a thread is a vertex that has either in-degree two or out-degree two.
For an integer $k$, a {\it $k$-alternating path} is a thread that contains exactly $k$ pivots.
If for each $i \in {\mathbb N}$, $G_i$ is the digraph obtained from an $i$-alternating path by attaching two leaves to each end of the path, then it is easy to see that $\{G_i: i \in {\mathbb N}\}$ forms an infinite antichain with respect to the strong immersion relation.

In fact, the alternating paths already form an infinite antichain if vertices are allowed to be labelled.
The main result of this paper proves the converse statement: forbidding long alternating paths is sufficient to ensure well-quasi-ordering even when vertices are labelled.

\begin{theorem} \label{main_label}
Let $k$ be a positive integer.
Let $(Q,\leq_Q)$ be a well-quasi-order.
For every $i \in {\mathbb N}$, let $D_i$ be a digraph with loops allowed and with no $k$-alternating path, and let $\phi_i: V(D_i) \rightarrow Q$ be a function.
Then there exist $1 \leq j < j'$ and a strong immersion embedding $\eta$ from $D_j$ to $D_{j'}$ such that for every $v \in V(D_j)$, $\phi_j(v) \leq_Q \phi_{j'}(\eta(v))$.
\end{theorem}

Note that Theorem \ref{main_label} is optimal when $Q$ contains two non-equivalent elements.
That is, it is the case when there exist elements $x$ and $y$ of $Q$ with $x \not \leq_Q y$.
For every $i \in {\mathbb N}$, let $D_i$ be an $i$-alternating path, and let $\phi_i$ be the function that maps the ends of $D_i$ to $x$ and maps all other vertices to $y$.
Clearly, there exist no strong immersion embedding from $D_i$ to $D_j$ preserving the labels on the vertices, for any $i \neq j$.

We remark that even though Theorem \ref{main_label} is optimal, it is known that some class of digraphs with arbitrarily long alternating paths are well-quasi-ordered by the strong immersion relation.
For example, Chudnovsky and Seymour \cite{cs} proved that tournaments are well-quasi-ordered by the strong immersion relation.
Note that the class of digraphs in Theorem \ref{main_label} is closed under taking subgraphs, but the class of tournaments is not.
The result for tournaments was extended to semicomplete digraphs by Barbero, Paul and Pilipczuk \cite{bpp}.

\subsection{Organization of the paper}

We shall prove Theorem \ref{main_label} by induction on $k$.
The proof of Theorem \ref{main_label} uses a strengthening of an idea in the work of Chudnovsky, the second author, Oum, Seymour and Wollan (see \cite{m}) for butterfly minors.
Roughly speaking, it shows that if $D$ is a digraph with no $k$-alternating path, then one can delete at most $f(k)$ vertices to kill all $(k-1)$-alternating paths in $D$ not contained in a ``series-parallel digraph with two roots''.
This suggests that we have to prove well-quasi-ordering results on those series-parallel digraphs with two roots with respect to the strong immersion relation preserving the roots.
In general, proving well-quasi-ordering for strong immersion preserving certain ``roots'' is required in many circumstances of this paper. 
It is significantly more complicated and requires more tricks than the analogous work for butterfly minors, even when dealing with the case of series-parallel digraphs.

This paper is organized as follows.
In Section \ref{sec:prelim wqo} we review some well-known results about well-quasi-ordering that will be used in this paper.
In Section \ref{sec:1-alt} we prove the case $k=1$ of Theorem \ref{main_label}.
We introduce the notion of series-parallel triples in Section \ref{sec:sp_triple}.
It is the formal form of the ``series-parallel digraphs with two roots'' mentioned above.
In Section \ref{sec:sp_tree}, we prove a well-quasi-ordering result for the tree-like digraphs formed by repeatedly gluing those series-parallel triples.
It is a crucial step toward the result for well-quasi-ordering series-parallel triples preserving roots which will be proved in Section \ref{sec:wqo_sp_triples}.
Then in Section \ref{sec:sp_sep}, we prove the tools that allow us to kill all $(k-1)$-alternating paths not hidden in a series-parallel triple mentioned above, and study the relationships between all series-parallel triples. 
Finally, we prove Theorem \ref{main_label} in Section \ref{sec:long_alt}.

\subsection{Notation}
For a graph (or digraph, respectively) $G$ and a subset $S$ of $V(G)$, we denote $G[S]$ by the graph (or digraph, respectively) induced on $S$; if $T \subseteq V(G)$, then $G-T$ is defined to be $G[V(G)-T]$; for a vertex $v$ of $G$, $G-v$ is defined to be $G-\{v\}$.

Let $f$ be a function with domain $X$.
If $S$ is a subset of $X$, then $f(S)$ is defined to be the set $\{f(s): s \in S\}$.
If $S$ is a sequence $(s_1,s_2,...,s_k)$ over $X$, then $f(S)$ is defined to be the sequence $(f(s_1),f(s_2),...,f(s_k))$.

For every positive integer $k$, we define $[k]$ to be the set $\{1,2,...,k\}$.

\section{Preliminary about well-quasi-ordering} \label{sec:prelim wqo}

In this section, we review some known useful tools about well-quasi-ordering.

Let $(Q_1,\preceq_1)$ and $(Q_2,\preceq_2)$ be well-quasi-orders. 
We say that $(Q,\preceq)$ is the {\it well-quasi-order obtained by taking the disjoint union of $(Q_1,\preceq_1)$ and $(Q_2,\preceq_2)$} if $Q$ is a disjoint union of a copy of $Q_1$ and a copy of $Q_2$ such that for $x,y \in Q$, $x \preceq y$ if and only if either $x,y \in Q_1$ with $x \preceq_1 y$, or $x,y \in Q_2$ with $x \preceq_2 y$.
We say that $(Q',\preceq')$ is the {\it well-quasi-order obtained by the Cartesian product of $(Q_1,\preceq_1)$ and $(Q_2,\preceq_2)$} if $Q'=Q_1 \times Q_2$ such that for $(x_1,y_1),(x_2,y_2) \in Q'$, $(x_1,y_1) \preceq' (x_2,y_2)$ if and only if $x_1 \preceq_1 x_2$ and $y_1 \preceq_2 y_2$.

Let $(Q,\leq_Q)$ be a well-quasi-order.
We say that the $(Q',\preceq')$ is the {\it well-quasi-order obtained from $(Q,\leq_Q)$ by Higman's Lemma} if $Q'$ is the set of finite sequences over $Q$ such that for elements $(a_1,a_2,...,a_m)$ and $(b_1,b_2,...,b_n)$ of $Q'$, $(a_1,a_2,...,a_m) \preceq' (b_1,b_2,...,b_n)$ if and only if there exists a strictly increasing function $\iota: [m] \rightarrow [n]$ such that $a_i \leq_Q b_{\iota(i)}$ for every $i \in [m]$.
Note that $(Q',\preceq')$ is indeed a well-quasi-order, as shown by a famous result of Higman \cite{h}.

Another known result that we will use in this paper is a strengthening of Kruskal's Tree Theorem proved by Kriz \cite{k}.
We need the following definition to formally state the theorem.

A {\it homeomorphic embedding} from a digraph $H$ possibly with loops to a digraph $G$ possibly with loops is a function $\eta$ satisfying the following.
	\begin{itemize}
		\item $\eta$ maps $V(H)$ to $V(G)$ injectively.
		\item $\eta$ maps each loop of $H$ with end $v$ to a directed cycle of $G$ passing through $\eta(v)$; $\eta$ maps each non-loop edge $H$ with tail $x$ and head $y$ to a directed path in $G$ from $\eta(x)$ to $\eta(y)$.
		\item If $e_1,e_2$ are distinct edges, then $\eta(e_1) \cap \eta(e_2) = \eta(e_1 \cap e_2)$.
		\item If a vertex $v$ of $H$ is not incident with an edge $e$ of $H$, then $\eta(v) \not \in \eta(e)$.
	\end{itemize}
A {\it rooted tree} is a directed graph whose underlying graph is a tree such that all but exactly one vertex have in-degree one. 
We denote the first infinite ordinal number by $\omega$.
We will only need the following special case of Kriz's theorem.

\begin{theorem}[\cite{k}] \label{wqo gap}
Let $(Q, \preceq)$ be a well-quasi-order.
For each positive integer $i$, let $T_i$ be a rooted tree, $\phi_i: V(T_i) \rightarrow Q$ and $\mu_i: E(T_i) \rightarrow {\mathbb N} \cup \{0,\omega\}$.
Then there exist $1 \leq i < j$ such that there exists a homeomorphic embedding $\eta$ from $T_i$ to $T_j$ such that the following hold.
	\begin{enumerate}
		\item For every $v \in V(T_i)$, $\phi_i(v) \preceq \phi_j(\eta(v))$.
		\item For every $e \in E(T_i)$, if $f$ is an edge in $\eta(e)$, then $\mu_i(e) \leq \mu_j(f)$.
	\end{enumerate}
\end{theorem}

\section{1-alternating paths} \label{sec:1-alt}

For a digraph $D$, a {\it source} in $D$ is a vertex of in-degree 0, and a {\it sink} in $D$ is a vertex of out-degree 0.

\begin{lemma} \label{1-alt structure}
Let $D$ be a digraph whose underlying graph is connected.
If $D$ has no 1-alternating path, then either $\lvert V(D) \rvert \leq 2$, or $D$ is obtained by a directed path or a directed cycle by duplicating edges arbitrarily many times.
\end{lemma}

\begin{pf}
We may assume that $D$ contains at least three vertices, for otherwise we are done.
Let $P$ be a thread in $D$ with maximum length.
Since the underlying graph of $D$ is connected and has at least three vertices, $P$ contains at least three vertices.
Denote $P$ by $v_1v_2...v_k$, where $k=\lvert V(P) \rvert$.
By symmetry, we may assume that $v_1$ is a source of $P$.

By the maximality of $P$, $v_1$ and $v_k$ have no neighbor in $D$ not contained in $P$.
Since $D$ has no 1-alternating path, $P$ is a directed path, and $v_i$ has no neighbor in $D$ not contained in $P$ for every $2 \leq i \leq k-1$.
Hence $P$ contains all vertices of $D$.

Let $e \in E(D)-E(P)$ with tail $v_i$ and head $v_j$.
Since $D$ is loopless, $i \neq j$.
Since $D$ has no 1-alternating path, either $j=i+1$, or $(i,j)=(k,1)$.
This proves the lemma.
\end{pf}

\begin{lemma} \label{1-alt wqo}
Let $(Q,\preceq)$ be a well-quasi-order.
For each $i \in {\mathbb N}$, let $D_i$ be a directed graph with no 1-alternating path, and let $\phi_i:V(D_i) \rightarrow Q$.
Then there exist $1 \leq i < j$ and a strong immersion embedding $\eta$ from $D_i$ to $D_j$ such that $\phi_i(v) \preceq \phi_j(\eta(v))$ for every $v \in V(D_i)$.
\end{lemma}

\begin{pf}
By Lemma \ref{1-alt structure}, each $D_i$ either contains at most two vertices or can be obtained from a directed path or a directed cycle by duplicating edges arbitrarily many times.
It is easy if there are infinitely many indices $i$ such that $D_i$ containing at most two vertices.
So we may assume that every $D_i$ contains at least three vertices, and either every $D_i$ is obtained from a directed path by duplicating edges arbitrarily many times, or every $D_i$ is obtained from a directed cycle by duplicating edges arbitrarily many times.

For each $i$, let $W_i$ be a Hamiltonian directed path of $D_i$.
For each $i$, let $x_i,y_i$ be the ends of $W_i$ such that $W_i$ is from $x_i$ to $y_i$, and let $\ell_i$ be the number of directed edges in $D_i$ between $y_i$ and $x_i$.
Let $(Q_1,\preceq_1)$ be the well-quasi-order obtained from $(Q,\preceq)$ and $({\mathbb N} \cup \{-1,0\}, \leq)$ by taking Cartesian product.
For each $i$, let $\phi_i':V(D_i) \rightarrow Q_1$ such that $\phi_i'(v)=(\phi_i(v),-1)$ for every $v \in V(D_i)-\{x_i,y_i\}$, and $\phi_i'(v)=(\phi_i(v),\ell_i)$ for $v \in \{x_i,y_i\}$.
For each $i$ and each $e \in E(W_i)$, define $\mu_i(e)$ to be the number of edges of $D_i'$ with tail and head equal to $e$.

Since each $W_i$ is a directed path, it is a rooted tree rooted at $x_i$.
By Theorem \ref{wqo gap}, there exist $i,j$ with $1 \leq i <j$ and a homeomorphic embedding $\eta$ from $W_i$ to $W_j$ such that $\phi_i'(v) \preceq_1 \phi_j'(\eta(v))$ for every $v \in V(D_i)$, and $\mu_i(e) \leq \mu_j(f)$ for every $e \in E(W_i)$ and $f \in E(\eta(e))$.

Since $\ell_i \geq 0>-1$ and $\phi_i'(x_i) \preceq_1 \phi_j'(\eta(x_i))$, $\eta(x_i) \in \{x_j,y_j\}$.
Similarly, $\eta(y_i) \in \{x_j,y_j\}$.
Since $W_i$ is from $x_i$ to $y_i$, and $W_j$ is from $x_j$ to $y_j$, we know that $\eta(x_i)=x_j$ and $\eta(y_i)=y_j$.
In addition, since $\phi_i'(x_i) \preceq_1 \phi_j'(\eta(x_i))=\phi_j'(x_j)$, $\ell_i \leq \ell_j$.
So there are $\ell_j \geq \ell_i$ directed edges in $D_j$ from $y_j$ to $x_j$.
Moreover, for any directed edge $e=(x_e,y_e)$ in $W_i$, since $\mu_i(e) \leq \mu_j(f)$ for every $f \in E(\eta(e))$, we know that there are at least $\mu_i(e)$ edge-disjoint directed paths in $D_j$ from $\eta(x_e)$ to $\eta(y_e)$ internally disjoint from $\eta(V(D_i))$, so there exists an injection $\eta_e$ from the set of edges of $D_i$ from $x_e$ to $y_e$ to the set of those paths in $D_j$.

Define $\eta^*$ to be a function with domain $V(D_i) \cup E(D_i)$ such that 
	\begin{itemize}
		\item $\eta^*(v)=\eta(v)$ for every $v \in V(D_i)=V(W_i)$, 
		\item $\eta^*$ maps the edges of $D_i$ from $y_i$ to $x_i$ to edges of $D_j$ from $y_j=\eta^*(y_i)$ to $x_j=\eta^*(x_i)$ injectively, and
		\item for each edge $f$ of $D_i$ not from $y_i$ to $x_i$, $\eta^*(f)=\eta_e(f)$, where $e$ is the edge of $W_i$ having the same tail and head as $f$.
	\end{itemize}
Then $\eta^*$ is a strong immersion embedding from $D_i$ to $D_j$ such that $\phi_i(v) \preceq \phi_j(\eta^*(v))$ for every $v \in V(D_i)$.
\end{pf}

\section{Series-parallel triples} \label{sec:sp_triple}

A {\it separation} of a graph (or a directed graph, respectively) $G$ is an ordered pair $(A,B)$ of edge-disjoint subgraphs (or subdigraphs, respectively) such that $A \cup B=G$. 
The {\it order} of $(A,B)$ is $\lvert V(A \cap B) \rvert$.

A {\it series-parallel triple} $(D,s,t)$ is a triple where $D$ is a directed graph whose underlying graph is connected and $s,t$ are distinct vertices of $D$ such that every thread in $D$ from $s$ to $t$ is a directed path, and there exists no separation $(A,B)$ of $D$ of order at most one such that $s,t \in V(A)$ and $V(B)-V(A) \neq \emptyset$. 
A series-parallel triple $(D,s,t)$ is {\it one-way} if either every thread in $D$ is a directed path from $s$ to $t$, or every thread in $D$ is a directed path from $t$ to $s$.

It was shown in \cite[Lemma 5.2]{m} that a one-way series parallel triple can be constructed by a sequence of certain series operations and parallel operations.
So it justifies its name. 
The following simple lemma shows that a series-parallel triple can also be constructed by series and parallel operations even though it is not one-way.
It is likely a folklore result, but we include it in this paper for completeness.

\begin{lemma} \label{constr_sp}
If $(D,s,t)$ is a series-parallel triple, then either
	\begin{enumerate}
		\item $D$ consists of an edge with ends $s$ and $t$, or
		\item there exist series-parallel triples $(D_1,s_1,t_1)$ and $(D_2,s_2,t_2)$ with $\lvert E(D_1) \rvert < \lvert E(D) \rvert$ and $\lvert E(D_2) \rvert < \lvert E(D) \rvert$ such that either
			\begin{enumerate}
				\item $s=s_1$, $t=t_2$, and $D$ is obtained from the disjoint union of $D_1$ and $D_2$ by identifying $t_1$ and $s_2$, or
				\item $D$ is obtained from the disjoint union of $D_1$ and $D_2$ by identifying $s_1$ and $s_2$ into $s$ and identifying $t_1$ and $t_2$ into $t$. 
			\end{enumerate}
	\end{enumerate}
\end{lemma}

\begin{pf}
We may assume that $D$ contains at least two edges for otherwise we are done.
When $\lvert V(D) \rvert=2$, Statement 2(b) holds.
So we may assume that $D$ contains at least three vertices.

We first assume that there exists a separation $(A,B)$ of $D$ of order one such that $s \in V(A)-V(B)$ and $t \in V(B)-V(A)$.
Let $z$ be the vertex in $V(A \cap B)$.
Since $(D,s,t)$ is a series-parallel triple and $z \not \in \{s,t\}$, $D-z$ has exactly two components, where one is $A-z$ and the other is $B-z$, for otherwise there exists a separation $(A',B')$ of $D$ such that $\{s,t\} \subseteq V(A')$ and $B'$ contains a component of $D-z$ disjoint from $s$ and $t$, contradicting that $(D,s,t)$ is a series-parallel triple. 
Since every thread in $A$ between $s$ and $z$ can be made a thread in $D$ between $s$ and $t$ by concatenating a thread between $z$ and $t$, every thread in $A$ between $s$ and $z$ is a directed path.
And there exists no separation $(A',B')$ of $A$ of order at most one such that $\{s,z\} \subseteq V(A')$ and $V(B')-V(A') \neq \emptyset$, for otherwise $(A' \cup B, B')$ is a separation of $D$ of order at most one such that $\{s,t\} \subseteq V(A' \cup B)$ and $V(B')-V(A' \cup B) \neq \emptyset$, contradicting that $(D,s,t)$ is a series-parallel triple.
So $(A,s,z)$ is a series-parallel triple.
Similarly, $(B,z,t)$ is a series-parallel triple.
Hence Statement 2(a) holds.

Therefore we may assume that there exists no separation $(A,B)$ of $D$ of order one such that $s \in V(A)-V(B)$ and $t \in V(B)-V(A)$.
Since $\lvert V(D) \rvert \geq 3$, there exist two internally disjoint threads $P_1,P_2$ in $D$ between $s$ and $t$.
Since $(D,s,t)$ is a series-parallel graph, $P_1,P_2$ are directed paths in $D$, and there exists no thread in $D$ between $V(P_1)-\{s,t\}$ and $V(P_2)-\{s,t\}$.
Hence there exists a separation $(A,B)$ of $D$ such that $V(A \cap B)=\{s,t\}$ and $E(P_1) \subseteq E(A)$ and $E(P_2) \subseteq E(B)$.
Since $(D,s,t)$ is a series-parallel triple, there exists no separation $(A',B')$ of $D$ of order at most one such that $\{s,t\} \subseteq V(A')$ and $V(B')-V(A') \neq \emptyset$, so $(A,s,t)$ and $(B,s,t)$ are series-parallel triples.
So Statement 2(b) holds.
\end{pf}

\bigskip

For a series-parallel triple $(D,s,t)$, we say that 
	\begin{itemize}
		\item $(D,s,t)$ is {\it series-irreducible} if either $\lvert E(D) \rvert=1$, or Statement 2(a) in Lemma \ref{constr_sp} does not hold, and
		\item $(D,s,t)$ is {\it parallel-irreducible} if either $\lvert E(D) \rvert=1$, or Statement 2(b) in Lemma \ref{constr_sp} does not hold.
	\end{itemize}

\section{Series-parallel trees} \label{sec:sp_tree}

A {\it march} is a sequence with distinct entries.
A {\it general rooted digraph} is a pair $(D,\sigma)$, where $D$ is a digraph and $\sigma$ is a march over $V(D)$.
We call $\sigma$ the {\it root march} of a general rooted digraph $(D,\sigma)$.
A {\it rooted digraph} is a pair $(D,v)$, where $D$ is a directed graph and $v \in V(D)$, and we call $v$ the root of $D$.

For simplicity of notations, we do not distinguish the rooted digraph $(D,v)$ and the general rooted digraph $(D,(v))$; and we do not distinguish the series-parallel triple $(D,s,t)$ and the general rooted digraph $(D,(s,t))$.

A {\it strong immersion embedding} from a general rooted digraph $(H,\sigma_H)$ to a general rooted digraph $(G,\sigma_G)$ is a strong immersion embedding $\eta$ from $H$ to $G$ such that $\eta(\sigma_H)=\sigma_G$.
Note that it implies that $\sigma_H$ and $\sigma_G$ have the same length.

Let $(Q,\preceq)$ be a quasi-order.
Let $(D,\sigma)$ and $(D',\sigma')$ be general rooted digraphs.
Let $\phi:V(D) \rightarrow Q$ and $\phi':V(D') \rightarrow Q$ be functions.
We say that $((D',\sigma'),\phi')$ {\it simulates} $((D,\sigma),\phi)$ if there exists a strong immersion embedding $\eta$ from $(D,\sigma)$ to $(D',\sigma')$ such that $\phi(v) \preceq \phi'(\eta(v))$ for every $v \in V(D)$.

A set $\F$ of general rooted digraphs is {\it well-behaved} if for every infinite sequence of general rooted digraphs $(D_1,\sigma_1),(D_2,\sigma_2),... \in \F$, every well-quasi-order $(Q,\leq_Q)$ and functions $\phi_i: V(D_i) \rightarrow Q$ for each $i \geq 1$, there exist $1 \leq j <j'$ such that $((D_{j'},\sigma_{j'}),\phi_{j'})$ simulates $((D_j,\sigma_j),\phi_j)$. 

A {\it cut-vertex} of a graph $G$ is a vertex $v$ of $G$ such that $G-v$ has more components than $G$.
A {\it block} of a graph $G$ is a maximal subgraph $B$ such that $B$ does not contain any cut-vertex of $B$.
A {\it block} of a directed graph $D$ is a directed subgraph whose underlying graph is a block of the underlying graph of $D$.

For a rooted digraph $(D,r)$ in which the underlying graph of $D$ is connected, the {\it block-structure} of $(D,r)$ is a rooted tree $T$ such that the following hold.
	\begin{itemize}
		\item There exists a bipartition $\{L,C\}$ of $V(T)$.
		\item There exists a bijection $f_C$ from $C$ to the set that is the union of $\{r\}$ and the set of cut-vertices of the underlying graph of $D$.
		\item There exists a bijection $f_L$ from $L$ to the set of blocks of the underlying graph of $D$.
		\item For any $v \in C$ and $B \in L$, $v$ is adjacent in $T$ to $B$ if and only if $f_C(v) \in V(f_L(B))$.
		\item The vertex of $T$ mapped to $r$ by $f_C$ is the root of $T$.
	\end{itemize}
For a block $B'$ of the underlying graph of $D$, a {\it child block} of $B'$ is a block $B''$ of the underlying graph of $D$ such that $V(B') \cap V(B'') \neq \emptyset$ and the vertex of $T$ mapped to $B''$ by $f_L$ is a descendant of the vertex of $T$ mapped to $B'$ by $f_L$.
If $B''$ is a child block of $B'$, then we say that $B'$ is the {\it parent block} of $B''$.

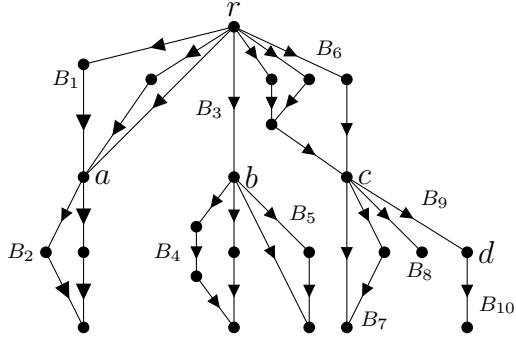
\begin{figure} 
\begin{tikzpicture}
\filldraw[black] (1.8,2) circle (2pt) node[anchor=south] {$r$};
\filldraw[black] (-0.2,1.5) circle (2pt);
\filldraw[black] (0.7,1.3) circle (2pt);
\filldraw[black] (-0.2,0) circle (2pt) node[anchor=west] {$a$};
\draw (-0.2,0) -- (-0.2,1.5) node[very thick,currarrow,pos=0.5,xscale=-1,sloped,scale=1] {};
\draw (-0.2,1.5) -- (1.8,2) node[thick,currarrow,pos=0.5,xscale=-1,sloped,scale=1] {};
\draw (0.7,1.3) -- (1.8,2) node[thick,currarrow,pos=0.5,xscale=-1,sloped,scale=1] {};
\draw (-0.2,0) -- (0.7,1.3) node[thick,currarrow,pos=0.5,xscale=-1,sloped,scale=1] {};
\draw (-0.2,0) -- (1.8,2) node[thick,currarrow,pos=0.5,xscale=-1,sloped,scale=1] {};
	\draw (-0.1,1.3) node[anchor=east] {{\scriptsize $B_1$}};

\filldraw[black] (-0.7,-1) circle (2pt);
\filldraw[black] (-0.2,-1) circle (2pt);
\filldraw[black] (-0.2,-2) circle (2pt);
\draw (-0.2,-1) -- (-0.2,0) node[very thick,currarrow,pos=0.5,xscale=-1,sloped,scale=1] {};
\draw (-0.2,-2) -- (-0.7,-1) node[very thick,currarrow,pos=0.5,xscale=1,sloped,scale=1] {};
\draw (-0.2,-2) -- (-0.2,-1) node[very thick,currarrow,pos=0.5,xscale=-1,sloped,scale=1] {};
\draw (-0.2,0) -- (-0.7,-1) node[currarrow,pos=0.5,xscale=-1,sloped,scale=1] {};
	\draw (-0.7,-1) node[anchor=east] {\scriptsize {$B_2$}};

\filldraw[black] (1.8,0) circle (2pt) node[anchor=west] {$b$};
\draw (1.8,2) -- (1.8,0) node[currarrow,pos=0.5,xscale=1,sloped,scale=1] {};
	\draw (1.8,0.9) node[anchor=east] {{\scriptsize $B_3$}};

\filldraw[black] (1.3,-0.66) circle (2pt);
\filldraw[black] (1.3,-1.32) circle (2pt);
\filldraw[black] (1.8,-1) circle (2pt);
\filldraw[black] (1.8,-2) circle (2pt);
\draw (1.8,0) -- (1.3,-0.66) node[currarrow,pos=0.5,xscale=-1,sloped,scale=1] {};
\draw (1.8,0) -- (1.8,-1) node[currarrow,pos=0.5,xscale=1,sloped,scale=1] {};
\draw (1.8,-1) -- (1.8,-2) node[currarrow,pos=0.5,xscale=1,sloped,scale=1] {};
\draw (1.3,-0.66) -- (1.3,-1.32) node[currarrow,pos=0.5,xscale=1,sloped,scale=1] {};
\draw (1.3,-1.32) -- (1.8,-2) node[currarrow,pos=0.5,xscale=1,sloped,scale=1] {};
	\draw (1.25,-1) node[anchor=east] {{\scriptsize $B_4$}};

\filldraw[black] (2.8,-1) circle (2pt);
\filldraw[black] (2.8,-2) circle (2pt);
\draw (1.8,0) -- (2.8,-1) node[currarrow,pos=0.5,xscale=1,sloped,scale=1] {};
\draw (2.8,-1) -- (2.8,-2) node[currarrow,pos=0.5,xscale=1,sloped,scale=1] {};
\draw (1.8,0) -- (2.8,-2) node[currarrow,pos=0.5,xscale=1,sloped,scale=1] {};
\draw (3.05,-0.5) node[anchor=east] {{\scriptsize $B_5$}};

\filldraw[black] (2.3,1.3) circle (2pt);
\filldraw[black] (2.8,1.3) circle (2pt);
\filldraw[black] (2.3,0.7) circle (2pt);
\filldraw[black] (3.3,1.3) circle (2pt);
\filldraw[black] (3.3,0) circle (2pt) node[anchor=west] {$c$};
\draw (1.8,2) -- (2.3,1.3) node[currarrow,pos=0.5,xscale=1,sloped,scale=1] {};
\draw (1.8,2) -- (2.8,1.3) node[currarrow,pos=0.5,xscale=1,sloped,scale=1] {};
\draw (2.3,1.3) -- (2.3,0.7) node[currarrow,pos=0.5,xscale=1,sloped,scale=1] {};
\draw (2.8,1.3) -- (2.3,0.7) node[currarrow,pos=0.5,xscale=-1,sloped,scale=1] {};
\draw (2.3,0.7) -- (3.3,0) node[currarrow,pos=0.5,xscale=1,sloped,scale=1] {};
\draw (1.8,2) -- (3.3,1.3) node[currarrow,pos=0.5,xscale=1,sloped,scale=1] {};
\draw (3.3,1.3) -- (3.3,0) node[currarrow,pos=0.5,xscale=1,sloped,scale=1] {};
\draw (3.4,1.7) node[anchor=east] {{\scriptsize $B_6$}};

\filldraw[black] (3.3,-2) circle (2pt);
\filldraw[black] (3.8,-1) circle (2pt);
\draw (3.3,0) -- (3.3,-2) node[currarrow,pos=0.5,xscale=1,sloped,scale=1] {};
\draw (3.3,0) -- (3.8,-1) node[currarrow,pos=0.5,xscale=1,sloped,scale=1] {};
\draw (3.8,-1) -- (3.3,-2) node[currarrow,pos=0.5,xscale=-1,sloped,scale=1] {};
\draw (4,-1.9) node[anchor=east] {{\scriptsize$B_7$}};

\filldraw[black] (4.3,-1) circle (2pt);
\draw (3.3,0) -- (4.3,-1) node[currarrow,pos=0.5,xscale=1,sloped,scale=1] {};
\draw (4.6,-1.25) node[anchor=east] {{\scriptsize $B_8$}};

\filldraw[black] (4.9,-1) circle (2pt) node[anchor=west] {$d$};
\filldraw[black] (4.9,-2) circle (2pt);
\draw (3.3,0) -- (4.9,-1) node[currarrow,pos=0.5,xscale=1,sloped,scale=1] {};
\draw (4.9,-1) -- (4.9,-2) node[currarrow,pos=0.5,xscale=1,sloped,scale=1] {};
\draw (4.8,-0.3) node[anchor=east] {{\scriptsize $B_9$}};
\draw (5.7,-1.7) node[anchor=east] {{\scriptsize $B_{10}$}};

\draw (-6,-1.7) node[anchor=east] {};

\end{tikzpicture}
\caption{An ${\mathcal F}$-series-parallel tree. $B_i$ is a block for each $i \in [10]$, and $r,a,b,c,d$ are cut-vertices. $B_1,B_3,B_6,B_9$ are middle blocks.}
\label{fig:sp_tree}
\end{figure}

Let $\F$ be a set of rooted digraphs.
A rooted digraph $(D,r)$ is a {\it $\F$-series-parallel tree} if the underlying graph of $D$ is connected, and for every block $B$ of $D$, the following hold.
	\begin{itemize}
		\item If $B$ is a block of $D$ containing $r$, then $(B,r) \in \F$.
		\item If $B$ is a block of $D$ not containing $r$, then $(B,v) \in \F$, where $v$ is the cut-vertex of the underlying graph of $D$ contained in $B$ and the parent block of $B$.
		\item For every cut-vertex $v$ of the underlying graph of $D$, every thread in $D$ from $r$ to $v$ is a directed path from $r$ to $v$.
		\item If $r \in V(B)$, then $B$ contains at most one cut-vertex of the underlying graph of $D$ that is not $r$; if $r \not \in V(B)$, then $B$ contains at most one cut-vertex of the underlying graph of $D$ that is not contained in the parent block of $B$.
	\end{itemize}
Note that if $(D,r)$ is a $\F$-series-parallel tree, then the vertices of the block-structure of $(D,r)$ corresponding to blocks are of degree at most two.
Observe that for each block $B$ of $D$ in which $B$ has a child block, $(B,x,y)$ is a series-parallel triple, where $x$ is either $r$ or the cut-vertex contained in $B$ and the parent block of $B$, and $y$ is the cut-vertex contained in $B$ and a child block of $B$.
In this case, we call $(B,x,y)$ a {\it middle block} of $(D,r)$.
See Figure \ref{fig:sp_tree} for an example.

A {\it splitter} of a series-parallel triple $(B,x,y)$ is an ordered partition $[X,Y]$ of $V(B)$ such that $x \in X$, $y \in Y$ and the number of edges with one end in $X$ and one end in $Y$ equals the maximum number of edge-disjoint threads in $B$ between $x$ and $y$.
We define $B_X$ to be the digraph obtained from $B$ by identifying $Y$ into a single vertex $y_Y$ and deleting all resulting loops, so $(B_X,x,y_Y)$ is a series-parallel triple.
Similarly, we define $B_Y$ to be the digraph obtained from $B$ by identifying $X$ into a single vertex $x_X$ and deleting all resulting loops, so $(B_Y,x_X,y)$ is a series-parallel triple.
Each of $(B_X,x,y_Y)$ and $(B_Y,x_X,y)$ is called a {\it truncation} of $(B,x,y)$ (with respect to $[X,Y]$).
Note that if $(B,x,y)$ is a one-way series-parallel triple, then every truncation of $(B,x,y)$ is a one-way series-parallel triple.
If there exists a function $\phi$ with domain $V(B)$, then let $\phi_X$ be the function with domain $V(B_X)$ such that $\phi_X(y_Y)=\phi(y)$ and $\phi_X(v)=\phi(v)$ for every $v \in V(B_x)-\{y_Y\}$, and let $\phi_Y$ be the function with domain $V(B_Y)$ such that $\phi_Y(x_X)=\phi(x)$ and $\phi_Y(v)=\phi(v)$ for every $v \in V(B_x)-\{x_X\}$.

\begin{figure} 
\begin{tikzpicture}
\filldraw[black] (12,6) circle (2pt) node[anchor=south] {{\scriptsize $(0,\phi(r))$}};
\filldraw[black] (7.5,0) circle (2pt) node[anchor=east] {{\scriptsize $(1,\phi(a))$}};
\filldraw[black] (12,0) circle (2pt) node[anchor=east] {{\scriptsize $(1,\phi(b))$}};
\filldraw[black] (16,0) circle (2pt) node[anchor=east] {{\scriptsize $(1,\phi(c))$}};
\filldraw[black] (20,-6) circle (2pt) node[anchor=east] {{\scriptsize $(1,\phi(d))$}};

\draw[black] (7.45,4.45) rectangle (7.55,4.55) node[anchor=east] {{\tiny $(3,(({B_1}_{X_{B_1}},r,a_{Y_{B_1}),\phi|_{X_{B_1}}})$}};
\filldraw[black] (7.43,2.93) rectangle (7.57,3.07) node[anchor=east] {{\tiny $(2,((B_1,r,a),\phi|_{B_1}))$}};
\draw[black] (7.45,1.45) rectangle (7.55,1.55) node[anchor=east] {{\tiny $(4,(({B_1}_{Y_{B_1}},r_{Y_{B_1}},a),\phi|_{Y_{B_1}}))$}};

\draw (12,6) -- (7.5,4.55) node[currarrow,pos=0.5,xscale=-1,sloped,scale=1] {};
\draw (9.75,5.25) node[anchor=south] {{\tiny $\omega$}};
\draw (7.5,4.45) -- (7.5,3) node[currarrow,pos=0.5,xscale=1,sloped,scale=1] {};
\draw (7.5,3.725) node[anchor=west] {{\tiny $3$}};
\draw (7.5,3) -- (7.5,1.55) node[currarrow,pos=0.5,xscale=1,sloped,scale=1] {};
\draw (7.5,2.225) node[anchor=west] {{\tiny $3$}};
\draw (7.5,1.45) -- (7.5,0) node[currarrow,pos=0.5,xscale=1,sloped,scale=1] {};
\draw (7.5,0.725) node[anchor=west] {{\tiny $\omega$}};

\draw[black] (11.95,4.45) rectangle (12.05,4.55) node[anchor=east] {{\tiny $(3,(({B_3}_{X_{B_3}},r,a_{Y_{B_1}),\phi|_{X_{B_3}}})$}};
\filldraw[black] (11.93,2.93) rectangle (12.07,3.07) node[anchor=east] {{\tiny $(2,((B_3,r,a),\phi|_{B_3}))$}};
\draw[black] (11.95,1.45) rectangle (12.05,1.55) node[anchor=east] {{\tiny $(4,(({B_3}_{Y_{B_3}},r_{Y_{B_3}},a),\phi|_{Y_{B_3}}))$}};

\draw (12,6) -- (12,4.55) node[currarrow,pos=0.5,xscale=1,sloped,scale=1] {};
\draw (12,5.25) node[anchor=west] {{\tiny $\omega$}};
\draw (12,4.45) -- (12,3) node[currarrow,pos=0.5,xscale=1,sloped,scale=1] {};
\draw (12,3.725) node[anchor=west] {{\tiny $1$}};
\draw (12,3) -- (12,1.55) node[currarrow,pos=0.5,xscale=1,sloped,scale=1] {};
\draw (12,2.225) node[anchor=west] {{\tiny $1$}};
\draw (12,1.45) -- (12,0) node[currarrow,pos=0.5,xscale=1,sloped,scale=1] {};
\draw (12,0.725) node[anchor=west] {{\tiny $\omega$}};

\draw[black] (15.95,4.45) rectangle (16.05,4.55) node[anchor=west] {{\tiny $(3,(({B_6}_{X_{B_6}},r,a_{Y_{B_6}),\phi|_{X_{B_6}}})$}};
\filldraw[black] (15.93,2.93) rectangle (16.07,3.07) node[anchor=west] {{\tiny $(2,((B_6,r,a),\phi|_{B_6}))$}};
\draw[black] (15.95,1.45) rectangle (16.05,1.55) node[anchor=west] {{\tiny $(4,(({B_6}_{Y_{B_6}},r_{Y_{B_6}},a),\phi|_{Y_{B_6}}))$}};

\draw (12,6) -- (16,4.55) node[currarrow,pos=0.5,xscale=1,sloped,scale=1] {};
\draw (14,5.23) node[anchor=south] {{\tiny $\omega$}};
\draw (16,4.45) -- (16,3) node[currarrow,pos=0.5,xscale=1,sloped,scale=1] {};
\draw (16,3.725) node[anchor=east] {{\tiny $2$}};
\draw (16,3) -- (16,1.55) node[currarrow,pos=0.5,xscale=1,sloped,scale=1] {};
\draw (16,2.225) node[anchor=east] {{\tiny $2$}};
\draw (16,1.45) -- (16,0) node[currarrow,pos=0.5,xscale=1,sloped,scale=1] {};
\draw (16,0.725) node[anchor=east] {{\tiny $\omega$}};

\draw[black] (19.95,-1.55) rectangle (20.05,-1.45) node[anchor=east] {}; 
\draw (19.95,-1.6) node[anchor=east] {{\tiny $(3,(({B_9}_{X_{B_9}},r,a_{Y_{B_9}),\phi|_{X_{B_9}}})$}};
\filldraw[black] (19.93,-3.07) rectangle (20.07,-2.93) node[anchor=east] {{\tiny $(2,((B_9,r,a),\phi|_{B_9}))$}};
\draw[black] (19.95,-4.55) rectangle (20.05,-4.45) node[anchor=east] {{\tiny $(4,(({B_9}_{Y_{B_9}},r_{Y_{B_9}},a),\phi|_{Y_{B_9}}))$}};

\draw (16,0) -- (20,-1.45) node[currarrow,pos=0.5,xscale=1,sloped,scale=1] {};
\draw (18,-0.73) node[anchor=south] {{\tiny $\omega$}};
\draw (20,-1.55) -- (20,-3) node[currarrow,pos=0.5,xscale=1,sloped,scale=1] {};
\draw (20,-2.275) node[anchor=west] {{\tiny $1$}};
\draw (20,-3) -- (20,-4.45) node[currarrow,pos=0.5,xscale=1,sloped,scale=1] {};
\draw (20,-3.775) node[anchor=west] {{\tiny $1$}};
\draw (20,-4.55) -- (20,-6) node[currarrow,pos=0.5,xscale=1,sloped,scale=1] {};
\draw (20,-5.275) node[anchor=west] {{\tiny $\omega$}};

\filldraw[black] (7.43,-1.43) rectangle (7.57,-1.57) node[anchor=east] {{\tiny $(5,((B_2,a),\phi|_{B_2}))$}};
\draw (7.5,0) -- (7.5,-1.5) node[currarrow,pos=0.5,xscale=1,sloped,scale=1] {};
\draw (7.5,-0.75) node[anchor=west] {{\tiny $\omega$}};

\filldraw[black] (8.93,-3.07) rectangle (9.07,-2.93) node[anchor=east] {{\tiny $(5,((B_4,b),\phi|_{B_4}))$}};
\filldraw[black] (11.93,-3.07) rectangle (12.07,-2.93) node[anchor=east] {{\tiny $(5,((B_5,b),\phi|_{B_5}))$}};
\draw (12,0) -- (9,-3) node[currarrow,pos=0.5,xscale=-1,sloped,scale=1] {};
\draw (10.5,-1.5) node[anchor=west] {{\tiny $\omega$}};
\draw (12,0) -- (12,-3) node[currarrow,pos=0.5,xscale=1,sloped,scale=1] {};
\draw (12,-1.5) node[anchor=west] {{\tiny $\omega$}};

\filldraw[black] (12.93,-6.07) rectangle (13.07,-5.93) node[anchor=east] {{\tiny $(5,((B_7,c),\phi|_{B_7}))$}};
\filldraw[black] (15.93,-6.07) rectangle (16.07,-5.93) node[anchor=east] {{\tiny $(5,((B_8,c),\phi|_{B_8}))$}};
\draw (16,0) -- (13,-6) node[currarrow,pos=0.5,xscale=-1,sloped,scale=1] {};
\draw (14.5,-3) node[anchor=east] {{\tiny $\omega$}};
\draw (16,0) -- (16,-6) node[currarrow,pos=0.5,xscale=1,sloped,scale=1] {};
\draw (16,-3) node[anchor=east] {{\tiny $\omega$}};

\filldraw[black] (19.93,-7.57) rectangle (20.07,-7.43) node[anchor=east] {{\tiny $(5,((B_{10},d),\phi|_{B_{10}}))$}};
\draw (20,-6) -- (20,-7.5) node[currarrow,pos=0.5,xscale=1,sloped,scale=1] {};
\draw (20,-6.75) node[anchor=west] {{\tiny $\omega$}};

\end{tikzpicture}
	\caption{The ${\mathcal S}$-portrait of the ${\mathcal F}$-series-parallel tree in Figure \ref{fig:sp_tree}, assuming $\Se$ is given. Solid circles and rectangles are vertices in the block-structure corresponding to cut-vertices and blocks, respectively. Empty rectangles are the vertices obtained by subdividing edges.}
	\label{fig:portrait}
\end{figure}
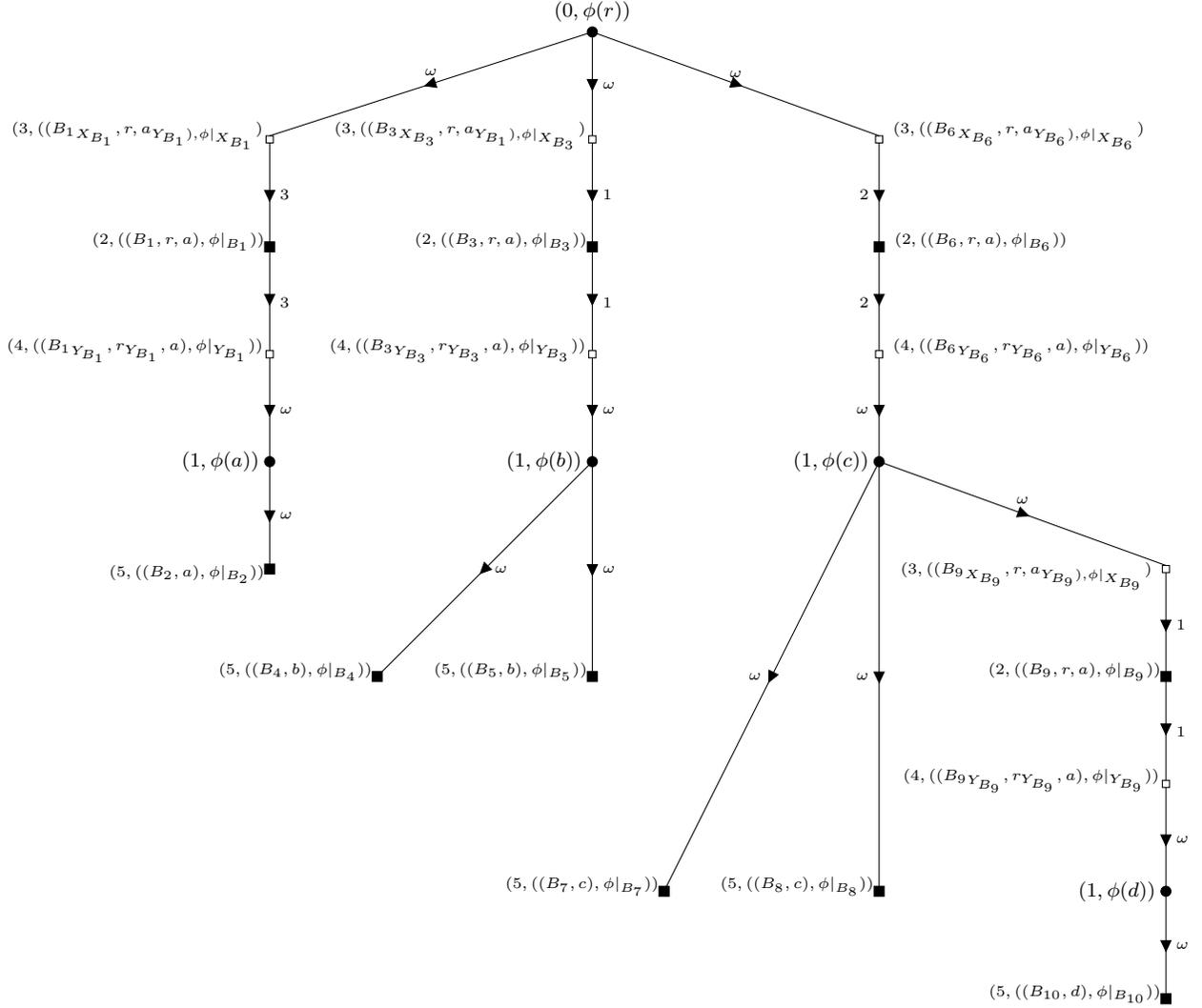

Let $\F$ be a family of rooted digraphs.
Let $(D,r)$ be a $\F$-series-parallel tree, and let $\phi$ be a function with domain $V(D)$.
For each middle block $(B,x,y)$ of $(D,r)$, we choose a splitter $[X_B,Y_B]$ of $(B,x,y)$.
Let $\Se$ be the set of $[X_B,Y_B]$ over all middle blocks $(B,x,y)$ of $(D,r)$.
The {\it $\Se$-portrait} of $((D,r),\phi)$ is a pair $(T,\psi)$ such that the following hold.
	\begin{itemize}
		\item $T$ is a tree, and $\psi$ is a function with domain $V(T) \cup E(T)$.
		\item $T$ is obtained from the block-structure of $(D,r)$ by subdividing each edge that is not incident with an non-root leaf once.
		\item $\psi(r)=(0,\phi(r))$.
		\item $\psi$ maps each node $t$ of $T$ corresponding to a cut-vertex of the underlying graph of $D$ to $(1,\phi(t))$.
		\item $\psi$ maps each node $t$ of $T$ corresponding to a middle block $(B,x,y)$ of $(D,r)$ to $(2,((B,x,y),\phi|_{V(B)}))$. 
		\item $\psi$ maps each node $t$ of $T$ that is obtained by subdividing an edge whose head corresponds to a middle block $(B,x,y)$ with splitter $(X_B,Y_B) \in \Se$ to $(3,((B_{X_B},x,y_{Y_B}),\phi|_{X_B}))$.
		\item $\psi$ maps each node $t$ of $T$ that is obtained by subdividing an edge whose tail corresponds to a middle block $(B,x,y)$ with splitter $(X_B,Y_B) \in \Se$ to $(4,((B_{Y_B},x_{X_B},y),\phi|_{Y_B}))$.
		\item $\psi$ maps each node $t$ of $T$ that corresponds to a block $B$ of $D$ with no child block to $(5,((B,x),\phi|_{V(B)}))$, where either $x=r$ or $x$ is the cut-vertex contained in $B$ and the parent block of $B$.
		\item $\psi$ maps each edge of $T$ incident with a node corresponding to a cut-vertex or $r$ to $\omega$.
		\item $\psi$ maps each edge of $T$ incident with a node corresponding to a middle block $(B,x,y)$ to the number of edges with one end in $X_B$ and one end in $Y_B$.
	\end{itemize}
See Figure \ref{fig:portrait} for an example.

\begin{lemma} \label{series parallel tree wqo}
Let $\F$ be a well-behaved family of rooted digraphs.
Let $\F'$ be the set of one-way series-parallel triples $(D,s,t)$ such that $(D,s) \in \F$. 
Let $\F''$ be the set consisting of all series-parallel triples that are truncations of members of $\F'$.
If $\F'$ and $\F''$ are well-behaved, then the set of $\F$-series-parallel trees is well-behaved.
\end{lemma}

\begin{pf}
Let $(Q,\leq_Q)$ be a well-quasi-order.
For $i \geq 1$, let $(D_i,r_i)$ be a $\F$-series-parallel tree, and let $\phi_i: V(D_i) \rightarrow Q$.
For each $i \geq 1$ and middle block $(B,x,y)$ of $(D_i,r_i)$, let $[X_B,Y_B]$ be a splitter of $(B,x,y)$.
For each $i \geq 1$, let $\Se_i=\{[X_B,Y_B]: (B,x,y)$ is a middle block of $(D_i,r_i)\}$, and let $(T_i,\psi_i)$ be the $\Se_i$-portrait of $((D_i,r_i),\phi_i)$.

Let $Q_1$ be the set consisting of the pairs $((D,s,t),\phi)$ of a series-parallel triple and a function such that there exists $i$ such that either $(D,s,t)$ is a middle block of $(D_i,r_i)$ and $\phi=\phi_i|_{V(D)}$, or $(D,s,t)$ is a truncation of a middle block $(B,x,y)$ of $(D_i,r_i)$ with respect to $[X_B,Y_B]$ and $\phi: V(D) \rightarrow Q$ is the function obtained from $\phi_i$ defined in the truncation. 
So the simulation relation, denoted by $\preceq_1$, is a quasi-order defined on $Q_1$.
Note that for every $((D,s,t),\phi) \in Q_1$, $(D,s,t) \in \F' \cup \F''$.
Since $\F'$ and $\F''$ are well-behaved, $(Q_1,\preceq_1)$ is a well-quasi-order.
Let $(Q_2,\preceq_2)$ be the well-quasi-order obtained by the disjoint union of $(Q,\leq_Q)$ and $(Q_1,\preceq_1)$.

Let $Q_3$ be the set consisting of the pairs $((D,r),\phi)$ such that there exists $i \in {\mathbb N}$ such that $D$ is a block of $D_i$ with no child block, $r$ is $r_i$ (if $r_i \in V(D)$) or the cut-vertex of the underlying graph of $D_i$ contained in $D$ (if $r_i \not \in V(D)$), and $\phi=\phi_i|_{V(B)}$.
Let $\preceq_3$ be the simulation relation defined on $Q_3$.
Since $\F$ is well-behaved, $(Q_3,\preceq_3)$ is a well-quasi-order.
Let $(Q_4,\preceq_4)$ be the well-quasi-order obtained from the disjoint union of $(Q_2,\preceq_2)$ and $(Q_3,\preceq_3)$. 
Define $(Q',\preceq)$ to be the well-quasi-order obtained by the Cartesian product of $(\{0,1,2,3,4,5\},=)$ and $(Q_4,\preceq_4)$.

Note that the image of each $\psi_i|_{V(T_i)}$ is contained in $Q'$.
By Theorem \ref{wqo gap}, there exist $1 \leq j <j'$ and a homeomorphic embedding $\eta$ from $T_j$ to $T_{j'}$ such that $\psi_j(v) \preceq \psi_{j'}(\eta(v))$ for every $v \in V(T_j)$, and $\psi_j(e) \leq \psi_{j'}(e')$ for every $e \in E(T_j)$ and $e' \in E(\eta(e))$.

Note that by the definition of $\psi_j$ and $\psi_{j'}$, for each middle block $(B,x,y)$ of $(D_j,r_j)$, there exists a node $t$ of $T_j$ such that $t$ corresponds to $(B,x,y)$, and $\eta(t)$ corresponds to a middle block of $(D_{j'},r_{j'})$.
For simplicity, for each middle block $(B,x,y)$ of $(D_j,r_j)$, we write $\eta(B,x,y)$ to denote the middle block of $(D_{j'},r_{j'})$ corresponding to $\eta(t)$, where $t$ is the node of $T_j$ corresponding to $(B,x,y)$, and write $\eta(B)$ to denote the first entry of $\eta(B,x,y)$.
Similarly, we write $\eta(B_{X_B},x,y_{Y_B})$ and $\eta(B_{Y_B},x_{X_B},y)$ to denote those series-parallel triples corresponding to $\eta(t)$, where $t$ is the node of $T_j$ with $\psi_j(t) = (3,((B_{X_B},x,y_{Y_B}),\phi_j|_{X_B}))$ and $\psi_j(t)=(4,((B_{Y_B},x_{X_B},y),\phi_j|_{Y_B}))$, respectively, and we write $\eta(B_{X_B})$ and $\eta(B_{Y_B})$ to denote the middle blocks of $(D_{j'},r_{j'})$ such that the first entries of $\eta(B_{X_B},x,y_{Y_B})$ and $\eta(B_{Y_B},x_{X_B},y)$, respectively, are obtained from $\eta(B_{X_B})$ and $\eta(B_{Y_B})$ by identifying vertices, respectively.
And for each middle block $(B,x,y)$ of $(D_j,r_j)$, we denote the corresponding strong immersion embedding that witness $\psi_j(t) \preceq \psi_{j'}(\eta(t))$ as $\eta_B, \eta_{X_B},\eta_{Y_B}$, respectively, where $t$ is the node of $T_j$ corresponding to $B, B_{X_B}, B_{Y_B}$, respectively.

We say that a middle block $(B,x,y)$ of $(D_j,r_j)$ is {\it tight} if the three nodes of $T_j$ corresponding to $(B_{X_B},x,y_{Y_B})$, $(B,x,y)$ and $(B_{Y_B},x_{X_B},y)$ are mapped by $\eta$ to a path in $T_{j'}$ on three vertices; otherwise we say $(B,x,y)$ is {\it loose}.
Note that by the definition of $\psi_j$ and $\psi_{j'}$, for every middle block $(B,x,y)$ of $(D_j,r_j)$, $(B,x,y)$ is tight if and only if $\eta(B_{X_B})=\eta(B_{Y_B})$.

\medskip

\noindent{\bf Claim 1:} Let $(B,x,y)$ be a loose middle block of $(D_j,r_j)$.
Let $s$ be the maximum number of edge-disjoint threads in $B$ from $x$ to $y$. 
Let $S_X$ be an $s$-element subset of the set of edges of $\eta(B_{X_B})$ between $X_{\eta(B_{X_B})}$ and $Y_{\eta(B_{X_B})}$. 
Let $S_Y$ be an $s$-element subset of the set of edges of $\eta(B_{Y_B})$ between $X_{\eta(B_{Y_B})}$ and $Y_{\eta(B_{Y_B})}$. 
Let $f$ be a bijection between $S_X$ and $S_Y$.
Then there exist $s$ edge-disjoint directed paths in $D_{j'}$ between $X_{\eta(B_{X_B})}$ and $Y_{\eta(B_{Y_B})}$ internally disjoint from $X_{\eta(B_{X_B})} \cup Y_{\eta(B_{Y_B})}$ such that each path contains $e$ and $f(e)$ for some $e \in S_X$.

\medskip

\noindent{\bf Proof of Claim 1:}
Let $W_1,W_2,...,W_k$ (for some integer $k \geq 2$) be the blocks of $D_{j'}$ such that $W_1=\eta(B_{X_B})$, $W_k=\eta(B_{Y_B})$ and every thread in $D_{j'}$ from $V(W_1)$ to $V(W_k)$ intersects $W_i$ for every $1 \leq i \leq k$.
For each $i \in [k]$, let $x_i$ and $y_i$ be the distinct vertices such that each of them is either equal to $r_{j'}$ or a cut-vertex of the underlying graph of $D_{j'}$ contained in $W_i$; and we assume that $x_i$ is closer to $r_{j'}$ than $y_i$.
Let $u_1,u_2,...u_s$ be the ends of the edges in $S_X$ contained in $X_{\eta(B_{X_B})}$.
Let $v_1,v_2,...v_s$ be the ends of the edges in $S_Y$ contained in $Y_{\eta(B_{Y_B})}$.
Note that $u_1,u_2,...,u_s$ are not necessarily distinct, and $v_1,v_2,...,v_s$ are not necessarily distinct.

Note that the two edges of $T_j$ incident with the node of $T_j$ corresponding to $(B,x,y)$ are mapped to $s$ by $\psi_j$.
So every edge of $T_{j'}$ incident with a node of $T_{j'}$ corresponding to one of $W_1,W_2,...,W_k$ is mapped to a number at least $s$ by $\psi_{j'}$.
Hence there exist $s$ edge-disjoint directed paths $M_1,M_2,...,M_s$ in $D_{j'}$ between $y_1$ to $x_k$ internally disjoint from $y_1$ and $x_k$.

Since $(W_1,x_1,y_1)$ is a one-way series-parallel triple, and $[X_{\eta(B_{X_B})},Y_{\eta(B_{X_B})}]$ is a splitter of $(W_1,x_1,y_1)$, there exist edge-disjoint directed paths in $W_1$ from $x_1$ to $y_1$ such that each path intersects exactly one edge between $X_{\eta(B_{X_B})}$ and $Y_{\eta(B_{Y_B})}$.
So $s$ of them intersects $S_X$.
Hence the subpaths $U_1,U_2,...,U_s$ of those $s$ paths are $s$ edge-disjoint directed paths in $W_1[\{u_1,u_2,...,u_s\} \cup Y_{\eta(B_{X_B}})]$ between $\{u_1,u_2,...,u_s\}$ and $y_1$.
Similarly, there exist $s$ edge-disjoint directed paths $U'_1,U'_2,...,U_s'$ in $W_k[X_{\eta(B_{Y_B})} \cup \{v_1,v_2,...,v_s\}]$ between $x_k$ and $\{v_1,v_2,...,v_s\}$.

By symmetry, we may denote the elements of $S_X$ by $e_1,e_2,...,e_s$ and the elements of $S_Y$ by $e_1',e_2',...,e_s'$ such that for every $i \in [s]$, $f(e_i)=e_i'$, $U_i$ contains $e_i$ and $U_i'$ contains $e_i'$.
Hence $U_1 \cup M_1 \cup U_1', U_2 \cup M_2 \cup U_2', ..., U_s \cup M_s \cup U_s'$ are desired directed paths in $D_{j'}$.
$\Box$

\medskip

\noindent{\bf Claim 2:} Let $(B,x,y)$ be a loose middle block of $(D_j,r_j)$.
Let $S=\{e_1,e_2,...,e_{\lvert S \rvert}\}$ be the set of edges of $B$ between $X_B$ and $Y_B$. 
For each $i \in [\lvert S \rvert]$, let $u_i$ be the end of $e_i$ in $X_B$, and let $v_i$ be the end of $e_i$ in $Y_B$.
Let $W_1,W_2,...,W_k$ (for some integer $k \geq 2$) be the blocks of $(D_{j'},r_{j'})$ such that $W_1=\eta(B_{X_B})$, $W_k=\eta(B_{Y_B})$, and every thread in $D_{j'}$ from $V(W_1)$ to $V(W_k)$ intersects $W_i$ for every $1 \leq i \leq k$.
Then there exist $\lvert S \rvert$ edge-disjoint directed paths $P_{e_1},P_{e_2},...,P_{e_{\lvert S \rvert}}$ in $\bigcup_{i=1}^kW_i$ internally disjoint from the image of $\eta_{X_B}|_{X_B}$ and $\eta_{Y_B}|_{Y_B}$ such that for each $i \in [\lvert S \rvert]$, $P_{e_i}$ is between $\eta_{X_B}(u_i)$ and $\eta_{Y_B}(v_i)$ containing $\eta_{X_B}(e_i) \cup \eta_{Y_B}(e_i)$.

\medskip

\noindent{\bf Proof of Claim 2:}
Note that every edge in $S$ is an edge of $B_{X_B}$.
So there are $\lvert S \rvert$ edges of $\eta(B_{X_B})$ incident with the third entry of $\eta(B_{X_B},x,y_{Y_B})$ contained in $\bigcup_{i=1}^{\lvert S \rvert}\eta_{X_B}(e_i)$, and those edges are between $X_{\eta(B_{X_B})}$ and $Y_{\eta(B_{X_B})}$.
Let $S_X$ be the set consisting of those $\lvert S \rvert$ edges.
Similarly, there exists a set $S_Y$ consisting of $\lvert S \rvert$ edges of $\eta(B_{Y_B})$ between $X_{\eta(B_{Y_B})}$ and $Y_{\eta(B_{Y_B})}$ contained in $\bigcup_{i=1}^{\lvert S \rvert}\eta_{Y_B}(e_i)$.
By Claim 1, there exist edge-disjoint directed paths $Z_1,Z_2,...,Z_{\lvert S \rvert}$ in $D_{j'}$ between $X_{\eta(B_{X_B})}$ and $Y_{\eta(B_{Y_B})}$ internally disjoint from $X_{\eta(B_{X_B})} \cup Y_{\eta(B_{Y_B})}$ such that for every $i \in [\lvert S \rvert]$, $Z_i$ intersects $\eta_{X_B}(e_i)$ and $\eta_{Y_B}(e_i)$.
For each $i \in [\lvert S \rvert]$, define $P_{e_i} = \eta_{X_B}(e_i)[X_{\eta(B_{X_B})}] \cup Z_i \cup \eta_{Y_B}(e_i)[Y_{\eta(B_{Y_B})}]$.
Then $P_{e_1},P_{e_2},...,P_{e_{\lvert S \rvert}}$ are desired directed paths.
$\Box$

\medskip

For a block $B$ of $D_j$ with no child block, let $\eta_B$ be the strong immersion embedding from $B$ to $B'$ witnessing $\psi_j(t) \preceq \psi_{j'}(\eta(t))$, where $t$ is the node of $T_j$ corresponding to $B$, and $B'$ corresponds to $\eta(t)$.

Define $\eta^*$ to be a function with domain $V(D_j) \cup E(D_j)$ such that the following statements hold.
	\begin{itemize}
		\item If $v$ is a cut-vertex of $D_j$ or $v=r_j$, then define $\eta^*(v)=\eta(v)$.
		\item If $v$ is a vertex belonging to a block $B$ of $D_j$ with no child block or belonging a tight middle block $(B,x,y)$ of $(D_j,r_j)$, and $v$ is not a cut-vertex of $D_j$ or $r_j$, then define $\eta^*(v)=\eta_B(v)$.
		\item If $v$ is a vertex belonging to a loose middle block $(B,x,y)$ of $(D_j,r_j)$, and $v$ is not a cut-vertex of $D_j$ or $r_j$, then $\eta^*(v)=\eta_{X_B}(v)$ when $v \in X_{B}$, and $\eta^*(v)=\eta_{Y_B}(v)$ when $v \in Y_{B}$.
		\item If $e$ is an edge belonging to a block $B$ of $D_j$ with no child block or belonging to a tight middle block $(B,x,y)$, then define $\eta^*(e)=\eta_B(e)$.
		\item If $e$ is an edge belonging to a loose middle block $(B,x,y)$ with both ends contained in $X_B$ (and $Y_B$, respectively), then $\eta^*(e)=\eta_{B_{X_B}}(e)$ (and $\eta^*(e)=\eta_{B_{Y_B}}(e)$, respectively).
		\item If $e$ is an edge belonging to a loose middle block $(B,x,y)$ between $X_B$ and $Y_B$, then $\eta^*(e)=P_e$, where $P_e$ is the directed path mentioned in Claim 2.
	\end{itemize}
See Figure \ref{fig:eta_loose_block} for an example.

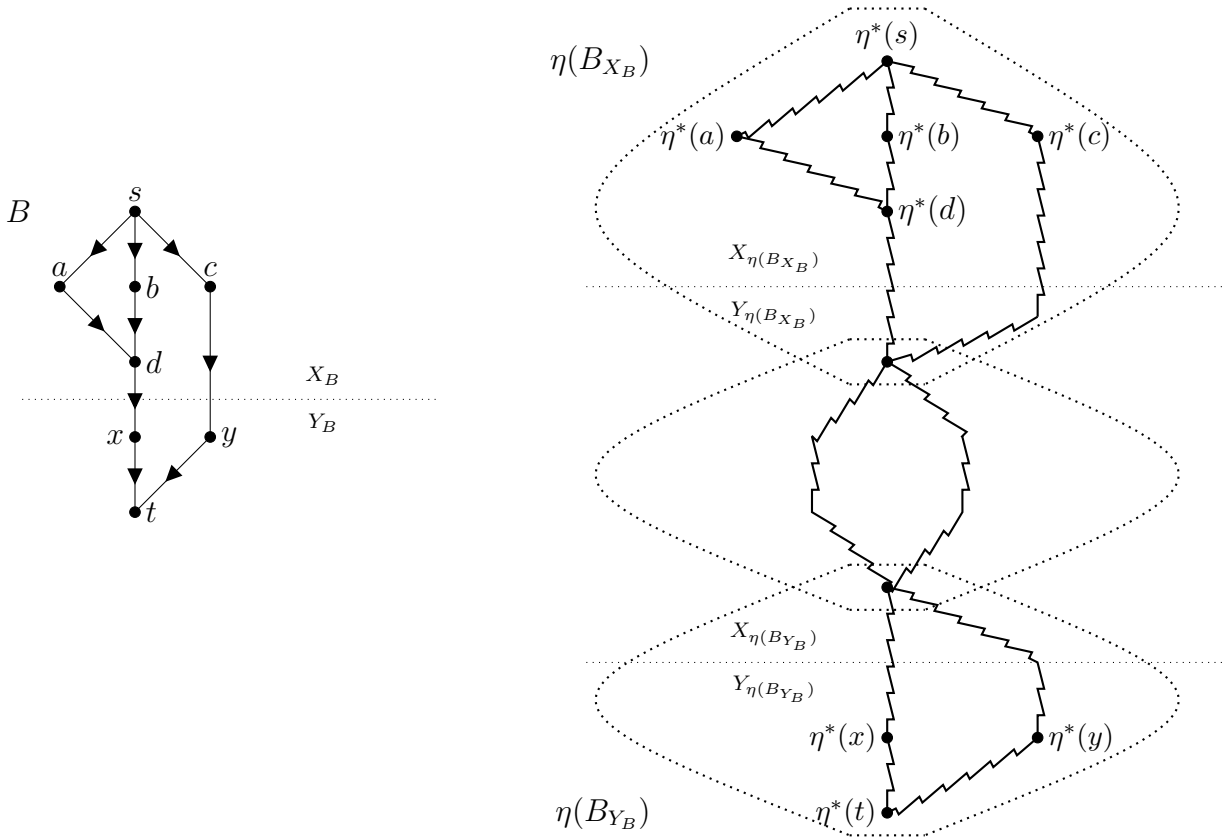
\begin{figure}
\begin{tikzpicture}
\draw (-2,4) node[anchor=west] {{ $B$}};

\filldraw[black] (0,4) circle (2pt) node[anchor=south] {$s$};
\filldraw[black] (-1,3) circle (2pt) node[anchor=south] {$a$};
\filldraw[black] (0,3) circle (2pt) node[anchor=west] {$b$};
\filldraw[black] (1,3) circle (2pt) node[anchor=south] {$c$};
\filldraw[black] (0,2) circle (2pt) node[anchor=west] {$d$};

\filldraw[black] (0,1) circle (2pt) node[anchor=east] {$x$};
\filldraw[black] (1,1) circle (2pt) node[anchor=west] {$y$};
\filldraw[black] (0,0) circle (2pt) node[anchor=west] {$t$};

\draw (0,4) -- (-1,3) node[very thick,currarrow,pos=0.5,xscale=-1,sloped,scale=1] {};
\draw (0,4) -- (0,3) node[very thick,currarrow,pos=0.5,xscale=1,sloped,scale=1] {};
\draw (0,4) -- (1,3) node[very thick,currarrow,pos=0.5,xscale=1,sloped,scale=1] {};
\draw (-1,3) -- (0,2) node[very thick,currarrow,pos=0.5,xscale=1,sloped,scale=1] {};
\draw (0,3) -- (0,2) node[very thick,currarrow,pos=0.5,xscale=1,sloped,scale=1] {};

\draw (0,2) -- (0,1) node[very thick,currarrow,pos=0.5,xscale=1,sloped,scale=1] {};
\draw (1,3) -- (1,1) node[very thick,currarrow,pos=0.5,xscale=1,sloped,scale=1] {};
\draw (0,1) -- (0,0) node[very thick,currarrow,pos=0.5,xscale=1,sloped,scale=1] {};
\draw (1,1) -- (0,0) node[very thick,currarrow,pos=0.5,xscale=-1,sloped,scale=1] {};

\draw[dotted] (-1.5,1.5) -- (4, 1.5);
	\draw (2.5,1.55) node[anchor=south] {{\scriptsize $X_B$}};
	\draw (2.5,1.45) node[anchor=north] {{\scriptsize $Y_B$}};

\draw (7,6) node[anchor=east] {{ $\eta(B_{X_B})$}};

	\filldraw[black] (10,6) circle (2pt) node[anchor=south] {{\small $\eta^*(s)$}};
	\filldraw[black] (8,5) circle (2pt) node[anchor=east] {{\small $\eta^*(a)$}};
	\filldraw[black] (10,5) circle (2pt) node[anchor=west] {{\small $\eta^*(b)$}};
	\filldraw[black] (12,5) circle (2pt) node[anchor=west] {{\small $\eta^*(c)$}};
	\filldraw[black] (10,4) circle (2pt) node[anchor=west] {$\eta^*(d)$};

\draw[thick, snake=saw] (10,6) -- (8,5);
\draw[thick, snake=saw] (10,6) -- (10,5);
\draw[thick, snake=saw] (10,6) -- (12,5);
\draw[thick, snake=saw] (8,5) -- (10,4);
\draw[thick, snake=saw] (10,5) -- (10,4);

\draw[dotted] (6,3) -- (14.5, 3);
	\draw (8.5,3.05) node[anchor=south] {{\scriptsize $X_{\eta(B_{X_B})}$}};
	\draw (8.5,2.95) node[anchor=north] {{\scriptsize $Y_{\eta(B_{X_B})}$}};

	\draw[thick,dotted] (9.5,6.7) .. controls (5,4) .. (9.5, 1.7);
	\draw[thick,dotted] (9.5,6.7) -- (10.5, 6.7);
	\draw[thick,dotted] (10.5,6.7) .. controls (15,4) .. (10.5, 1.7);
	\draw[thick,dotted] (9.5,1.7) -- (10.5, 1.7);

	\filldraw[black] (10,2) circle (2pt) node[anchor=south] {};
\draw[thick, snake=saw] (10,4) -- (10,2);
\draw[thick, snake=saw] (12,5) -- (12,2.6);
\draw[thick, snake=saw] (12,2.6) -- (10,2);

\filldraw[black] (10,-1) circle (2pt) node[anchor=south] {};
	\draw[thick, snake=saw] (10,2) -- (9,1) -- (9, 0) -- (10,-1);
	\draw[thick, snake=saw] (10,2) -- (11,1) -- (11,0) -- (10,-1);

	\draw[thick,dotted] (9.5,2.3) .. controls (5,0.5) .. (9.5, -1.3);
	\draw[thick,dotted] (9.5,2.3) -- (10.5, 2.3);
	\draw[thick,dotted] (10.5,2.3) .. controls (15,0.5) .. (10.5, -1.3);
	\draw[thick,dotted] (9.5,-1.3) -- (10.5,-1.3);

	\draw[dotted] (6,-2) -- (14.5, -2);
	\draw (8.5,-1.98) node[anchor=south] {{\scriptsize $X_{\eta(B_{Y_B})}$}};
	\draw (8.5,-2.03) node[anchor=north] {{\scriptsize $Y_{\eta(B_{Y_B})}$}};

\filldraw[black] (10,-3) circle (2pt) node[anchor=east] {{\small $\eta^*(x)$}};
\filldraw[black] (12,-3) circle (2pt) node[anchor=west] {{\small $\eta^*(y)$}};
\filldraw[black] (10,-4) circle (2pt) node[anchor=east] {{\small $\eta^*(t)$}};

	\draw[thick, snake=saw] (10,-1) -- (10,-3);
	\draw[thick, snake=saw] (10,-1) -- (12,-2) -- (12,-3);
	\draw[thick, snake=saw] (10,-3) -- (10,-4);
	\draw[thick, snake=saw] (12,-3) -- (10,-4);

	\draw (7,-4) node[anchor=east] {{ $\eta(B_{Y_B})$}};
	\draw[thick,dotted] (9.5,-0.7) .. controls (5,-2.5) .. (9.5, -4.3);
	\draw[thick,dotted] (9.5,-0.7) -- (10.5, -0.7);
	\draw[thick,dotted] (10.5,-0.7) .. controls (15,-2.5) .. (10.5, -4.3);
	\draw[thick,dotted] (9.5,-4.3) -- (10.5,-4.3);
\end{tikzpicture}
\caption{An example about how $\eta$ maps vertices and edges in a loose middle block $B$ to vertices and directed paths, respectively. Each zigzag line indicates a directed path with the obvious direction.}
\label{fig:eta_loose_block}
\end{figure}

Clearly, $\eta^*|_{V(D_j)}$ is injective and for every $v \in V(D_j)$, $\phi_j(v) \leq_Q \phi_{j'}(\eta^*(v))$.
Since the edges of $T_j$ and $T_{j'}$ mapped to $\omega$ by $\psi_j$ and $\psi_{j'}$ are exactly the edges incident with $r_j$, $r_{j'}$ or some nodes of $T_j$ and $T_{j'}$ corresponding to cut-vertices of $D_j$ and $D_{j'}$, it is straightforward to verify that $\eta^*$ is a strong immersion embedding from $(D_j,r_j)$ to $(D_{j'},r_{j'})$.
This proves the lemma.
\end{pf}

\section{One-way series-parallel triples} \label{sec:wqo_sp_triples}

The goal of this section is to prove Lemma \ref{series-parallel triple wqo} which shows that series-parallel triples are well-quasi-ordered.
The strategy is to decompose a given series-parallel triple in a ``series way'' or ``parallel way'' as stated in Lemma \ref{constr_sp} to reduce the ``complexity'', and prove well-quasi-ordering by induction on the ``complexity''.
It could be helpful if the readers first read the related definitions and the statements of the lemmas in this section without going into the proofs to get a big picture of the entire procedure.

For a set $\F$ of one-way series-parallel triples, the {\it parallel-extension} of $\F$ is a set $\F'$ of one-way series-parallel triples such that for every $(D,s,t) \in \F'$, there exist a positive integer $\ell$ and members $(D_1,s_1,t_1), (D_2,s_2,t_2), ..., (D_\ell,s_\ell,t_\ell)$ of $\F$ such that $D$ is obtained from the disjoint union of $D_1,D_2,...,D_\ell$ by identifying $s_1,s_2,...,s_\ell$ into $s$ and identifying $t_1,t_2,...,t_\ell$ into $t$.

\begin{lemma} \label{sp parallel wqo}
Let $\F$ be a well-behaved set of one-way series-parallel triples.
Let $\F'$ be the parallel-extension of $\F$.
Then $\F'$ is well-behaved.
\end{lemma}

\begin{pf}
Let $(Q,\preceq)$ be a well-quasi-order.
For each $i \in {\mathbb N}$, let $(D_i,s_i,t_i)$ be a member of $\F'$, and let $\phi_i: V(D_i) \rightarrow Q$.
For each $i \in {\mathbb N}$, since $(D_i,s_i,t_i) \in \F'$, there exist $\ell_i \in {\mathbb N}$ and members $(D_{i,1},s_{i,1},t_{i,1}), (D_{i,2},s_{i,2},t_{i,2}),...,(D_{i,\ell_i},s_{\ell_i},t_{\ell_i})$ of $\F$ such that $D_i$ is obtained from the disjoint union of $D_{i,1},D_{i,2},...,D_{i,\ell_i}$ by identifying $s_{i,1},s_{i,2},...,s_{i,\ell_i}$ into $s_i$ and identifying $t_{i,1},t_{i,2},...,t_{i,\ell_i}$ into $t_i$.

For each $i \in {\mathbb N}$, let $a_i$ be the sequence $((D_{i,1},s_{i,1},t_{i,1}),\phi_i|_{V(D_{i,1})}), ((D_{i,2},s_{i,2},t_{i,2}),\phi_i|_{V(D_{i,2})}), \allowbreak ..., ((D_{i,\ell_i},s_{i,\ell_i},t_{i,\ell_i}),\phi_i|_{V(D_{i,\ell_i})})$.
Since $\F$ is well-behaved, by Higman's Lemma, there exist $1 \leq j <j'$ and a strictly increasing function $f: [\ell_j] \rightarrow [\ell_{j'}]$ such that for every $i \in [\ell_j]$, $((D_{j',f(i)},s_{j',f(i)},t_{j',f(i)}),\phi_{j'}|_{V(D_{j',f(i)})})$ simulates $((D_{j,i},s_{j,i},t_{j,i}),\phi_j|_{V(D_{j,i})})$.
Hence $((D_{j'},s_{j'},t_{j'}),\phi_{j'})$ simulates $((D_j,s_j,t_j),\phi_j)$.
Therefore, $\F'$ is well-behaved.
\end{pf}

\bigskip

For a set $\F$ of one-way series-parallel triples, the {\it series-extension} of $\F$ is a set $\F'$ of one-way series-parallel triples such that for every $(D,s,t) \in \F'$, there exist a positive integer $\ell$ and members $(D_1,s_1,t_1), (D_2,s_2,t_2), ..., (D_\ell,s_\ell,t_\ell)$ of $\F$ such that $D$ is obtained from the disjoint union of $D_1,D_2,...,D_\ell$ by for each $i \in [\ell-1]$, identifying $t_i$ and $s_{i+1}$.

\begin{lemma} \label{sp series wqo}
Let $\F$ be a well-behaved set of one-way series-parallel triples.
Let $\F_1$ be the series-extension of $\F$.
Let $\F_2$ be the set consisting of all series-parallel triples that are truncations of members of $\F$.
If $\F_2$ is well-behaved, then $\F_1$ is well-behaved.
\end{lemma}

\begin{pf}
Let $(Q,\preceq)$ be a well-quasi-order.
For each $i \in {\mathbb N}$, let $(D_i,s_i,t_i)$ be a member of $\F_1$, and let $\phi_i: V(D_i) \rightarrow Q$.
For each $i \in {\mathbb N}$, since $(D_i,s_i,t_i) \in \F_1$, there exist $\ell_i \in {\mathbb N}$ and members $(D_{i,1},s_{i,1},t_{i,1}), (D_{i,2},s_{i,2},t_{i,2}),...,(D_{i,\ell_i},s_{\ell_i},t_{\ell_i})$ of $\F$ such that $D_i$ is obtained from the disjoint union of $D_{i,1},D_{i,2},...,D_{i,\ell_i}$ by for each $j \in [\ell_i-1]$, identifying $t_j$ and $s_{j+1}$.
To prove this lemma, it suffices to prove that there exist $1 \leq j <j'$ such that $((D_{j'},s_{j'},t_{j'}),\phi_{j'})$ simulates $((D_j,s_j,t_j),\phi_j)$.

Since each $(D_i,s_i,t_i)$ is one-way by the definition of a series-extension, by symmetry and possibly removing some members in the sequence, we may assume that for each $i \in {\mathbb N}$, every thread in $D_i$ between $s_i$ and $t_i$ is a directed path in $D_i$ from $s_i$ to $t_i$.
Let $\F'$ be the set consisting of the rooted digraphs $(D,s)$ such that $(D,s,t) \in \F$ for some $t \in V(D)$.
Note that for each $i \in {\mathbb N}$, $(D_i,s_i)$ is a $\F'$-series-parallel tree.

Since $\F$ is a well-behaved set of one-way series-parallel triples, $\F'$ is a well-behaved set of rooted digraphs.
Since $\F$, $\F'$ and $\F_2$ are well-behaved, the set of $\F'$-series parallel trees is well-behaved by Lemma \ref{series parallel tree wqo}.

Let $(Q',\preceq')$ be the well-quasi-order obtained by the Cartesian product of $(Q,\preceq)$ and $([2],=)$.
For each $i \in {\mathbb N}$, let $\phi'_i: V(Q') \rightarrow V(D_i)$ such that $\phi'_i(t_i)=(\phi_i(t_i),2)$, and $\phi'_i(v)=(\phi_i(v),1)$ for every $v \in V(D_i)-\{t_i\}$.
Since the set of $\F'$-series parallel trees is well-behaved, there exist $1 \leq j <j'$ such that there exists a strong immersion embedding $\eta$ from $(D_j,s_j)$ to $(D_{j'},s_{j'})$ such that $\phi'_j(v) \preceq' \phi'_{j'}(\eta(v))$ for every $v \in V(D_j)$.

By the definition of $\phi'_j$ and $\phi'_{j'}$, $\eta(t_j)=t_{j'}$.
So $((D_{j'},s_{j'},t_{j'}),\phi_{j'})$ simulates $((D_j,s_j,t_j),\phi_j)$.
This proves the lemma.
\end{pf}

\bigskip

Define $\A_{0}$ to be the set of one-way series-parallel triples $(D,s,t)$ such that $D$ consists of an edge. 
Define $\A_{0,0}=\A_0$.
For any nonnegative integers $k$ and $i$, we define the following.
	\begin{itemize}
		\item Define $\A_{k,2i+1}$ to be the parallel-extension of $\A_{k,2i}$. 
		\item Define $\A_{k,2i+2}$ to be the series-extension of $\A_{k,2i+1}$. 
		\item Define $\A_{k+1}$ to be the set of one-way series-parallel triples $(D,s,t)$ such that there exists no $(k+1)$-alternating path in $D$ with one end $s$ or one end $t$.
		\item Define $\A_{k+1,0}$ to be the set of one-way series-parallel triples $(D,s,t)$ such that either
			\begin{itemize}
				\item every $(k+1)$-alternating path in $D$ with one end $s$ intersects $t$, and there exists no $(k+1)$-alternating path in $D$ with one end $t$, or
				\item every $(k+1)$-alternating path in $D$ with one end $t$ intersects $s$, and there exists no $(k+1)$-alternating path in $D$ with one end $s$.
			\end{itemize}
	\end{itemize}

\begin{lemma} \label{sp full constr}
For every nonnegative integer $k$, $\A_{k+1} \subseteq \A_{k,4}$.
\end{lemma}

\begin{pf}
By Lemma \ref{1-alt structure}, $\A_1 \subseteq \A_{0,2} \subseteq \A_{0,4}$.
So we may assume $k \geq 1$.

\medskip

\noindent{\bf Claim 1:} Every series-irreducible one-way series-parallel triple in $\A_{k+1}$ belongs to $\A_{k,3}$.

\medskip

\noindent{\bf Proof of Claim 1:}
Suppose to the contrary that there exists a series-irreducible one-way series-parallel triple $(D,s,t)$ in $\A_{k+1}-\A_{k,3}$.
If $\lvert E(D) \rvert=1$, then $(D,s,t) \in \A_{k,0} \subseteq \A_{k,3}$, a contradiction.
By Lemma \ref{constr_sp}, $(D,s,t)$ is not parallel-irreducible.
So there exist $\ell \in {\mathbb N}$ with $\ell \geq 2$ and parallel-irreducible one-way series-parallel triples $(D_1,s_1,t_1),...,(D_\ell,s_\ell,t_\ell)$ such that $D$ is obtained from a disjoint union of $D_1,...D_\ell$ by identifying $s_1,...,s_\ell$ into $s$ and identifying $t_1,...,t_\ell$ into $t$.
Since $(D,s,t) \not \in \A_{k,3}$, there exists $i \in [\ell]$ such that $(D_i,s_i,t_i) \not \in \A_{k,2}$.
By symmetry, we may assume that $(D_1,s_1,t_1) \not \in \A_{k,2}$.
In particular, $\lvert E(D_1) \rvert \geq 2$.
Since $(D_1,s_1,t_1)$ is parallel-irreducible, by Lemma \ref{constr_sp}, there exist $\ell_1 \in {\mathbb N}$ with $\ell_1 \geq 2$ and series-irreducible one-way series-parallel triples $(D_{1,1},s_{1,1},t_{1,1}),...,(D_{1,\ell_1},s_{1,\ell_1},t_{1,\ell_1})$ such that $D$ is obtained from a disjoint union of $D_{1,1},...D_{1,\ell_1}$ by for each $j \in [\ell_1-1]$, identifying $t_{1,j}$ with $s_{1,j+1}$.
Since $(D_1,s_1,t_1) \not \in \A_{k,2}$, there exists $j \in [\ell_1]$ such that $(D_{1,j},s_{1,j},t_{1,j}) \not \in \A_{k,1}$.
Since $\ell_1 \geq 2$, by possibly switching $s$ and $t$, switching $s_j$ and $t_j$, and replacing $j$ by $\ell_1+1-j$, we may assume that $s_{1,j} \neq s=s_1$.

If there exists a $k$-alternating path $P$ in $D_{1,j}$ with one end $s_{1,j}$ disjoint from $t_{1,j}$, then by concatenating $P$ with a thread in $D_2$ from $t_2=t \not \in V(P)$ to $s_2=s$ and a thread in $D_1$ from $s$ to $s_{1,j}$, we obtain a $(k+1)$-alternating path in $D$ with one end $t$, so $(D,s,t) \not \in \A_{k+1}$, a contradiction.

So every $k$-alternating path $P$ in $D_{1,j}$ with one end $s_{1,j}$ intersects $t_{1,j}$.
If there exists a $k$-alternating path $P$ in $D_{1,j}$ with one end $t_{1,j}$, then by concatenating $P$ with a thread in $D_2$ from $s_2=s \not \in V(D_{1,j})$ to $t_2=t$ and a thread in $D_1$ from $t$ to $t_{1,j}$, we obtain a $(k+1)$-alternating path in $D$ with one end $s$, so $(D,s,t) \not \in \A_{k+1}$, a contradiction.
So no $k$-alternating path $P$ in $D_{1,j}$ has one end $t_{1,j}$.
Hence $(D_{1,j},s_{1,j},t_{1,j}) \in \A_{k,0} \subseteq \A_{k,1}$, a contradiction.
$\Box$

\medskip

Now we prove that $\A_{k+1} \subseteq \A_{k,4}$.

Suppose to the contrary that there exists a series-parallel triple $(D',s',t') \in \A_{k+1}-\A_{k,4}$.
Since $\A_{k,3} \subseteq \A_{k,4}$, $(D',s',t')$ is not series-irreducible by Claim 1.
By Lemma \ref{constr_sp}, there exist $\ell' \in {\mathbb N}$ with $\ell' \geq 2$ and series-irreducible one-way series-parallel triples $(D_1',s_1',t_1'),...,(D_{\ell'}',s_{\ell'}',t_{\ell'}')$ such that $D'$ is obtained from the disjoint union of $D_1',...,D_{\ell'}'$ by for each $i \in [\ell'-1]$ identifying $t'_i$ with $s'_{i+1}$.
Since $(D',s',t') \not \in \A_{k,4}$, there exists $i^* \in [\ell']$ such that $(D'_{i^*},s'_{i^*},t'_{i^*}) \not \in \A_{k,3}$.
Since $(D'_{i^*},s_{i^*}',t_{i^*}')$ is series-irreducible, it is not in $\A_{k+1}$ by Claim 1.
So by symmetry, we may assume that there exists a $(k+1)$-alternating path in $D'_{i^*}$ with one end $t_{i^*}'$.
But then we can extend it to a $(k+1)$-alternating path in $D'$ with one end $t'$, contradicting that $(D',s',t') \in \A_{k+1}$.
This proves the lemma.
\end{pf}

\begin{lemma} \label{truncate no boundary alternating}
Let $k$ be a positive integer.
Let $(D,s,t)$ be a one-way series-parallel triple. 
Let $[X,Y]$ be a splitter of $(D,s,t)$.
	\begin{enumerate}
		\item If $(D,s,t) \in \A_k$, then every truncation of $(D,s,t)$ with respect to $[X,Y]$ belongs to $\A_k$.
		\item If $(D,s,t) \in \A_{k,0}$, then every truncation of $(D,s,t)$ with respect to $[X,Y]$ belongs to $\A_{k,0}$.
	\end{enumerate}
\end{lemma}

\begin{pf}
Let $(D_X,s,t_Y)$ be the series-parallel triple such that $D_X$ is obtained from $D$ by identifying all vertices in $Y$ into a vertex $t_Y$ and deleting all resulting loops.
By symmetry, it suffices to prove that if $(D,s,t) \in \A_k$ then $(D_X,s,t_Y) \in \A_k$, and if $(D,s,t) \in \A_{k,0}$, then $(D_X,s,t_Y) \in \A_{k,0}$.

Let $w$ be the number of edges of $D$ between $X$ and $Y$.
Since $[X,Y]$ is a splitter, there exist $w$ edge-disjoint threads in $D$ from $s$ to $t$ such that every edge between $X$ and $Y$ belongs to exactly one of those threads.
For each edge $e$ between $X$ and $Y$, we denote the aforementioned thread containing $e$ by $P_e$.
Since $(D,s,t)$ is one-way, $P_e$ is a directed path for each $e$.

Let $P$ be a $k$-alternating path in $D_X$ with one end $s$ or one end $t_Y$.
So $P$ contains at most two edges incident with $t_Y$.
Note that every edge incident with $t_Y$ is an edge of $D$ between $X$ and $Y$.
If $P$ contains no edge incident with $t_Y$ or contains two edges incident with $t_Y$ that has a common end in $Y$, then $P$ is a $k$-alternating in $D$ with one end $s$. 
If $P$ contains exactly one edge incident with $t_Y$, then one can concatenate $P$ with a thread in $Y$ to obtain a $k$-alternating path in $D$ with one end $t$. 
If $P$ contains exactly two edges $e_1,e_2$ incident with $t_Y$, and the ends of these two edges in $Y$ are distinct, then $P \cup P_{e_1} \cup P_{e_2}$ contains a $k$-alternating path in $D$ with one end $s$, since $t_Y$ is not an end of $P$. 

Hence it is straightforward to verify that if $(D,s,t) \in \A_k$, then $(D_X,s,t_Y) \in \A_k$; and if $(D,s,t) \in \A_{k,0}$, then $(D_X,s,t_Y) \in \A_{k,0}$.
\end{pf}

\begin{lemma} \label{sp wqo alt leng}
Let $k$ be a nonnegative integer.
If $\A_k$ is well-behaved, then $\A_{k,0}$ is well-behaved.
\end{lemma}

\begin{pf}
When $k=0$, $\A_{k,0}=\A_k$.
So we may assume $k \geq 1$.

Let $(Q,\preceq)$ be a well-quasi-order.
For each $i \in {\mathbb N}$, let $(D_i,s_i,t_i)$ be a member of $\A_{k,0}$, and let $\phi_i: V(D_i) \rightarrow Q$.
To prove this lemma, it suffices to prove that there exist $1 \leq j <j'$ such that $((D_{j'},s_{j'},t_{j'}),\phi_{j'})$ simulates $((D_j,s_j,t_j),\phi_j)$.

By symmetry and possibly removing some $(D_i,s_i,t_i)$, we may assume that for each $i \in {\mathbb N}$, every thread in $D_i$ is a directed path from $s_i$ to $t_i$, every $k$-alternating path in $D_i$ with one end $s_i$ intersects $t_i$, and there exists no $k$-alternating path in $D_i$ with one end $t_i$.
Let $\F = \{(D,s):$ there exists $t \in V(D)-\{s\}$ such that $(D,s,t) \in \A_k\}$.

\medskip

\noindent{\bf Claim 1:} For each $i \in {\mathbb N}$, $(D_i-t_i,s_i)$ is a $\F$-series-parallel tree.

\medskip

\noindent{\bf Proof of Claim 1:}
Let $i$ be a fixed positive integer.
Since $(D_i,s_i,t_i)$ is a series-parallel triple, there exists no separation $(A_1,A_2)$ of $D_i$ of order at most one such that $\{s_i,t_i\} \subseteq V(A_1)$ and $V(A_2)-V(A_1) \neq \emptyset$.
So for every block $B$ of $D_i-t_i$ that has no child block, $t_i$ is adjacent in $D_i$ to a vertex in $V(B)-\{v\}$, where $v=s_i$ if $s_i \in V(B)$, and $v$ is the vertex contained in $B$ and its parent block if $s_i \not \in V(B)$.
Hence no block $B$ of $D_i-t_i$ has two child blocks intersecting $B$ at different vertices, for otherwise there exists a thread in $D_i$ from $s_i$ to $t_i$ that is not a directed path, a contradiction.
Similarly, for every cut-vertex $v$ of the underlying graph of $D_i-t_i$, every thread in $D_i-t_i$ from $s_i$ to $v$ is a directed path from $s_i$ to $v$.
Hence $(D_i-t_i,s_i)$ is a $\F_0$-series-parallel tree for some set $\F_0$ of one-way series-parallel triples.
Since every $k$-alternating path in $D_i$ with one end $s_i$ intersects $t_i$, and there exists no $k$-alternating path in $D_i$ with one end $t_i$, we know $\F_0 \subseteq \F$.
$\Box$

\medskip

Let $(Q',\preceq')$ be the well-quasi-order obtained by the Cartesian product of $(Q,\preceq)$, $([2],=)$, $({\mathbb N}, \leq )$ and $(Q,\preceq)$.
For each $i \in {\mathbb N}$, let $\phi'_i: V(D_i-t_i) \rightarrow Q'$ such that for every $v \in V(D_i-t_i)$, if $v$ is not adjacent in $D_i$ to $t_i$, then $\phi'_i(v)=(\phi_i(v),1,1,\phi_i(t_i))$, otherwise $\phi'_i(v)=(\phi_i(v),2,d_v,\phi_i(t_i))$, where $d_v$ is the number of edges of $D_i$ between $v$ and $t_i$.

Let $\F_1$ be the set consisting of all series-parallel triples that are truncations of members of $\A_k$.
By Lemma \ref{truncate no boundary alternating}, $\F_1 \subseteq \A_k$.
Since $\A_k$ is well-behaved, $\F$ and $\F_1$ are well-behaved.
By Lemma \ref{series parallel tree wqo}, there exist $1 \leq j <j'$ and a strong immersion embedding $\eta$ from $(D_j-t_j,s_j)$ to $(D_{j'}-t_{j'},s_{j'})$ such that for every $v \in V(D_j-t_j)$, $\phi'_j(v) \preceq' \phi'_{j'}(\eta(v))$.
By the definition of $\phi'_j$ and $\phi'_{j'}$, for every neighbor $u$ of $t_j$ in $D_j$, $\eta(u)$ is a neighbor of $t_{j'}$ in $D_{j'}$, and the number of edges of $D_{j'}$ between $\eta(u)$ and $t_{j'}$ is at least the number of edges of $D_j$ between $u$ and $t_j$.
Since $t_j$ and $t_{j'}$ are sinks, one can extend $\eta$ to be a strong immersion embedding $\eta^*$ from $(D_j,s_j,t_j)$ to $(D_{j'},s_{j'},t_{j'})$ such that $\eta^*(s_j)=s_{j'}$, $\eta^*(t_j)=t_{j'}$, and for every $v \in V(D_j)$, $\phi_j(v) \preceq \phi_{j'}(\eta(v))$.
This proves the lemma.
\end{pf}

\begin{lemma} \label{sp truncation close}
Let $k$ be a positive integer.
Let $a$ be a nonnegative integer.
Let $\F$ be the set of series-parallel triples that are truncations of members of $\A_{k,a}$.
Then $\F \subseteq \A_{k,a}$.
\end{lemma}

\begin{pf}
Let $(D,s,t) \in \A_{k,a}$.
Let $[X,Y]$ be a splitter of $(D,s,t)$.
To prove this lemma, it suffices to show that every truncation of $(D,s,t)$ respect to $[X,Y]$ belongs to $\A_{k,a}$.

Let $(D_X,s,t_Y)$ be the series-parallel triple such that $D_X$ is obtained from $D$ by identifying $Y$ into the vertex $t_Y$ and deleting resulting loops.
By symmetry, it suffices to prove that $(D_X,s,t_Y) \in \A_{k,a}$.
We shall prove it by induction on $a$.
The case $a=0$ follows from Statement 2 of Lemma \ref{truncate no boundary alternating}.
So we may assume that $a \geq 1$, and this lemma holds when $a$ is smaller.

We first assume that $a$ is odd.
So there exist $\ell \in {\mathbb N}$ and $(D_1,s_1,t_1),...,(D_\ell,s_\ell,t_\ell) \in \A_{k,a-1}$ such that $D$ is obtained from the disjoint union of $D_1,...,D_\ell$ by identifying $s_1,s_2,...,s_\ell$ into $s$ and identifying $t_1,t_2,...,t_\ell$ into $t$.
For each $i \in [\ell]$, let $X_i=X \cap V(D_i)$ and $Y_i=Y \cap V(D_i)$, and let $(D_{i,X},s_i,t_{i,Y})$ be the series-parallel triple such that $D_{i,X}$ is obtained from $D_i$ by identifying $Y_i$ into a vertex $t_{i,Y}$ and deleting resulting loops.
Note that for each $i \in [\ell]$, $[X_i,Y_i]$ is a splitter of $(D_i,s_i,t_i)$, so $(D_{i,X},s_i,t_{i,Y}) \in \A_{k,a-1}$ by the induction hypothesis.
And $D_X$ is obtained from a disjoint union of $D_{1,X},...,D_{\ell,X}$ by identifying $s_1,...,s_\ell$ into $s$ and identifying $t_{1,Y},...,t_{\ell,Y}$ into $t_Y$.
So $(D_X, s,t_Y) \in \A_{k,a}$.

Hence we may assume that $a$ is even.
So there exist $\ell' \in {\mathbb N}$ and $(D'_1,s'_1,t'_1),...,(D'_{\ell'},s'_{\ell'},t'_{\ell'}) \in \A_{k,a-1}$ such that $D$ is obtained from the disjoint union of $D'_1,...,D'_{\ell'}$ by for each $i \in [\ell'-1]$, identifying $t'_i$ and $s'_{i+1}$.
Since $[X,Y]$ is a splitter, there exists $\ell^* \in [\ell']$ such that all edges of $D$ between $X$ and $Y$ are edges of $D_{\ell^*}$.
Let $X'=X \cap V(D_{\ell^*}')$ and $Y' = Y \cap V(D'_{\ell^*})$.
Note that $[X',Y']$ is a splitter of $(D'_{\ell^*},s'_{\ell^*},t_{\ell^*}')$.
Let $(D'_{\ell^*,X},s'_{\ell^*},t'_{\ell^*,Y})$ be the series-parallel triple such that $D'_{\ell^*,X}$ is obtained from $D'_{\ell^*}$ by identifying $Y'$ into a vertex $t'_{\ell^*,Y}$ and deleting resulting loops.
By the induction hypothesis, $(D'_{\ell^*,X},s'_{\ell^*},t'_{\ell^*,Y}) \in \A_{k,a-1}$.
Note that $D_X$ is obtained from a disjoint union of $D'_{1},...,D'_{\ell^*-1}, D'_{\ell^*,X}$ by for each $i \in [\ell^*-1]$, identifying $t'_i$ with $s'_{i+1}$.
So $(D_X, s,t'_Y) \in \A_{k,a}$.
\end{pf}

\begin{lemma} \label{sp wqo 1}
For any nonnegative integers $k$ and $a$, if $\A_{k,0}$ is well-behaved, then $\A_{k,a}$ is well-behaved.
\end{lemma}

\begin{pf}
We shall prove this lemma by induction on $a$.
When $a=0$, $\A_{k,0}$ is well-behaved.
So we may assume that $a \geq 1$, and $\A_{k,a-1}$ is well-behaved.
If $a$ is odd, then $\A_{k,a}$ is well-behaved by Lemma \ref{sp parallel wqo}.
If $a$ is even, then by Lemmas \ref{sp truncation close} and \ref{sp series wqo}, $\A_{k,a}$ is well-behaved.
\end{pf}

\begin{lemma} \label{sp wqo 2}
For every nonnegative integer $k$, $\A_k$ is well-behaved.
\end{lemma}

\begin{pf}
We shall prove this lemma by induction on $k$.
When $k=0$, $\A_k$ is clearly well-behaved.
So we may assume that $k \geq 1$, and $\A_{k-1}$ is well-simulated.
By Lemma \ref{sp wqo alt leng}, $\A_{k-1,0}$ is well-behaved.
By Lemma \ref{sp wqo 1}, $\A_{k-1,4}$ is well-behaved.
By Lemma \ref{sp full constr}, $\A_k \subseteq \A_{k-1,4}$, so $\A_k$ is well-behaved.
\end{pf}

\begin{lemma} \label{series-parallel triple wqo}
Let $(Q,\preceq)$ be a well-quasi-order.
Let $k$ be a positive integer.
For each $i \in {\mathbb N}$, let $(D_i,s_i,t_i)$ be a one-way series-parallel triple such that $D_i$ does not contain a $k$-alternating path, and let $\phi_i: V(D_i) \rightarrow Q$.
Then there exist $1 \leq j <j'$ and a strong immersion embedding $\eta$ from $(D_{j},s_{j},t_{j})$ to $(D_{j'},s_{j'},t_{j'})$ such that $\phi_j(x) \preceq \phi_{j'}(\eta(x))$ for every $x \in V(D_j)$.
\end{lemma}

\begin{pf}
For each $i \in {\mathbb N}$, since $D_i$ has no $k$-alternating path, $(D_i,s_i,t_i) \in \A_k$.
So this lemma immediately follows from Lemma \ref{sp wqo 2}.
\end{pf}

\section{Series-parallel separations} \label{sec:sp_sep}

The goal of this section is to prove Lemmas \ref{ep path} and \ref{collection sp sep}.
Roughly speaking, they show that if a digraph with no $(k+1)$-alternating path is not of a very special form, then we can find a cross-free collection $\Se$ of separations and a set $Z$ of vertices of bounded size such that every $k$-alternating path either intersects $Z$ or is contained in the series-parallel part of a separation in $\Se$.

We remark that an obvious simplification of the statement mentioned in the previous paragraph is false.
That is, there exists no function $f$ such that every digraph with no $(k+1)$-alternating path can be made a digraph with no $k$-alternating path by deleting at most $f(k)$ vertices.
Consider the series-parallel triple $(D,s,t)$ obtained by the series-concatenation of any number of copies of the one-way series-parallel triple whose underlying graph is $K_{2,3}$ and the roots are the two vertices of degree two.
It is not hard to see that $D$ has no $3$-alternating path.
But there is a 2-alternating path in $D$ contained in each oriented $K_{2,3}$.
As the number of copies of $K_{2,3}$ can be arbitrarily large, $D$ can contain arbitrarily many disjoint 2-alternating paths.
Hence 2-alternating paths in $D$ cannot be killed by deleting a bounded number of vertices.

We need some simple lemmas in order to prove Lemma \ref{ep path}.

For positive integers $m$ and $n$, the {\it $m \times n$-grid} is the graph with vertex-set $[m] \times [n]$ such that any vertices $(x,y)$ and $(x',y')$ are adjacent if and only if $\lvert x-x' \rvert + \lvert y-y' \rvert =1$.

\begin{lemma}
Let $t$ be a positive integer.
Let $D$ be a digraph.
If some subgraph of the underlying graph of $D$ is isomorphic to a subdivision of the $2 \times (2t+1)$-grid, then $D$ contains a $t$-alternating path.
\end{lemma}

\begin{pf}
Since some subgraph of the underlying graph of $D$ is isomorphic to a subdivision of the $2 \times (2t+1)$-grid, there exist two disjoint threads $R_1$ and $R_2$ and $2t+1$ disjoint threads $P_1,P_2,...,P_{2t+1}$ from $V(R_1)$ to $V(R_2)$ internally disjoint from $V(R_1) \cup V(R_2)$ such that for each $i \in [2]$, $R_i$ passes through $x_{i,1},x_{i,2},...,x_{i,2t+1}$ in the order listed, where $V(R_i) \cap \bigcup_{j=1}^{2t+1}V(P_j)=\{x_{i,1},x_{i,2},...,x_{i,2t+1}\}$. 
For each $i \in [2]$ and $k \in [t]$, let $R_{i,k}$ be the subthread of $R_i$ between $x_{i,k}$ and $x_{i,2t+1}$.

To prove this lemma, it suffices to prove that for each $k \in [t+1]$ and $v \in \{x_{1,2k-1}, x_{2,2k-1}\}$, there exists a $(t+1-k)$-alternating path in $D[V(R_{1,2k-1}) \cup V(R_{2,2k-1}) \cup \bigcup_{j=2k-1}^{2t+1}P_j]$ with one end $v$.
We shall prove it by induction on $t+1-k$.

When $t+1-k=0$, we have that $k=t+1$ and $P_{2t+1}$ is a 0-alternating path in $D[V(R_{1,2k-1}) \cup V(R_{2,2k-1}) \cup \bigcup_{j=2k-1}^{2t+1}P_j]$ with ends $x_{1,2t+1}$ and $x_{2,2t+1}$.
So we may assume that $k \in [t]$, and the claim holds when $t+1-k$ is smaller.

By symmetry, we may assume that $v=x_{1,2k-1}$.
By the induction hypothesis, for each $i \in [2]$, there exists a $(t-k)$-alternating path $P^*_i$ in $D[V(R_{1,2k+1}) \cup V(R_{2,2k+1}) \cup \bigcup_{j=2k+1}^{2t+1}P_j]$ with an end $x_{i,2k+1}$.

Let $W_1$ be subthread of $R_1$ between $x_{1,2k-1}$ and $x_{1,2k+1}$.
Let $W_2$ be the thread obtained from $P_{2k-1} \cup P_{2k}$ by concatenating the subthread of $R_2$ between $x_{2,2k-1}$ and $x_{2,2k}$ and the subthread of $R_1$ between $x_{1,2k}$ and $x_{1,2k+1}$.
Let $W_3$ be the thread obtained from the subthread of $R_1$ between $x_{1,2k-1}$ and $x_{1,2k}$ by concatenating $P_{2k}$ and the subthread of $R_2$ between $x_{2,2k}$ and $x_{2,2k+1}$.
Then it is straight forward to show that there exists $i^* \in [3]$ such that $W_{i^*}$ does not induce a directed path in $D$.
Hence $W_{i^*} \cup P^*_1 \cup P^*_2$ contains a $(t-k+1)$-alternating path with end $x_{1,2k-1}$. 
Moreover, $W_{i^*} \cup P^*_1 \cup P^*_2$ is contained in $D[V(R_{1,2k-1}) \cup V(R_{2,2k-1}) \cup \bigcup_{j=2k-1}^{2t+1}P_j]$.
This proves the lemma.
\end{pf}

\bigskip

The {\it tree-width} of a graph $G$ is the minimum $w$ such that $G$ is a subgraph of a graph with no induced cycle of length at least four and with no clique of size $w+2$.
The following is a restatement of a lemma about the Erd\H{o}s-P\'{o}sa property proved in \cite{t_ep}.

\begin{lemma}[{{\cite[Proposition 2.1]{t_ep}}}] \label{tw ep}
Let $w$ be a positive integer.
Let $G$ be a graph with tree-width at most $w$.
Let $\F$ be a collection of subsets of $V(G)$ such that for every $S \in \F$, $G[S]$ is connected.
Let $k$ be a positive integer.
If there do not exist $k$ pairwise disjoint members of $\F$, then there exists $Z \subseteq V(G)$ with $\lvert Z \rvert \leq (k-1)(w+1)$ such that $Z \cap S \neq \emptyset$ for every $S \in \F$.
\end{lemma}

Let $D$ be a digraph.
A {\it series-parallel 2-separation} of $D$ is a separation $(A,B)$ of $D$ such that $(A, s,t)$ is a one-way series-parallel triple, where $V(A\cap B) =\{s,t\}$.
The following lemma is a slight strengthening of \cite[Theorem 5.6]{m}, and our proof is a modification of \cite[Theorem 5.6]{m}.

\begin{lemma} \label{ep path}
For every positive integer $t$, there exists an integer $f(t)$ such that the following holds.
If $D$ is a digraph whose underlying graph is 2-connected, and $D$ does not contain a $(t+1)$-alternating path, then there exists $Z \subseteq V(D)$ with $\lvert Z \rvert \leq f(t)$ such that for every $t$-alternating path $P$ in $D$, either
	\begin{enumerate}
		\item $V(P) \cap Z \neq \emptyset$, or
		\item there exists a series-parallel 2-separation $(A,B)$ of $D$ with $P \subseteq A$.
	\end{enumerate}
\end{lemma}

\begin{pf}
Let $t$ be a positive integer.
Since $D$ does not contain a $t$-alternating path, the underlying graph of $D$ does not contain a subdivision of the $2 \times (2t+1)$-grid.
Since the $2 \times (2t+1)$-grid is subcubic, there exists an integer $w$ such that the underlying graph of $D$ has tree-width at most $w$ by the Grid Minor Theorem \cite{rs_V}.

Define $f(t)=4(w+1)$.
Note that $w$ only depends on $t$, so $f(t)$ only depends on $t$.

For a $t$-alternating path $P$ in $D$,
	\begin{itemize}
		\item if $t$ is odd, then let $m_P$ be the $\lceil \frac{t}{2} \rceil$-th pivot of $P$. 
		\item if $t$ is even, then let $m_P$ and $m_P'$ be the $\frac{t}{2}$-th and $(\frac{t}{2}+1)$-th pivots of $P$, respectively, and by symmetry, we may assume that $m_P$ is a sink and $m'_P$ is a source.
	\end{itemize}

\medskip

\noindent{\bf Claim 1:} Let $P_1$ and $P_2$ be disjoint $t$-alternating paths in $D$.
Then for every thread $P$ in $D$ intersecting $V(P_1)$ and $V(P_2)$ internally disjoint from $V(P_1) \cup V(P_2)$,
	\begin{itemize}
		\item if $t$ is odd, then there exists $i \in [2]$ such that 
			\begin{itemize}
				\item $V(P_i) \cap V(P)=\{m_{P_i}\}$, 
				\item the vertex in $V(P_{3-i}) \cap V(P)$ belongs to the sub-thread of $P_{3-i}$ between the $(\lceil \frac{t}{2} \rceil-1)$-th pivot and the $(\lceil \frac{t}{2} \rceil+1)$-th pivot, and
				\item if $V(P_{3-i}) \cap V(P) \neq \{m_{P_{3-i}}\}$, then $P$ is a directed path;
			\end{itemize}
		\item if $t$ is even, then $P$ is a directed path, and there exists $i \in [2]$ such that $V(P) \cap V(P_i) = \{m_{P_i}\}$ and $V(P) \cap V(P_{3-i})=\{m'_{P_i}\}$. 
	\end{itemize}

\medskip

\noindent{\bf Proof of Claim 1:}
For each $i \in [2]$, let $v_i$ be the end of $P$ in $V(P_i)$, and let $Q_i$ be a sub-thread of $P_i$ between $v_i$ and an end of $P_i$ such that the number of pivots is as large as possible.
Note that if $t$ is even, then $Q_i$ contains at least $\frac{t}{2}$ pivots, and the equality holds only when $v_i$ is contained in the sub-thread of $P_i$ between $m_{P_i}$ and $m'_{P_i}$; if $t$ is odd, then $Q_i$ contains at least $\lceil \frac{t}{2} \rceil-1$ pivots, and the equality holds only when $v_i=m_{P_i}$.

We first assume that $t$ is odd.
If $v_1 \neq m_{P_1}$ and $v_2 \neq m_{P_2}$, then the number of pivots in $Q_1 \cup P \cup Q_2$ is at least $2 \cdot \lceil \frac{t}{2} \rceil \geq t+1$, a contradiction.
So there exists $i \in [2]$ such that $v_i=m_{P_i}$.
That is, $V(P_i) \cap V(P) = \{m_{P_i}\}$.
Similarly, $Q_{3-i}$ contains at most $\lceil \frac{t}{2} \rceil$ pivots, so $v_{3-i}$ belongs to the sub-thread of $P_{3-i}$ between the $(\lceil \frac{t}{2} \rceil-1)$-th and the $(\lceil \frac{t}{2} \rceil+1)$-th pivots.
If $V(P_{3-i}) \cap V(P) \neq \{m_{P_{3-i}}\}$, then the number of pivots in $Q_1 \cup P \cup Q_2$ is at least $2\lceil \frac{t}{2} \rceil-1$ plus the number of pivots of $P$, so $P$ has no pivots.
This proves the case when $t$ is odd.

Now we assume that $t$ is even.
Then the number of pivots of $Q_1 \cup P \cup Q_2$ is at least $2\cdot \frac{t}{2}$ plus the number of pivots of $P$.
Since there exists no $(t+1)$-alternating path in $D$, $P$ is a directed path, and for every $i \in [2]$, $v_i$ is contained in the sub-thread of $P_i$ between $m_{P_i}$ and $m'_{P_i}$.
Suppose to the contrary that there exists $j \in [2]$ such that $v_j \not \in \{m_{P_j},m'_{P_j}\}$.
So $v_j$ is an internal vertex of a directed subpath of $P_j$ between $m_{P_j}$ and $m'_{P_j}$.
Hence we can choose $Q_j$ such that $v_j$ is a pivot of $Q_j \cup P$.
Therefore, $Q_j \cup P \cup Q_{3-j}$ contains at least $t+1$ pivots, a contradiction.
So for every $i \in [2]$, $v_i \in \{m_{P_i},m_{P_i'}\}$.
For $i \in [2]$, if $m_{P_i} \in V(P)$, then since $m_{P_i}$ is a sink in $P_i$, $m_{P_i}$ is the source of $P$, for otherwise there exists a $(t+1)$-alternating path in $P_1 \cup P \cup P_2$; similarly, if $m'_{P_i} \in V(P)$, then since $m'_{P_i}$ is a source in $P_i$, $m'_{P_i}$ is the sink of $P$.
Hence there exists $i \in [2]$ such that $V(P) \cap V(P_i) = \{m_{P_i}\}$ and $V(P) \cap V(P_{3-i})=\{m'_{P_i}\}$.
$\Box$

\medskip

\noindent{\bf Claim 2:} Let $P_1$ and $P_2$ be two disjoint $t$-alternating paths in $D$.
Then there exist a separation $(A,B)$ of $D$ of order two such that $P_1 \subseteq A$ and $P_2 \subseteq B$, and there exist two disjoint directed paths $Q_{P_1,P_2}$ and $Q_{P_2,P_1}$ where each of them intersects $V(P_1)$ and $V(P_2)$ and is internally disjoint from $V(P_1) \cup V(P_2)$ such that 
	\begin{itemize}
		\item if $t$ is odd, then 
			\begin{itemize}
				\item $V(Q_{P_1,P_2}) \cap V(P_1) = \{m_{P_1}\}$ and $V(Q_{P_1,P_2}) \cap V(P_2) \neq \{m_{P_2}\}$, and
				\item $V(Q_{P_2,P_1}) \cap V(P_1) \neq \{m_{P_1}\}$ and $V(Q_{P_2,P_1}) \cap V(P_2) = \{m_{P_2}\}$;
			\end{itemize}
		\item if $t$ is even, then 
			\begin{itemize}
				\item $Q_{P_1,P_2}$ is a directed path between $m_{P_1}$ and $m'_{P_2}$, and 
				\item $Q_{P_2,P_1}$ is a directed path between $m_{P_2}$ and $m'_{P_1}$. 
			\end{itemize}
	\end{itemize}

\medskip

\noindent{\bf Proof of Claim 2:}
By Claim 1, there do not exist three disjoint threads in $D$ between $V(P_1)$ and $V(P_2)$.
So there exist a separation $(A,B)$ of $D$ of order at most two such that $P_1 \subseteq A$ and $P_2 \subseteq B$.
Since the underlying graph of $D$ is 2-connected, the order of $(A,B)$ equals two, and there exist two disjoint threads $Q_{P_1,P_2}$ and $Q_{P_2,P_1}$ between $V(P_1)$ and $V(P_2)$ internally disjoint from $V(P_1) \cup V(P_2)$.

We first assume that $t$ is odd.
By Claim 1, each $Q_{P_1,P_2}$ and $Q_{P_2,P_1}$ intersects $\{m_{P_1},m_{P_2}\}$.
Since $Q_{P_1,P_2}$ and $Q_{P_2,P_1}$ are disjoint, this claim holds,

Now we assume that $t$ is even.
By Claim 1, each $Q_{P_1,P_2}$ and $Q_{P_2,P_1}$ is a directed path. 
By Claim 1 and symmetry, we may assume that for $i \in [2]$, $Q_{P_i,P_{3-i}}$ is a directed path between $m_{P_i}$ and $m'_{P_{3-i}}$.
So the claim holds.
$\Box$

\medskip

\noindent{\bf Claim 3:} If $P_1,P_2,P_3,P_4,P_5$ are five disjoint $t$-alternating paths in $D$, then there exist $i \in [5]$ and a series-parallel 2-separation $(A_i,B_i)$ of $D$ with $P_i \subseteq A_i$.

\noindent{\bf Proof of Claim 3:}
Suppose to the contrary that for every $i \in [5]$, there exists no series-parallel 2-separation $(A_i,B_i)$ of $D$ with $P_i \subseteq A_i$.

By Claim 2, there exists a separation $(A,B)$ of $D$ of order two such that $P_1 \subseteq A$ and $P_2 \subseteq B$.
By assumption, $(A,B)$ and $(B,A)$ are not series-parallel 2-separations of $D$.
Since at most two of $P_3,P_4,P_5$ intersect $V(A \cap B)$, by symmetry, we may assume that $V(P_3) \cap V(B) = \emptyset$ and $P_3 \subseteq A- V(B)$.
Let $s$ and $t$ be the vertices in $V(A) \cap V(B)$.
Let $Q_{P_1,P_3}$ and $Q_{P_3,P_1}$ be the disjoint paths mentioned in Claim 2 (by taking $P_1,P_2$ in the statement of Claim 2 by $P_1,P_3$, respectively).

Suppose that $Q_{P_3,P_1}$ intersects both $s$ and $t$.
Then replacing the sub-thread of $Q_{P_3,P_1}$ between $s$ and $t$ by any thread in $B$ between $s$ and $t$, we obtain a thread in $D$ from $V(P_1-m_{P_1})$ to $m_{P_3}$ internally disjoint from $V(P_1) \cup V(P_3)$, so it must be a directed path by Claim 1.
But $(B,A)$ is not a series-parallel 2-separation of $D$, and the underlying graph of $D$ is 2-connected.
So $(B,s,t)$ is not a one-way series-parallel triple.
Hence there exists a thread in $B$ between $s$ and $t$ such that replacing the sub-thread of $Q_{P_3,P_1}$ between $s$ and $t$ by it does not create a directed path, a contradiction.

So $\lvert V(Q_{P_3,P_1}) \cap \{s,t\} \rvert \leq 1$.
Hence $Q_{P_3,P_1} \subseteq A$.
In addition, by Claim 1, there exists no thread in $D-\{m_{P_1},m_{P_2}\}$ between $V(P_1-m_{P_1})$ and $V(P_2-m_{P_2})$.

Suppose there exists a thread $P$ in $D-\{m_{P_1},m_{P_2}\}$ between $V(P_3)$ and $V(P_2-m_{P_2})$ internally disjoint from $V(P_2) \cup V(P_3)$.
By Claim 1, the end of $P$ in $V(P_3)$ is $m_{P_3}$.
Since there exists no thread in $D-\{m_{P_1},m_{P_2}\}$ between $V(P_1-m_{P_1})$ and $V(P_2-m_{P_2})$, $P$ is disjoint from $V(P_1)$.
So $V(Q_{P_3,P_1} \cup P) \cap V(P_1)$ consists of the end of $Q_{P_3,P_1}$ in $V(P_1)$.
Since $m_{P_3}$ is a common vertex in $Q_{P_3,P_1}$ and $P$, there exists a thread $P'$ in $Q_{P_3,P_1} \cup P$ from $V(P_1-m_{P_1})$ to $V(P_2-m_{P_2})$ internally disjoint from $V(P_1)$.
Since $Q_{P_3,P_1} \subseteq A$, $P'$ is internally disjoint from $V(P_2)$.
Hence $P'$ is a thread in $D$ from $V(P_1-m_{P_1})$ to $V(P_2-m_{P_2})$ internally disjoint from $V(P_1) \cup V(P_2)$, contradicting Claim 1.

Hence there exists no thread in $D-\{m_{P_1},m_{P_2}\}$ between $V(P_3)$ and $V(P_2-m_{P_2})$ internally disjoint from $V(P_2) \cup V(P_3)$.
So no component of $D-\{m_{P_1},m_{P_2}\}$ intersects both $V(P_3)$ and $V(P_2-m_{P_2})$.
Since there exists no thread in $D-\{m_{P_1},m_{P_2}\}$ between $V(P_1-m_{P_1})$ and $V(P_2-m_{P_2})$, no component of $D-\{m_{P_1},m_{P_2}\}$ intersects both $V(P_1-m_{P_1})$ and $V(P_2-m_{P_2})$.
Hence there exists a separation $(A',B')$ of $D$ of order two such that $V(A' \cap B') = \{m_{P_1},m_{P_2}\}$, $P_1 \cup P_3 \subseteq A'$ and $P_2 \subseteq B'$.
Since $P_2 \subseteq B'$, $(B',m_{P_1},m_{P_2})$ is not a one-way series-parallel triple.

Let $Q_{P_2,P_3}$ and $Q_{P_3,P_2}$ be the disjoint paths mentioned in Claim 2 (by taking $P_1,P_2$ in the statement of Claim 2 by $P_2,P_3$, respectively).

Suppose $V(Q_{P_2,P_3}) \cap V(P_1) = \emptyset$.
Since $m_{P_1} \in V(A' \cap B') \cap V(P_1)$, $\lvert V(Q_{P_2,P_3}) \cap V(A' \cap B') \rvert \leq 1$.
Since $V(Q_{P_2,P_3}) \cap V(A') \supseteq V(Q_{P_2,P_3}) \cap V(P_3) \neq \emptyset$, $Q_{P_2,P_3} \subseteq A'$.
Since $(B',m_{P_1},m_{P_2})$ is not a one-way series-parallel triple, we can concatenate $Q_{P_2,P_3}$ with a thread in $B'$ between $m_{P_2}$ and $m_{P_1}$ to create a non-directed path thread in $D$ from $V(P_3-m_{P_3})$ to $m_{P_1}$ internally disjoint from $V(P_1) \cup V(P_3)$, contradicting Claim 1.

So $V(Q_{P_2,P_3}) \cap V(P_1) \neq \emptyset$.
Let $P''$ be the sub-thread of $Q_{P_2,P_3}$ between $V(P_3-m_{P_3})$ and $V(P_1)$ internally disjoint from $V(P_3) \cup V(P_1)$.
By Claim 1, the end of $P''$ in $V(P_1)$ is $m_{P_1}$.
Since $m_{P_2}$ is an end of $Q_{P_2,P_3}$ not in $V(P_1) \cup V(P_3)$, $m_{P_2} \not \in V(P'')$.
So $\lvert V(P'') \cap V(A' \cap B') \rvert \leq 1$.
Since $V(P'') \cap V(A') \supseteq V(P'') \cap V(P_1) \neq \emptyset$, $P'' \subseteq A'$.
Since $Q_{P_1,P_2}$ is between $m_{P_1} \in V(A' \cap B')$ and $V(P_2-m_{P_2}) \subseteq V(B')$ internally disjoint from $V(P_2)$, $Q_{P_1,P_2} \subseteq B'$.
Hence $P'' \cup Q_{P_1,P_2}$ is a thread in $D$ from $V(P_3-m_{P_3})$ to $V(P_2-m_{P_2})$ internally disjoint from $V(P_3) \cup V(P_2)$, contradicting Claim 1.
This proves the claim.
$\Box$

\medskip

Let $\F = \{V(P): P$ is a $t$-alternating path in $D$ such that there exists no series-parallel 2-separation $(A,B)$ of $D$ with $P \subseteq A\}$.
By Claim 3, $\F$ does not contain five pairwise disjoint members.
Recall that the tree-width of the underlying graph of $D$ is at most $w$.
By Lemma \ref{tw ep}, there exists $Z \subseteq V(D)$ with $\lvert Z \rvert \leq 4(w+1)=f(t)$ such that $Z \cap S \neq \emptyset$ for every $S \in \F$.
This proves the lemma.
\end{pf}

\bigskip

A series-parallel 2-separation $(A,B)$ of a digraph $D$ is {\it maximal} if there exists no series-parallel 2-separation $(A',B')$ of $D$ with $A \subset A'$.

\begin{lemma} \label{maximal_sp_sep_0}
Let $D$ be a digraph whose underlying graph is 2-connected.
Assume that there do not exist distinct vertices $s,t$ and one-way series-parallel triples $(X,s,t)$ and $(Y,t,s)$ such that $D=X \cup Y$.
If $(A_1,B_1)$ and $(A_2,B_2)$ are distinct maximal series-parallel 2-separations of $D$, then for every $i \in [2]$, $V(A_i \cap B_i) \subseteq V(A_{3-i})$ or $V(A_i \cap B_i) \subseteq V(B_{3-i})$.
\end{lemma}

\begin{pf}
By symmetry, it suffices to show that $V(A_2 \cap B_2) \subseteq V(A_1)$ or $V(A_2 \cap B_2) \subseteq V(B_1)$.

For each $i \in [2]$, let $s_i$ and $t_i$ be the vertices in $A_i \cap B_i$.
Since $(A_1,s_1,t_1)$ is one-way, by symmetry, we may assume that every thread in $A_1$ between $s_1$ and $t_1$ is a directed path from $s_1$ to $t_1$.

Suppose to the contrary that $V(A_2 \cap B_2)-V(A_1) \neq \emptyset$ and $V(A_2 \cap B_2)-V(B_1) \neq \emptyset$.
By symmetry, we may assume $t_2 \in V(A_1)-V(B_1)$ and $s_2 \in V(B_1)-V(A_1)$.
In particular, $\{s_1,t_1\} \cap \{s_2,t_2\}=\emptyset$.

Let $P_2$ and $P_2'$ be threads in $A_2$ and $B_2$ between $s_2$ and $t_2$, respectively.
Since $s_2 \in V(B_1)-V(A_1)$ and $t_2 \in V(A_1)-V(B_1)$, $\lvert (V(P_2)-\{s_2,t_2\}) \cap \{s_1,t_1\} \rvert=1=\lvert (V(P_2')-\{s_2,t_2\}) \cap \{s_1,t_1\} \rvert$.
So $\{s_1,t_1\} \cap V(A_2)-V(B_2) \neq \emptyset \neq \{s_1,t_1\} \cap V(B_2)-V(A_2)$.
Hence for every thread $P$ in $A_1$ from $s_1$ to $t_1$, $P$ is a thread in $A_1$ between $V(A_2)-V(B_2)$ and $V(B_2)-V(A_2)$, so $P$ contains $V(A_2) \cap V(B_2) \cap V(A_1)=\{t_2\}$, and hence $P$ passes through $s_1,t_2,t_1$ in the order listed.
Since every thread in $A_1$ from $s_1$ to $t_1$ is a directed path from $s_1$ to $t_1$ and contains $t_2$, every thread in $A_1$ from $t_2$ to $\{s_1,t_1\}$ is either a directed path from $s_1$ to $t_2$ or a directed path from $t_2$ to $t_1$.

Since $(A_2,B_2)$ is a series-parallel 2-separation, $P_2$ is a directed path.
Since $P_2$ contains a subpath between $\{s_1,t_1\}$ and $t_2$, $P_2$ contains a vertex $v \in \{s_1,t_1\}$, where $v=s_1$ if $P_2$ is a directed path from $s_2$ to $t_2$, and $v=t_1$ otherwise.

Since $(A_2,s_2,t_2)$ is one-way and $P_2$ is an arbitrary thread in $A_2$ between $s_2$ and $t_2$, we know that every thread in $A_2$ from $s_2$ to $t_2$ contains $v$.
In particular, $v \in V(A_2)-V(B_2)$.
Let $u$ be the vertex in $\{s_1,t_1\}-\{v\}$.
Since $\{s_1,t_1\} \cap V(A_2)-V(B_2) \neq \emptyset \neq \{s_1,t_1\} \cap V(B_2)-V(A_2)$, we know $u \in V(B_2)-V(A_2)$.

Note that $V(A_1 \cup A_2) \cap V(B_1 \cap B_2) \subseteq \{s_1,t_1,s_2,t_2\}$.
Since $v \in V(A_2)-V(B_2)$ and $t_2 \in V(A_1)-V(B_1)$, $V(A_1 \cup A_2) \cap V(B_1 \cap B_2) \subseteq \{u,s_2\}$.

Let $Q$ be a thread in $A_1 \cup A_2$ from $s_2$ to $u$.
Since $s_2 \in V(B_1)-V(A_1)$ and $Q \subseteq A_1 \cup A_2$, the edge of $Q$ incident with $s_2$ is in $A_2-V(A_1)$.
Since $u \in V(B_2)-V(A_2)$ and $Q \subseteq A_1 \cup A_2$, the edge of $Q$ incident with $u$ is in $A_1-V(A_2)$.
So some internal vertex of $Q$ is in $V(A_2) \cap \{s_1,t_1\}=\{v\}$ and some internal vertex of $Q$ is in $V(A_1) \cap \{s_2,t_2\}=\{t_2\}$.
Hence $Q$ passes through $s_2,v,t_2,u$ in the order listed.
Since the subthread of $Q$ between $s_2$ and $t_2$ is in $A_2$, it is a directed path. 
Similarly, the subthread of $Q$ between $s_1$ and $t_1$ is a directed path. 
Since these two subtreads of $Q$ share a thread in $A_1$ between $v$ and $t_2$, and every thread in $A_1$ between $v$ and $t_2$ is a directed path, $Q$ is a directed path between $s_2$ and $u$.
Moreover, since all threads in $A_1$ between $v$ and $t_2$ are directed paths with the same direction, the direction of $Q$ is independent with the choice of $Q$.

Therefore, $(A_1 \cup A_2,s_2,u)$ is a one-way series-parallel triple.
So $(A_1 \cup A_2,B_1 \cap B_2)$ is a series-parallel 2-separation of $D$.
But $s_2 \in V(A_1 \cup A_2)-V(A_1)$, so $A_1 \subset A_1 \cup A_2$, contradicting the maximality of $A_1$.
\end{pf}

\begin{lemma} \label{maximal_sp_sep}
Let $D$ be a digraph whose underlying graph is 2-connected.
Assume that there do not exist distinct vertices $s,t$ and one-way series-parallel triples $(X,s,t)$ and $(Y,t,s)$ such that $D=X \cup Y$.
If $(A_1,B_1)$ and $(A_2,B_2)$ are distinct maximal series-parallel 2-separations of $D$, then $A_1 \subseteq B_2$ and $A_2 \subseteq B_1$.
\end{lemma}

\begin{pf}
By symmetry, it suffices to show $A_1 \subseteq B_2$.

For each $i \in [2]$, let $s_i$ and $t_i$ be the vertices in $A_i \cap B_i$.
By symmetry, we may assume that for every $i \in [2]$, every thread in $A_i$ between $s_i$ and $t_i$ is a directed path from $s_i$ to $t_i$.

\medskip

\noindent{\bf Claim 1:} $V(A_2) \cap V(B_2) \subseteq V(B_1)$.

\noindent{\bf Proof of Claim 1:}
Suppose to the contrary that $V(A_2) \cap V(B_2) - V(B_1) \neq \emptyset$.
By Lemma \ref{maximal_sp_sep_0}, $V(A_2) \cap V(B_2) \subseteq V(A_1)$.
Since the underlying graph of $D$ is 2-connected, the underlying graph of $B_1$ is connected.
Since $V(A_2) \cap V(B_2) \subseteq V(A_1)$, either $B_1 \subseteq A_2$ or $B_1 \subseteq B_2$.

Suppose that $B_1 \subseteq B_2$.
Then $A_1 \supseteq A_2$.
Since $(A_1,B_1) \neq (A_2,B_2)$, $A_1 \neq A_2$, so $A_1 \supset A_2$, contradicting the maximality of $A_2$.

So $B_1 \subseteq A_2$. 
Since the underlying graph of $D$ is 2-connected, there exist two disjoint threads $P_1,P_2$ in $A_1$ from $\{s_1,t_1\}$ to $\{s_2,t_2\}$.
By symmetry, we may assume that $P_1$ contains $s_2$, and $P_2$ contains $t_2$.
Since each $P_i$ intersects $\{s_2,t_2\}$ in at most one vertex, it is contained in $A_2$ or in $B_2$.
Since each $P_i$ intersects $\{s_1,t_1\} \subseteq V(B_1) \subseteq V(A_2)$, it is contained in $A_2$.
Since $B_1 \subseteq A_2$, $P_1 \cup P_2 \cup B_1 \subseteq A_2$.

For every thread $P$ in $B_1$ between $s_1$ and $t_1$, $P_1 \cup P \cup P_2$ is a thread in $A_2$ between $s_2$ and $t_2$, so $P_1 \cup P \cup P_2$ is a directed path from $s_2$ to $t_2$, and hence $P$ is a directed path in $B_1$ between $s_1$ and $t_1$ such that $P$ is a directed path from $s_1$ to $t_1$ if and only if $P_1$ contains $s_1$.
Therefore, $(B_1,t_1,s_1)$ is a one-way series-parallel triple.
But $(A_1,s_1,t_1)$ is a one-way series-parallel triple such that $A_1 \cup B_1=D$, a contradiction.
$\Box$

\medskip

Since the underlying graph of $D$ is 2-connected, the underlying graph of $A_1$ is connected.
Since $V(A_2) \cap V(B_2) \subseteq V(B_1)$ by Claim 1, either $A_1 \subseteq A_2$ or $A_1 \subseteq B_2$.
If $A_1 \subseteq A_2$, then since $(A_1,B_1) \neq (A_2,B_2)$, $A_1 \subset A_2$, a contradiction.
So $A_1 \subseteq B_2$.
This proves the lemma.
\end{pf}

\bigskip

A {\it series-parallel cover} of a digraph $D$ is a collection $\Se$ of separations of $D$ satisfying the following.
	\begin{itemize}
		\item For every $(A,B) \in \Se$, $(A,B)$ is a maximal series-parallel 2-separation of $D$. 
		\item If $(A_1,B_1)$ and $(A_2,B_2)$ are distinct members of $\Se$, then $A_1 \subseteq B_2$ and $A_2 \subseteq B_1$.
		\item For every series-parallel 2-separation $(A',B')$ of $D$, there exists $(A,B) \in \Se$ such that $A' \subseteq A$. 
	\end{itemize}

\begin{lemma} \label{collection sp sep}
Let $D$ be a digraph whose underlying graph is 2-connected.
If there do not exist distinct vertices $s,t$ and one-way series-parallel triples $(X,s,t)$ and $(Y,t,s)$ such that $D=X \cup Y$, then there exists a series-parallel cover of $D$. 
\end{lemma}

\begin{pf}
Let $\Se$ be the collection of all maximal series-parallel 2-separations of $D$.
So $\Se$ satisfies Statement 1 in the definition of a series-parallel cover.
By Lemma \ref{maximal_sp_sep} , $\Se$ satisfies Statement 2 in the definition of a series-parallel cover.
For every series-parallel 2-separation $(A',B')$ of $D$, there exists a maximal series-parallel 2-separation $(A,B)$ of $D$ such that $A' \subseteq A$, so $A' \subseteq A$ for some $(A,B) \in \Se$.
Hence $\Se$ is a series-parallel cover. 
\end{pf}

\section{Longer alternating paths} \label{sec:long_alt}

In this section we prove Theorem \ref{main_label} (that is, digraphs with no $k$-alternating paths are well-quasi-ordered).
The strategy is to use the lemmas proved in Section \ref{sec:sp_sep} to reduce the ``complexity'' of the digraphs, and prove well-quasi-ordering by induction on the ``complexity''.
It could be helpful if the readers read the related definitions and statements of the lemmas before going into their proofs to get a big picture of the entire process.

\begin{lemma} \label{3 cut-vertices}
Let $D$ be a digraph whose underlying graph is 2-connected.
Let $r,x,y$ be three distinct vertices.
Then either there exists a 1-alternating path between $r$ and $\{x,y\}$, or there exist $z \in \{x,y\}$ such that there exist a directed path from $r$ to $z$ and a directed path from $z$ to $r$.
\end{lemma}

\begin{pf}
Since the underlying graph is 2-connected, there exist two threads $P,Q$, where $P$ is between $r$ and $x$, and $Q$ is between $r$ and $y$, and their intersection is $\{r\}$.
By the 2-connectedness, there exists a thread $R$ in $D-\{r\}$ from $V(P)-\{r\}$ to $V(Q)-\{r\}$.

We are done if $P$ or $Q$ is not a directed path.
So we may assume that $P$ and $Q$ are directed paths.
Similarly, we may assume that $R$ is a directed path, for otherwise we are done.

We first assume that $r$ is the sink of both $P$ and $Q$ or the source of both $P$ and $Q$.
By symmetry, we may assume that $r$ is the sink of $P$ and $Q$.
By symmetry, we may assume that the sink of $R$ is in $V(P)-\{r\}$.
Then there exists a 1-alternating path from $r$ to $x$.

So we may assume that $r$ is the sink of one of $P$ and $Q$, and is the source of the other.
By symmetry, we may assume that $r$ is the sink of $Q$ and the source of $P$.
If $R$ has source in $V(Q)-\{r\}$, then there exists a 1-alternating path in $D$ from $r$ to $x$.
So we may assume that $R$ has source in $V(P)-\{r\}$.
If $V(Q) \cap V(R) \neq \{y\}$, then there exists a 1-alternating path between $r$ to $y$; otherwise there exists a directed path from $r$ to $y$.
Note that $Q$ is a directed path from $y$ to $r$, so we are done.
\end{pf}

\bigskip

Let $t,k$ be nonnegative integers.
We define $\F_{t,k}$ to be the set consisting of the rooted digraphs $(D,r)$ satisfying the following.
	\begin{itemize}
		\item The underlying graph of $D$ is connected.
		\item $r$ is not a cut-vertex of the underlying graph of $D$.
		\item There exists no $(t+1)$-alternating path in $D$.
		\item No block of $D$ contains a $t$-alternating path.
		\item There exists no $k$-alternating path in $D$ having $r$ as an end.
	\end{itemize}
Note that every one-vertex thread is a $0$-alternating path, so $\F_{t,0}=\emptyset$.
We define the following.
	\begin{itemize}
		\item Define $\F_t$ to be the set consisting of the rooted digraphs $(D,r)$ such that the underlying graph of $D$ is connected, and there exist no $t$-alternating path in $D$.
		\item Define $\F_t'$ to be the set consisting of the rooted digraphs $(D,r)$ such that there exist no $t$-alternating path in $D$, and either 
			\begin{itemize}
				\item the underlying graph of $D$ is 2-connected, or 
				\item the underlying graph of $D$ is connected and contains at most two vertices.
			\end{itemize}
		\item Define $\F_t^*$ to be the set consisting of the rooted digraphs $(D,r)$ such that there exist no $t$-alternating path in $D$.
	\end{itemize}
Note that $\F_t' \subseteq \F_t \subseteq \F_t^*$ and $\emptyset = \F_{t,0} \subseteq \F_{t,1} \subseteq ... \subseteq \F_{t,t+1} = \bigcup_{k \geq 0}\F_{t,k}$.

\begin{lemma} \label{root conn suff}
Let $t$ be a nonnegative integer.
If $\F_{t,t+1}$ is well-behaved, then $\F_t \cup \F_t^*$ is well-behaved.
\end{lemma}

\begin{pf}
If $\F_t$ is well-behaved, then $\F^*_t$ is well-behaved by Higman's lemma.
So it suffices to prove that $\F_t$ is well-behaved.

Let $(Q,\leq_Q)$ be a well-quasi-order.
For $i \in {\mathbb N}$, let $(D_i,r_i) \in \F_t$ and $\phi_i: V(D_i) \rightarrow Q$.
For each component $C$ of $D_i-r_i$, let $(D_{i,C},r_i)$ be the rooted digraph, where $D_{i,C}$ is the sub-digraph of $D_i$ induced by $V(C) \cup \{r_i\}$.
Hence $r_i$ is not a cut-vertex of the underlying graph of $D_{i,C}$ for all $i,C$.
Since $D_i$ has no $t$-alternating path and $D_{i,C}$ is connected, $(D_{i,C},r_i) \in \F_{t,t+1}$. 

For each $i \in {\mathbb N}$, let $a_i$ be the sequence $((D_{i,C},r_i): C$ is a component of $D_i-r_i)$.
Since $\F_{t,t+1}$ is well-behaved, by Higman's lemma, there exist $1 \leq j <j'$, a function $\iota$ that maps components of $D_j-r_j$ to components of $D_{j'}-r_{j'}$ injectively, such that for each component $C$ of $D_j-r_j$, there exists a strong immersion embedding $\eta_C$ from $(D_{j,C},r_j)$ to $(D_{j',\iota(C)},r_{j'})$ such that $\eta_C(r_j)=r_{j'}$ and $\phi_j(v) \leq_Q \phi_{j'}(\eta_C(v))$ for every $v \in V(D_{j,C})$.
Then it is easy to construct a strong immersion embedding $\eta$ from $(D_j,r_j)$ and $(D_{j'},r_{j'})$ such that $\eta(r_j)=r_{j'}$ and $\phi_j(v) \leq_Q \phi_{j'}(\eta(v))$ for every $v \in V(D_j)$.
This proves the lemma.
\end{pf}

\bigskip

Let $(D,r)$ be a rooted digraph, and let $v$ be a cut-vertex of the underlying graph of $D$.
Assume that $v \neq r$.
So there exists a separation $(A,B)$ such that $V(A \cap B)=\{v\}$, $r \in V(B)-V(A)$, and $V(A)-V(B) \neq \emptyset$.
For each such $(A,B)$, if $v$ is not a cut-vertex of $A$, then we call the rooted digraph $(A,v)$ a {\it branch of $D$ at $v$}.

\begin{lemma} \label{2-conn is suff}
Let $t$ be a nonnegative integer.
If $\F'_t$ is well-behaved, then $\F_{t,k}$ is well-behaved for every nonnegative integer $k$.
\end{lemma}

\begin{pf}
We shall prove this lemma by induction on $k$.
Since $\F_{t,0}=\emptyset$, the lemma holds when $k=0$.
So we may assume that $k \geq 1$ and $\F_{t,k-1}$ is well-behaved.

Let $(Q,\leq_Q)$ be a well-quasi-order.
For $i \in {\mathbb N}$, let $(D_i,r_i) \in \F_{t,k}$, and let $\phi_i:V(D_i) \rightarrow Q$.

For $i \in {\mathbb N}$, let $S_i$ be the subset of $V(D_i)$ satisfying the following.
	\begin{itemize}
		\item Every vertex in $S_i$ is a cut-vertex of the underlying graph of $D_i$.
		\item For each $v \in S_i$, some branch of $D_i$ at $v$ belongs to $\F_t' \cup \F_{t,k-1}$.
		\item If some branch $(R,v)$ of $D_i$ at $v$ belongs to $\F_t' \cup \F_{t,k-1}$, then there exists $u \in S_i$ (not necessarily different from $v$) such that $R \subseteq R'$ for some branch $(R',u)$ of $D_i$ at $u$ with $(R',u) \in \F_t' \cup \F_{t,k-1}$.
		\item There exist no distinct $u,v \in S_i$ such that $u \in V(R)$ for some branch $(R,v)$ of $D_i$ at $v$ with $(R,v) \in \F_t' \cup \F_{t,k-1}$.
	\end{itemize}
Note that $r_i \not \in S_i$ since $r_i$ is not a cut-vertex of the underlying graph of $D_i$.
For each $i \in {\mathbb N}$, let $\Se_i = \{((B,v),\phi_i|_{V(B)}): (B,v)$ is a branch at $v$ for some $v \in S_i$ and $(B,v) \in \F_t' \cup \F_{t,k-1}\}$, and let $\Se'_i = \{((\bigcup_{((L,v),\phi_i|_{V(L)}) \in \Se_i}L,v),\phi_i|_{V(L)}): v \in S_i\}$.

Let $Q_1 = \bigcup_{i \geq 1}\Se'_i$.
Let $\preceq_1$ be the binary relation on $Q_1$ such that for any $(B_1,f_1),(B_2,f_2) \in Q_1$, $(B_1,f_1) \preceq_1 (B_2,f_2)$ if and only if there exists a strong immersion embedding $\eta$ from $B_1$ to $B_2$ such that $f_1(v) \preceq_1 f_2(\eta(v))$ for every $v \in V(B_1)$.
Since $\F'_t$ and $\F_{t,k-1}$ are well-behaved, $(Q_1,\preceq_1)$ is a well-quasi-order by Higman's Lemma.
Let $\perp$ be an element not in $Q_1$.
Let $(Q_2,\preceq_2)$ be the well-quasi-order obtained by the disjoint union of $(Q_1,\preceq_1)$ and $(\{\perp\},=)$.
Let $(Q_3,\preceq_3)$ be the well-quasi-order obtained  by the Cartesian product of $(Q,\leq_Q)$ and $(Q_2,\preceq_2)$.

For each $i \geq 1$, 
	\begin{itemize}
		\item let $D_i' = D_i-\bigcup_{v \in S_i}\bigcup_{((R,v),\phi_i|_{V(R)}) \in \Se'_i}(V(R)-\{v\})$, and 
		\item define $\phi_i': V(D_i') \rightarrow Q_3$ to be the function such that for every $v \in V(D_i')$,
			\begin{itemize}
				\item if $v \in S_i$, then $\phi'_i(v) = (\phi_i(v), ((R,v),\phi_i|_{V(R)}))$, where $(R,v)$ is the unique member of $\Se'_i$ such that its second entry is $v$, and 
				\item if $v \in V(D_i')-S_i$, then $\phi'_i(v) = (\phi_i(v),\perp)$.
			\end{itemize}
	\end{itemize}
Note that the underlying graph of each $D_i'$ is connected.

To prove this lemma, it suffices to prove that there exist $1 \leq j <j'$ and a strong immersion embedding $\eta$ from $(D_j',r_j)$ to $(D_{j'}',r_{j'})$ such that $\phi'_j(v) \preceq_3 \phi'_{j'}(\eta(v))$ for every $v \in V(D_j')$.

If $i$ is an index such that the underlying graph of $D_i'$ is 2-connected or has at most two vertices, then the underlying graph of $D_i'$ is a block of the underlying graph of $D_i$, and since $(D_i,r_i) \in \F_{t,k}$, $(D_i',r_i) \in \F_t'$.
Since $\F_t'$ is well-behaved, if there are infinitely many indices $i$ satisfying the previous property, then we are done.
So by removing finitely many terms in the sequence, we may assume that for each $i \in {\mathbb N}$, the underlying graph of each $D_i'$ is not 2-connected and has at least three vertices. 

\medskip

\noindent{\bf Claim 1:} For every $i \in {\mathbb N}$ and every cut-vertex $x$ of the underlying graph of $D_i'$, every thread in $D_i'$ from $r_i$ to $x$ is a directed path.

\medskip

\noindent{\bf Proof of Claim 1:}
Suppose to the contrary that there exists a thread $P$ in $D_i'$ from $r_i$ to $x$ such that $P$ is not a directed path.
Let $H'$ be a branch of $D_i'$ at $x$.
Since $x$ is a cut-vertex of the underlying graph of $D_i'$, $x$ is a cut-vertex of the underlying graph of $D_i$, so there exists a branch $H$ of $D_i$ at $x$ containing $H'$.
By the definition of $S_i$, $(H,x) \not \in \F_{t,k-1}$.
Hence there exists a $(k-1)$-alternating path $P'$ in $H$ having $x$ as an end.
So $P \cup P'$ is a $k$-alternating path in $D_i$ having $r_i$ as an end, contradicting that $(D_i,r_i) \in \F_{t,k}$.
$\Box$

\medskip

\noindent{\bf Claim 2:} For every $i \in {\mathbb N}$ and every cut-vertex $x$ of the underlying graph of $D_i'$, either all directed paths in $D_i'$ between $r_i$ and $x$ are from $r_i$ to $x$, or all directed paths in $D_i'$ between $r_i$ and $x$ are from $x$ to $r_i$.

\medskip

\noindent{\bf Proof of Claim 2:}
Suppose to the contrary that there exist a directed path $P_1$ from $r_i$ to $x$ and a directed path $P_2$ from $x$ to $r_i$.
Let $H'$ be a branch of $D_i'$ at $x$.
So there exists a branch $H$ of $D_i$ at $x$ such that $H$ contains $H'$.
By the definition of $S_i$, $(H,x) \not \in \F_{t,k-1}$.
Hence there exists a $(k-1)$-alternating path $P$ in $H$ having $x$ as an end.
If $k \geq 2$, then $P$ contains at least one edge, so $P_1 \cup P$ or $P_2 \cup P$ is a $k$-alternating path in $D_i$ having $r_i$ as an end, contradicting that $(D_i,r_i) \in \F_{t,k}$.
Hence $k=1$.
Since the underlying graph of $D_i$ is connected and $(H,x)$ is a branch at $x$, there exists an one-edge directed path $P'$ of $H$ having $x$ as an end.
Then $P_1 \cup P'$ or $P_2 \cup P'$ is a $k$-alternating path in $D_i$ having $r_i$ as an end, contradicting that $(D_i,r_i) \in \F_{t,k}$.
This proves the claim.
$\Box$

\medskip

\noindent{\bf Claim 3:} Every block of the underlying graph of $D_i'$ contains at most two cut-vertices of the underlying graph of $D_i'$, and the block of the underlying graph of $D_i'$ containing $r_i$ contains at most one cut-vertex of the underlying graph of $D_i'$.

\medskip

\noindent{\bf Proof of Claim 3:}
Suppose to the contrary that there exists a block $B$ of the underlying graph of $D_i'$ such that either $B$ contains $r_i$ and two cut-vertices $x,y$ of the underlying graph of $D_i'$, or $B$ contains three cut-vertices $r,x,y$, where $r$ is the cut-vertex contained in the parent block of $B$.
Since $r_i$ is not a cut-vertex of the underlying graph of $D_i$, $B$ contains at least three vertices, so $B$ is 2-connected.
If $B$ contains $r_i$, then we let $r=r_i$.
By Lemma \ref{3 cut-vertices}, either there exists a 1-alternating path from $r$ to $x$ or $y$, or there exists $z \in \{x,y\}$ such that there exist a directed path from $r$ to $z$ and a directed path from $z$ to $r$.
This contradicts Claims 1 and 2.
$\Box$

\medskip

For each $i \geq 1$, by Claims 1-3, either all threads in $D_i'$ between $r_i$ and a cut-vertex of the underlying graph of $D_i'$ are directed paths with source $r_i$, or all threads in $D_i'$ between $r_i$ and a cut-vertex of the underlying graph of $D_i'$ are directed paths with sink $r_i$.
Hence by possibly reversing the direction of all edges of $D_i$ and removing some $D_i$ from the sequence, we may assume that for each $i \in {\mathbb N}$ and for every cut-vertex $x$ of the underlying graph of $D_i'$, every thread in $D_i'$ between $r_i$ and $x$ is a directed path from $r_i$ to $x$.

\medskip

\noindent{\bf Claim 4:} For every $i \in {\mathbb N}$, $(D_i',r_i)$ is a $\F_t'$-series-parallel tree.

\medskip

\noindent{\bf Proof of Claim 4:}
Since $(D_i,r_i) \in \F_{t,k}$, no block of $D_i$ contains a $t$-alternating path.
So no block of $D_i'$ contains a $t$-alternating path.
Hence for every block $B$ of the underlying graph of $D_i'$, $(B,v) \in \F'_t$, where $v=r_i$ if $r_i \in V(B)$, and $v$ is the cut-vertex of the underlying graph of $D_i'$ contained in $B$ and the parent block of $B$.
Since the underlying graph of $D_i$ is connected, Claims 1-3 imply that $(D_i',r_i)$ is a $\F_t'$-series-parallel tree.
$\Box$

\medskip

Let $\F'$ be the set of one-way series-parallel triples $(B,x,y)$ such that $(B,x) \in \F_t'$ and $y \in V(B)-\{x\}$.
For every $(B,x,y) \in \F'$, since $(B,x) \in \F_t'$, there exist no $t$-alternating path in $B$, so $(B,x,y) \in \A_t$.
Hence $\F' \subseteq \A_t$.
Let $\F''$ be the set of all series-parallel triples that are truncations of members of $\F'$.
By Statement 1 of Lemma \ref{truncate no boundary alternating}, $\F'' \subseteq \A_t$.
By Lemma \ref{sp wqo 2}, $\F'$ and $\F''$ are well-behaved.
Since $\F_t'$ is well-behaved, this lemma follows from Claim 4 and Lemma \ref{series parallel tree wqo}.
\end{pf}

\begin{lemma} \label{no_cover}
Let $t$ be a positive integer.
Let $(Q,\leq_Q)$ be a well-quasi-order.
For $i \geq 1$, let $(D_i,r_i) \in \F_t'$ and $\phi_i: V(D_i) \rightarrow Q$.
If for each $i \geq 1$, $D_i$ has no series-parallel cover, then there exist $1 \leq j <j'$ and a strong immersion embedding $\eta$ from $(D_j,r_j)$ to $(D_{j'},r_{j'})$ such that $\phi_j(v) \leq_Q \phi_{j'}(\eta(v))$ for every $v \in V(D_j)$.
\end{lemma}

\begin{pf}
This lemma obviously holds if there are infinitely many indices $i$ such that $D_i$ contains at most two vertices.
So by removing finitely members in the sequence, we may assume that $D_i$ contains at least three vertices for each $i \geq 1$.
By the definition of $\F_t'$, the underlying graph of each $D_i$ is 2-connected.
For each $i \geq 1$, since $D_i$ has no series-parallel cover, by Lemma \ref{collection sp sep}, there exist one-way series-parallel triples $(X_i,s_i,t_i)$ and $(Y_i,t_i,s_i)$ for some distinct vertices $s_i,t_i \in V(D_i)$ such that $D_i = X_i \cup Y_i$.

For each $i \geq 1$, since $(D_i,r_i) \in \F_t'$, $D_i$ does not contain a $t$-alternating path, so $(X_i,s_i,t_i)$ and $(Y_i,t_i,s_i)$ belong to $\A_t$.
Let $\F=\{((X_i,s_i,t_i),\phi_i|_{V(X_i)}): i \geq 1\}$.
Let $\preceq_1$ be the simulation relation defined on $\F$.
By Lemma \ref{sp wqo 2}, $(\F,\preceq_1)$ is a well-quasi-order.
Let $(Q_2,\preceq_2)$ be the well-quasi-order obtained from $(Q,\leq_Q)$ and $(\F,\preceq_1)$ by taking Cartesian product.

For each $i \geq 1$, let $f_i: V(Y_i) \rightarrow Q_2$ be the function such that $f_i(v)=(\phi_i(v), \linebreak ((X_i,s_i,t_i),\phi_i|_{V(X_i)}))$ for every $v \in V(Y_i)$.
Since $\A_t$ is well-behaved, there exist $1 \leq j < j'$ and a strong immersion embedding $\eta_Y$ from $(Y_j,t_j,s_j)$ to $(Y_{j'},t_{j'},s_{j'})$ such that $f_j(v) \preceq_2 f_{j'}(\eta_Y(v))$ for every $v \in V(Y_j)$.
Note that $\eta_Y(s_j)=s_{j'}$ and $\eta_Y(t_j)=t_{j'}$ by the definition of strong immersion embedding of general rooted digraphs.
Since $f_j(s_j) \preceq_2 f_{j'}(\eta_Y(s_j))=f_{j'}(s_{j'})$, there exists a strong immersion embedding $\eta_X$ from $(X_j,s_j,t_j)$ to $(X_{j'},s_{j'},t_{j'})$ such that $\phi_j(v) \leq_Q \phi_{j'}(\eta_X(v))$ for every $v \in V(X_j)$.
Then combining $\eta_Y$ and $\eta_X$ results in a strong immersion embedding $\eta$ from $(D_j,r_j)$ to $(D_{j'},r_{j'})$ such that $\phi_j(v) \leq_Q \phi_{j'}(\eta(v))$ for every $v \in V(D_j)$.
\end{pf}

\bigskip

Let $\Se$ be a series-parallel cover of a digraph $D$.
Let $t$ be a positive integer.
The {\it $(\Se,t)$-compression} of $D$ is the digraph obtained from $D$ by for each $(A,B) \in \Se$ with $\lvert V(A)-V(B) \rvert \geq \min\{t,2\}$, 
	\begin{itemize}
		\item deleting $V(A)-V(B)$,
		\item deleting all edges of $A$ between the two vertices in $V(A \cap B)$,
		\item adding new vertices $v_{A,L},v_{A,M},v_{A,R}$ and new edges such that $v_{A,0}v_{A,L}v_{A,M}v_{A,R}v_{A,1}$ is a directed path from $v_{A,0}$ to $v_{A,1}$, where $v_{A,0}$ and $v_{A,1}$ are the two vertices in $V(A \cap B)$ such that every thread in $A$ is a directed path from $v_{A,0}$ to $v_{A,1}$, and
		\item duplicating $v_{A,0}v_{A,L}$, $v_{A,L}v_{A,M}$, $v_{A,M}v_{A,R}$ and $v_{A,R}v_{A,1}$ such that the following hold. 
			\begin{itemize}
				\item The number of edges between $v_{A,0}$ and $v_{A,L}$ equals the degree of $v_{A,0}$ in $A$.
				\item The number of edges between $v_{A,L}$ and $v_{A,M}$ equals the maximum number of edge-disjoint directed paths in $A$ from $v_{A,0}$ to $v_{A,1}$.
				\item The number of edges between $v_{A,M}$ and $v_{A,R}$ equals the maximum number of edge-disjoint directed paths in $A$ from $v_{A,0}$ to $v_{A,1}$.
				\item The number of edges between $v_{A,R}$ and $v_{A,1}$ equals the degree of $v_{A,1}$ in $A$.
			\end{itemize}
	\end{itemize}

\begin{lemma} \label{compression_hit}
Let $t$ be a positive integer. 
Let $D$ be a digraph.
Let $Z$ be a subset of $V(D)$ such that for every $t$-alternating path $P$ in $D$, either $V(P) \cap Z \neq \emptyset$, or $P \subseteq A$ for some series-parallel 2-separation $(A,B)$ of $D$.
Let $\Se$ be a series-parallel cover of $D$.
Let $D'$ be the $(\Se,t)$-compression of $D$.
Then there exists $Z' \subseteq V(D) \cap V(D')$ with $\lvert Z' \rvert \leq 2\lvert Z \rvert$ such that $V(P') \cap Z' \neq \emptyset$ for every $t$-alternating path $P'$ in $D'$.
\end{lemma}

\begin{pf}
For each $(A,B) \in \Se$ with $\lvert V(A)-V(B) \rvert \geq \min\{t,2\}$, let $v_{A,0},v_{A,L},v_{A,M},v_{A,R},v_{A,1}$ be the vertices mentioned in the definition of the $(\Se,t)$-compression.
Let $Z' = \{v \in Z: v \not \in V(A)-V(B)$ for every $(A,B) \in \Se$ with $\lvert V(A)-V(B) \rvert \geq \min\{t,2\}\} \cup \{v_{A,0},v_{A,1}: (A,B) \in \Se$ with $\lvert V(A)-V(B) \rvert \geq \min\{t,2\}$ and $Z-V(B) \neq \emptyset\}$.
Note that $Z' \subseteq V(D') \cap V(D)$ and $\lvert Z' \rvert \leq 2\lvert Z \rvert$.

Suppose to the contrary that there exists a $t$-alternating path $P'$ in $D'$ with $V(P') \cap Z'=\emptyset$.
We may assume that $\lvert E(P') \rvert$ is as small as possible.

Since for every $(A,B) \in \Se$ with $\lvert V(A)-V(B) \rvert \geq \min\{t,2\}$, the edges incident with $\{v_{A,L}, \allowbreak v_{A,R}\}$ are obtained by copying edges in a directed path $v_{A,0}v_{A,L}v_{A,M}v_{A,R}v_{A,1}$, so the minimality of $\lvert E(P') \rvert$ implies that 
	\begin{itemize}
		\item if $v_{A,L} \in V(P')$, then $P'$ contains an edge between $v_{A,0}$ and $v_{A,L}$,
		\item if $v_{A,R} \in V(P')$, then $P'$ contains an edge between $v_{A,R}$ and $v_{A,1}$,
		\item if $v_{A,M} \in V(P')$, then $P'$ contains a directed path $v_{A,0}v_{A,L}v_{A,M}v_{A,R}v_{A,1}$, and 
		\item if $P'$ contains both $v_{A,L}$ and $v_{A,R}$ but not $v_{A,M}$, then $v_{A,L'}$ and $v_{A,R}$ are the ends of $P'$, and $t \geq 2$ (since $P'$ is not a directed path).
	\end{itemize}
In addition, for every $(A,B) \in \Se$ with $\lvert V(A)-V(B) \rvert \geq \min\{t,2\}$, there exist vertices $u_A$ and $u_A'$ in $V(A)-V(B)$ such that $u_A$ is a neighbor of $v_{A,0}$ in $D$ and $u_A'$ is a neighbor of $v_{A,1}$ in $D$, and such that $u_A$ and $u_A'$ are distinct when $\lvert V(A)-V(B) \rvert \geq 2$.

Let $P$ be the thread in $D$ obtained from $P'$ by for each $(A,B) \in \Se$ with $\lvert V(A)-V(B) \rvert \geq \min\{t,2\}$, 
	\begin{itemize}
		\item if $v_{A,M} \in V(P')$, then deleting $v_{A,L},v_{A,M},v_{A,R}$ and adding a directed path in $A$ from $v_{A,0}$ to $v_{A,1}$, 
		\item if $v_{A,L} \in V(P')$ and $v_{A,M} \not \in V(P')$, then deleting $v_{A,L}$ and adding an edge of $D$ between $u_A$ and $v_{A,0}$, and
		\item if $v_{A,R} \in V(P')$ and $v_{A,M} \not \in V(P')$, then deleting $v_{A,R}$ and adding an edge of $D$ between $u_A'$ and $v_{A,1}$.
	\end{itemize}
Since $P'$ is a $t$-alternating path, $P$ is a $t$-alternating path in $D$.

Note that there exists no $(A,B) \in \Se$ with $\lvert V(A)-V(B) \rvert \geq \min\{t,2\}$ and $P \subseteq A$, for otherwise $P'$ is contained a directed path $v_{A,0}v_{A,L}v_{A,M}v_{A,R}v_{A,1}$, contracting $t \geq 1$.
Moreover, for each $(A,B) \in \Se$, if $\lvert V(A)-V(B) \rvert \leq \min\{t-1,1\}$ and $P \subseteq A$, then since $(A,B)$ is a series-parallel 2-separation, we know $t=1$, so $V(A)-V(B)=\emptyset$ and hence $P$ is a directed edge, a contradiction.
So there exists no $(A,B) \in \Se$ such that $P \subseteq A$.
If there exists a series-parallel 2-separation $(A',B')$ of $D$ with $P \subseteq A'$, then since $\Se$ is a series-parallel cover, there exists $(A,B) \in \Se$ with $P \subseteq A' \subseteq A$, a contradiction.
So the property of $Z$ implies that $V(P) \cap Z \neq \emptyset$.

Since $V(P') \cap Z' = \emptyset$, by the construction of $Z'$, there exists $z \in Z \cap V(P)-V(B_z)$ for some $(A_z,B_z) \in \Se$ with $\lvert V(A_z)-V(B_z) \rvert \geq \min\{t,2\}$.
But this implies that $\{v_{A_z,0},v_{A_z,1}\} \subseteq Z'$ and $\{v_{A_z,0},v_{A_z,1}\} \cap V(P') \neq \emptyset$, so $V(P') \cap Z' \neq \emptyset$, a contradiction.
\end{pf}

\begin{lemma} \label{apex wqo}
Let $\F$ be a well-behaved set of rooted digraphs, and let $s$ be a positive integer.
Let $\F'$ be the set consisting of the rooted digraphs $(D,r)$ satisfying that $(D-X,r') \in \F$ for some $X \subseteq V(D)$ with $r \in X$ and $\lvert X \rvert \leq s$ and for some $r' \in V(D)-X$.
Then $\F'$ is well-behaved.
\end{lemma}

\begin{pf}
Let $(Q,\leq_Q)$ be a well-quasi-order.
For $i \geq 1$, let $(D_i,r_i) \in \F'$ and let $\phi_i: V(D_i) \rightarrow Q$.

By the definition of $\F'$, for each $i \geq 1$, there exist $X_i \subseteq V(D_i)$ with $r_i \in X_i$ and $\lvert X_i \rvert \leq s$ and $r_i' \in V(D_i)-X_i$ such that $(D_i-X_i,r_i') \in \F$.
For each $i \in {\mathbb N}$, we denote $X_i$ by $\{u_{i,1},u_{i,2},...,u_{i,\lvert X_i \rvert}\}$, where $u_{i,1}=r_i$.
Since $\lvert X_i \rvert \leq s$ for all $i$, we may assume that $\lvert X_1 \rvert=\lvert X_i \rvert$ for all $i \geq 1$.
By Higman's Lemma, we may assume that for all $1 \leq a < b$, $D_a[X_a]$ is a subdigraph of $D_b[X_b]$ such that for each $j \in [\lvert X_1 \rvert]$, $u_{a,j}$ corresponds to $u_{b,j}$, and $\phi_a(u_{a,j}) \leq_Q \phi_b(u_{b,j})$. 

Let $(Q_1,\preceq_1)$ be the well-quasi-order obtained by the disjoint union of $({\mathbb N},\leq)$ and $(\{0\},=)$.
Let $(Q_2,\preceq_2)$ be the well-quasi-order that is obtained by the Cartesian product of $2\lvert X_1 \rvert$ copies of $(Q_1,\preceq_1)$.
Let $(Q_3,\preceq_3)$ be the well-quasi-order obtained by the Cartesian product $(Q,\leq_Q)$ and $(Q_2,\preceq_2)$.

For each $i \geq 1$ and $v \in V(D_i)-X_i$, define $\phi_i'(v)=(\phi_i(v),a_1,b_1,a_2,b_2,...,a_{\lvert X_i \rvert},b_{\lvert X_i \rvert})$, where for each $j \in [\lvert X_1 \rvert]$, $a_j$ is the number of edges of $D_i$ from $u_{i,j}$ to $v$, and $b_j$ is the number of edges of $D_i$ from $v$ to $u_{i,j}$.
Note that each $\phi'_i$ is a function from $V(D_i')$ to $Q_3$.

Since $\F$ is well-behaved, there exist $1 \leq j < j'$ and a strong immersion embedding $\eta'$ from $(D_j-X_j,r_j')$ to $(D_{j'}-X_{j'},r_{j'}')$ such that $\phi'_j(v) \preceq_3 \phi'_{j'}(\eta'(v))$ for every $v \in V(D_j)-X_j$.
Then it is easy to extend $\eta'$ to a strong immersion embedding $\eta$ from $(D_j,r_j)$ to $(D_{j'},r_{j'})$ such that $\phi_j(v) \leq_{Q} \phi_{j'}(\eta(v))$ for all $v \in V(D_j)$, and $\eta(u_{j,\ell})=u_{j',\ell}$ for all $\ell \in [\lvert X_1 \rvert]$.
This proves the lemma.
\end{pf}

\begin{lemma} \label{2-conn wqo}
For every positive integer $t$, $\F_t'$ is well-behaved.
\end{lemma}

\begin{pf}
We shall prove this lemma by induction on $t$.
By Lemma \ref{1-alt wqo}, $\F_1'$ is well-behaved.
So we may assume that $t \geq 2$ and $\F_{t-1}'$ is well-behaved.
By Lemma \ref{2-conn is suff}, $\F_{t-1,t}$ is well-behaved.
By Lemma \ref{root conn suff}, $\F_{t-1}^*$ is well-behaved.

Let $(Q,\leq_Q)$ be a well-quasi-order.
For $i \geq 1$, let $(D_i,r_i) \in \F_t'$ and $\phi_i: V(D_i) \rightarrow Q$.
It suffices to prove that there exist $1 \leq j <j'$ and a strong immersion embedding $\eta$ from $(D_j,r_j)$ to $(D_{j'},r_{j'})$ such that $\phi_j(v) \leq_Q \phi_{j'}(\eta(v))$ for every $v \in V(D_j)$.

We are done if there are infinitely many indices $i$ such that either $\lvert V(D_i) \rvert \leq 2$ or $D_i$ has no series-parallel cover by Lemma \ref{no_cover}.
So by removing finitely many members from the sequence, we may assume that for each $i \in {\mathbb N}$, the underlying graph of $D_i$ is 2-connected, and $D_i$ has a series-parallel cover $\Se_i$.

By Lemma \ref{ep path}, there exists a positive integer $N$ such that for each $i \in {\mathbb N}$, there exists $Z_i \subseteq V(D_i)$ with $\lvert Z_i \rvert \leq N$ such that for every $(t-1)$-alternating path $P$ in $D_i$, either $V(P) \cap Z_i \neq \emptyset$, or there exists a series-parallel 2-separation $(A,B)$ of $D_i$ with $P \subseteq A$.

For $i \geq 1$, define the following:
	\begin{itemize}
		\item define $D_i'$ to be the $(\Se_i,t)$-compression of $D_i$, 
		\item for each $(A,B) \in \Se_i$ with $\lvert V(A)-V(B) \rvert \geq \min\{t,2\}$, let $v_{A,0},v_{A,L},v_{A,M},v_{A,R},v_{A,1}$ be the vertices mentioned in the definition of the $(\Se_i,t)$-compression of $D_i$, and 
		\item if $r_i \in V(D_i) \cap V(D'_i)$, then let $r_i'=r_i$; otherwise, let $r_i'$ be an arbitrary vertex of $D'_i$.
	\end{itemize}
By Lemma \ref{compression_hit}, for each $i \geq 1$, there exists $Z_i' \subseteq V(D_i) \cap V(D_i')$ with $r_i' \in Z_i'$ and $\lvert Z_i' \rvert \leq 2N+1$ such that every $(t-1)$-alternating in $D_i'$ intersects $Z_i'$.

For $i \geq 1$, let $D_i''=D_i'-Z_i'$ and $r_i''$ be a vertex of $D_i''$.
For each $i \geq 1$, since every $(t-1)$-alternating path in $D_i'$ intersects $Z_i'$, we know $(D_i'',r_i'') \in \F_{t-1}^*$.

Let $\F^*$ be the family of rooted digraphs $(D,r)$ such that there exists $r \in Z \subseteq V(D)$ with $\lvert Z \rvert \leq 2N+1$ and such that $(D-Z,r') \in \F_{t-1}^*$ for some $r' \in V(D)-Z$.
Note that $(D_i',r_i') \in \F^*$ for each $i \geq 1$.
Since $\F^*_{t-1}$ is well-behaved, by Lemma \ref{apex wqo}, $\F^*$ is well-behaved.

For each $(A,B) \in \Se_i$ with $\lvert V(A)-V(B) \rvert \geq \min\{t,2\}$, let $[X_A,Y_A]$ be a splitter of $(A,v_{A,0},v_{A,1})$.
Recall that $(A_{X_A},v_{A,0},v_{{A,1}_{Y_A}})$ and $(A_{Y_A},v_{{A,0}_{X_A}},v_{A,1})$ are the truncations of $(A,v_{A,0},v_{A,1})$ with respect to $[X_A,Y_A]$.

By possibly further adding a new element into $Q$ and further label $r_i$ by using this element, we may assume that for each $i \geq 1$, $\phi_i(r_i)$ is an element in $Q$ incomparable with all other elements in $Q$, and $\phi_i(v) \neq \phi_i(r_i)$ for every $v \in V(D_i)$. 

Let $Q_1 = \{((D,s,t),\phi): (D,s,t) \in \A_{t}, \phi: V(D) \rightarrow Q\}$.
Let $\preceq_1$ be the simulation relation on $Q_1$.
By Lemma \ref{sp wqo 2}, $(Q_1,\preceq_1)$ is a well-quasi-order.
Let $(Q_2,\preceq_2)$ be the well-quasi-order obtained from $(Q_1,\preceq_1)$ and $([3],=)$ by taking Cartesian product.
Let $(Q_3,\preceq_3)$ be the well-quasi-order obtained from $(Q,\leq_Q)$ and $(Q_2,\preceq_2)$ by taking disjoint union.

For $i \geq 1$, define $\phi'_i$ to be a function with domain $V(D_i')$ as follows.
	\begin{itemize}
		\item If $v \in V(D'_i)-\{v_{A,L},v_{A,M},v_{A,R}: (A,B) \in \Se_i$ with $\lvert V(A)-V(B) \rvert \geq \min\{t,2\}\}$, then define $\phi_i'(v)=\phi_i(v)$.
		\item If $v = v_{A,M}$ for some $(A,B) \in \Se_i$ with $\lvert V(A)-V(B) \rvert \geq \min\{t,2\}$, then define $\phi'_i(v)=(((A,v_{A,0},v_{A,1}),\phi_i|_{V(A)}),1)$.
		\item If $v=v_{A,L}$ for some $(A,B) \in \Se_i$ with $\lvert V(A)-V(B) \rvert \geq \min\{t,2\}$, then define $\phi'_i(v)=(((A_{X_A},v_{A,0},v_{{A,1}_{Y_A}}),\phi_i|_{X_A}),2)$.
		\item If $v=v_{A,R}$ for some $(A,B) \in \Se_i$ with $\lvert V(A)-V(B) \rvert \geq \min\{t,2\}$, then define $\phi'_i(v)=(((A_{Y_A},v_{{A,0}_{X_A}},v_{A,1}),\phi_i|_{Y_A}),3)$.
	\end{itemize}
By Lemma \ref{truncate no boundary alternating}, the image of $\phi_i'$ for each $i$ is contained in $Q_3$.

Since $\F^*$ is well-behaved, there exist $1 \leq j < j'$ and a strong immersion embedding $\eta$ from $(D_j',r_j')$ to $(D_{j'}',r'_{j'})$ such that $\phi'_j(v) \preceq_3 \phi'_{j'}(\eta(v))$ for every $v \in V(D_j')$.
By the definition of $\phi'_i$, there exist injections $\iota_L$, $\iota_M$ and $\iota_R$ from $\{A: (A,B) \in \Se_j$ with $\lvert V(A)-V(B) \rvert \geq \min\{t,2\}\}$ to $\{A': (A',B') \in \Se_{j'}$ with $\lvert V(A')-V(B') \rvert \geq 2\}$ such that for every $(A,B) \in \Se_j$ with $\lvert V(A)-V(B) \rvert \geq \min\{t,2\}$, 
	\begin{itemize}
		\item $\eta(v_{A,L})=v_{\iota_L(A),L}$, $\eta(v_{A,M})=v_{\iota_M(A),M}$, $\eta(v_{A,R})=v_{\iota_R(A),R}$, and 
		\item there exist 
			\begin{itemize}
				\item a strong immersion embedding $\eta_{A,L}$ from $(A_{X_A},v_{A,0},v_{{A,1}_{Y_A}})$ to $(\iota_L(A)_{X_{\iota_L(A)}}, \allowbreak v_{\iota_L(A),0}, \allowbreak v_{{\iota_L(A),1}_{Y_{\iota_L(A)}}})$ such that $\phi_j(v) \leq_Q \phi_{j'}(\eta_{A,L}(v))$ for every $v \in X_A$, 
				\item a strong immersion embedding $\eta_{A,M}$ from $(A,v_{A,0},v_{A,1})$ to $(\iota_M(A), \allowbreak v_{\iota_M(A),0},v_{\iota_M(A),1})$ such that $\phi_j(v) \leq_Q \phi_{j'}(\eta_{A,M}(v))$ for every $v \in V(A)$, and 
				\item a strong immersion embedding $\eta_{A,R}$ from $(A_{Y_A},v_{{A,0}_{X_A}},v_{A,1})$ to $(\iota_R(A)_{Y_{\iota_L(A)}}, \linebreak v_{{\iota_R(A),0}_{X_{\iota_R(A)}}}, \allowbreak v_{\iota_R(A),1})$ such that $\phi_j(v) \leq_Q \phi_{j'}(\eta_{A,R}(v))$ for every $v \in Y_A$. 
			\end{itemize}
	\end{itemize}

For $(A,B) \in \Se_j$ with $\lvert V(A)-V(B) \rvert \geq \min\{t,2\}$, we say that $(A,B)$ is {\it loose} if $\iota_L(A) \neq \iota_R(A)$; otherwise we say that $(A,B)$ is {\it tight}.
Note that if $(A,B)$ is tight, then $\iota_M(A)=\iota_L(A)=\iota_R(A)$ and $\eta(v_{A,M})=v_{\iota_L(A),M}=v_{\iota_R(A),M}$.

Define $\eta^*$ to be a function whose domain is the union of $V(D_j)$ and a subset of $E(D_j)$ such that the following hold.
	\begin{itemize}
		\item If $v \in V(D_j) \cap V(D_j')$, then define $\eta^*(v)=\eta(v)$.
		\item If $v \in V(A)-V(B)$ for some $(A,B) \in \Se_j$ with $\lvert V(A)-V(B) \rvert \geq \min\{t,2\}$, then
			\begin{itemize}
				\item if $(A,B)$ is tight, then define $\eta^*(v) = \eta_{A,M}(v)$, and
				\item if $(A,B)$ is loose, then
					\begin{itemize}
						\item if $v \in X_A-\{v_{A,0}\}$, then define $\eta^*(v)= \eta_{A,L}(v)$, and
						\item if $v \in Y_A-\{v_{A,1}\}$, then define $\eta^*(v)=\eta_{A,R}(v)$.
					\end{itemize}
			\end{itemize}
		\item If $e \in E(D_j) \cap E(D_j')$ and $\eta(e) \subseteq D_{j'}$, then define $\eta^*(e)=\eta(e)$.
		\item If $e \in E(A)$ with both ends in $X_A-\{v_{A,0}\}$ for some loose $(A,B) \in \Se_j$, then define $\eta^*(e)=\eta_{A,L}(e)$.
		\item If $e \in E(A)$ with both ends in $Y_A-\{v_{A,1}\}$ for some loose $(A,B) \in \Se_j$, then define $\eta^*(e)=\eta_{A,R}(e)$.
		\item If $e \in E(A)$ with both ends in $V(A)-\{v_{A,0},v_{A,1}\}$ for some tight $(A,B) \in \Se_j$, then define $\eta^*(e)=\eta_{A,M}(e)$.
	\end{itemize}
Note that $\phi_j(v) \leq_Q \phi_{j'}(\eta^*(v))$ for every $v \in V(D_j)$.
Recall that $\phi_i(v)$ is incomparable with $\phi_i(r_i)$ for every $i \geq 1$ and $v \in V(D_i)$.
So $\eta^*(r_j)=r_{j'}$.
Moreover, it is straightforward to see that for every directed edge $e$ of $D_j$, say from $u$ to $v$, with $\eta^*(e)$ is defined, $\eta^*(e)$ is a directed path from $\eta^*(u)$ to $\eta^*(v)$ internally disjoint from $\eta^*(V(D_j))$; and for distinct directed edges $e,e'$ of $D_j$ with $\eta^*(e),\eta^*(e')$ defined, $\eta^*(e)$ and $\eta^*(e')$ are edge-disjoint.

To prove this lemma, it suffices to show that we can further define $\eta^*(e)$ for the rest of edges of $e \in E(D_j)$ to extend $\eta^*$ to a strong immersion embedding from $(D_j,r_j)$ to $(D_{j'},r_{j'})$.

Note that every edge $e$ of $D_j$ for which $\eta^*(e)$ was not defined satisfies one of the following:
	\begin{itemize}
		\item[(i)] $e \in E(D_j) \cap E(D_{j'})$ and $\eta(e) \not \subseteq D_{j'}$.
		\item[(ii)] $e \in E(A)$ for some loose $(A,B) \in \Se_j$, and $e$ is between $X_A$ and $Y_A$.
		\item[(iii)] $e \in E(A)$ for some loose $(A,B) \in \Se_j$, and $e$ is incident with exactly one of $v_{A,0}$ or $v_{A,1}$. 
		\item[(iv)] $e \in E(A)$ for some tight $(A,B) \in \Se_j$, and $e$ is incident with exactly one of $v_{A,0}$ or $v_{A,1}$.
		\item[(v)] $e \in E(A)$ for some tight $(A,B) \in \Se_j$, and $e$ is from $v_{A,0}$ to $v_{A,1}$.
	\end{itemize}

For each $(A',B') \in \Se_{j'}$ with $\lvert V(A')-V(B') \rvert \geq \min\{t,2\}$, we say that an edge $e \in E(D_j')$ is {\it $(A',B')$-free} if $\eta(e)$ contains all vertices in $\{v_{A',L},v_{A',M},v_{A',R}\}$ as internal vertices.

\medskip

\noindent{\bf Claim 1:} For each $(A',B') \in \Se_{j'}$ with $\lvert V(A')-V(B') \rvert \geq \min\{t,2\}$, there exists an injection from the set of $(A',B')$-free edges of $D_j'$ to a set of edge-disjoint directed paths in $A' \subseteq D_{j'}'$ from $v_{A',0}$ to $v_{A',1}$ internally disjoint from $\eta^*(V(D_j))$.

\noindent{\bf Proof of Claim 1:}
We may assume that the set of $(A',B')$-free edges of $D_j'$ is nonempty, for otherwise we are done.
Hence $\{v_{A',L},v_{A',M},v_{A',R}\}$ is disjoint from $\eta(V(D_j'))$.
So $V(A')-V(B')$ is disjoint from $\eta^*(D_j)$.
Let $k_{A'}$ be the number of edges of $D'_{j'}$ between $v_{A',L}$ and $v_{A',M}$.
So there are $k_{A'}$ edge-disjoint paths in $A'$ from $v_{A',0}$ to $v_{A',1}$ by the definition of the $(\Se_{j'},t)$-compression.
In addition, $\eta$ maps each $(A',B')$-free edge of $D_j'$ to a path containing an edge between $v_{A',L}$ and $v_{A',M}$, so there are at most $k_{A'}$ $(A',B')$-free edges.
Therefore, there exists an injection from the set of $(A',B')$-free edges of $D_j'$ to a set of edge-disjoint directed paths in $A'$ from $v_{A',0}$ to $v_{A',1}$ internally disjoint from $\eta^*(V(D_j))$.
$\Box$

\medskip

For each edge $e$ of $D_j'$, the {\it free-substitution} of $e$ is the directed path in $D_{j'} \cup D_{j'}'$ obtained from $\eta(e)$ by for each $(A',B') \in \Se_{j'}$ in which $e$ is $(A',B')$-free, replacing the subpath $v_{A',0}v_{A',L}v_{A',M}v_{A',R}v_{A',1}$ of $\eta(e)$ by the image of $e$ under the injection mentioned in Claim 1.
Note that the free-substitution of $e$ can be edge-partitioned into a directed path in $D_{j'}$ and at most two directed paths in $D_{j'}'-E(D_{j'})$.

For each edge $e$ of $D_j$ satisfying (i), define $\eta^*(e)$ to be the free-substitution of $e$. 
Since $e$ satisfies (i), $\eta^*(e)$ is a directed path in $D_{j'}$.
Then it is clear that the function $\eta^*$ defined so far does not violate any condition for being a strong immersion embedding from $(D_j,r_j)$ to $(D_{j'},r_{j'})$.
It remains to define $\eta^*(e)$ for edges $e$ satisfying (ii), (iii), (iv) or (v).

For each $(A,B) \in \Se_j$ with $\lvert V(A)-V(B) \rvert \geq \min\{t,2\}$, we define the following:
	\begin{itemize}
		\item if $(A,B)$ is loose, then let $\eta_{A,X}=\eta_{A,L}$ and $\eta_{A,Y}=\eta_{A,R}$, 
		\item if $(A,B)$ is tight, then let $\eta_{A,X}=\eta_{A,M}$ and $\eta_{A,Y}=\eta_{A,M}$, 
		\item let $\pi_{A,L}$ be a bijection from the set of edges of $A$ incident with $v_{A,0}$ to the set of edges of $D'_j$ between $v_{A,0}$ and $v_{A,L}$, (note that this bijection exists since these two sets have the same size by the definition of the $(\Se_j,t)$-compression), and 
		\item let $\pi_{A,R}$ be a bijection from the set of edges of $A$ incident with $v_{A,1}$ to the set of edges of $D'_j$ between $v_{A,1}$ and $v_{A,R}$. 
	\end{itemize}

Now we define $\eta^*(e)$ for edges $e$ satisfying (iii), (iv) or (v).
For each edge $e \in E(A)$ for some $(A,B) \in \Se_j$ with $\lvert V(A)-V(B) \rvert \geq \min\{t,2\}$, 
	\begin{itemize}
		\item if $e$ is incident with $v_{A,0}$ but not incident with $v_{A,1}$, then define $\eta^*(e)$ to be the directed path in $D_{j'}$ obtained from the free-substitution of $\pi_{A,L}(e)$ by replacing the edge $v_{\iota_L(A),0}v_{\iota_L(A),L}=v_{\iota_L(A),0}\eta(v_{A,L})$ by the directed path $\eta_{A,X}(e)$, 
		\item if $e$ is incident with $v_{A,1}$ but not incident with $v_{A,0}$, then define $\eta^*(e)$ to be the directed path in $D_{j'}$ obtained from the free-substitution of $\pi_{A,R}(e)$ by replacing the edge $v_{\iota_R(A),R}v_{\iota_R(A),1}=\eta(v_{A,R})v_{\iota_R(A),1}$ by the directed path $\eta_{A,Y}(e)$, and
		\item if $e$ is from $v_{A,0}$ to $v_{A,1}$ and $(A,B)$ is tight, then define $\eta^*(e)$ to be the directed path in $D_{j'}$ obtained by concatenating the following three directed paths in $D_{j'}$:
			\begin{itemize}
				\item the directed path obtained from the free-substitution of $\pi_{A,L}(e)$ by deleting the vertex $v_{\iota_L(A),L}=\eta(v_{A,L})$, 
				\item $\eta_{A,M}(e)$, and
				\item the directed path obtained from the free-substitution of $\pi_{A,R}(e)$ by deleting the vertex $v_{\iota_R(A),R}=\eta(v_{A,R})$. 
			\end{itemize}
	\end{itemize}
Then it is clear that the function $\eta^*$ defined so far does not violate any condition for being a strong immersion embedding from $(D_j,r_j)$ to $(D_{j'},r_{j'})$.
It suffices to define $\eta^*(e)$ for edges $e$ satisfying (ii).

A {\it middle path} in $D'_j$ is a directed path of the form $v_{A,L}v_{A,M}v_{A,R}$ for some $(A,B) \in \Se_j$ with $\lvert V(A)-V(B) \rvert \geq \min\{t,2\}$.
For each loose $(A,B) \in \Se_j$, there exists a bijection $\pi_{A,M}$ from the set of edges of $A$ between $X_A$ and $Y_A$ and a set of edge-disjoint middle paths contained in $D'_j[v_{A,L},v_{A,M},v_{A,R}]$, as these two sets have the same size.

For each $(A',B') \in \Se_{j'}$ with $\lvert V(A')-V(B') \rvert \geq \min\{t,2\}$, we say that a middle path $P$ in $D'_j$ is {\it $(A',B')$-semifree} if $\eta$ maps $P$ to a directed path containing both $v_{A',L}$ and $v_{A',R}$ as internal vertices.

\medskip

\noindent{\bf Claim 2:} For each $(A',B') \in \Se_{j'}$ with $\lvert V(A')-V(B') \rvert \geq \min\{t,2\}$, there exists an injection from the set of $(A',B')$-semifree internal paths in $D_j'$ to a set of edge-disjoint directed paths in $A' \subseteq D_{j'}'$ from $v_{A',0}$ to $v_{A',1}$ internally disjoint from $\eta^*(V(D_j))$.

\noindent{\bf Proof of Claim 2:}
We may assume that the set of $(A',B')$-semifree middle paths of $D_j'$ is nonempty, for otherwise we are done.
Hence $\{v_{A',L},v_{A',R}\}$ is disjoint from $\eta(V(D_j'))$.
So $V(A')-V(B')$ is disjoint from $\eta^*(D_j)$.
Let $k_{A'}$ be the number of edges of $D'_{j'}$ between $v_{A',L}$ and $v_{A',M}$.
So there are $k_{A'}$ edge-disjoint paths in $A'$ from $v_{A',0}$ to $v_{A',1}$ by the definition of the $(\Se_{j'},t)$-compression.
In addition, $\eta$ maps each $(A',B')$-semifree middle path of $D_j'$ to a path containing an edge between $v_{A',L}$ and $v_{A',M}$, so there are at most $k_{A'}$ $(A',B')$-semifree middle paths.
Therefore, there exists an injection from the set of $(A',B')$-semifree middle paths in $D_j'$ to a set of edge-disjoint directed paths in $A'$ from $v_{A',0}$ to $v_{A',1}$ internally disjoint from $\eta^*(V(D_j))$.
$\Box$

\medskip

For each middle path $P$ with edges $e_1,e_2$ in $D_j'$, the {\it semifree-substitution} of $P$ is the directed path in $D_{j'} \cup D_{j'}'$ obtained from $\eta(e_1) \cup \eta(e_2)$ by for each $(A',B') \in \Se_{j'}$ in which $P$ is $(A',B')$-semifree, replacing the subpath $v_{A',0}v_{A',L}v_{A',M}v_{A',R}v_{A',1}$ of $\eta(e_1) \cup \eta(e_2)$ by the image of $P$ under the injection mentioned in Claim 2.

If $e$ is an edge of $A$ between $X_A$ and $Y_A$ for some loose $(A,B) \in \Se_j$, then $\eta_{A,X}(e)$ contains exactly one edge between $X_{\iota_L(A)}$ and $Y_{\iota_L(A)}$.
Since there are $k_{A'}$ edge-disjoint paths in $A'$ between $v_{A',0}$ and $v_{A',1}$, where $A'=\iota_L(A)$ and $k_{A'}$ is the number of edges of $A'$ between $X_{A'}$ and $Y_{A'}$, we can extend $\eta_{A,X}(e)$, for each $e$ of $A$ between $X_A$ and $Y_B$, by concatenating a path in $A'[Y_{A'}]$ to obtain a directed path $P_{e,L}$ from $\eta_{A,X}(v_e)$ to $v_{A',1}$, where $v_e$ is the tail of $e$, such that if $e_1,e_2$ are distinct edges of $A$ between $X_A$ and $Y_A$, then $P_{e_1,L}$ and $P_{e_2,L}$ are edge-disjoint.

Similarly, for each loose $(A,B) \in \Se_j$ and edge $e$ of $A$ between $X_A$ and $X_B$, we can concatenate $\eta_{A,Y}(e)$ with a path in $A'[X_{A'}]$, where $A'=\iota_R(A)$, to obtain a directed path $P_{e,R}$ from $v_{A',0}$ to $\eta_{A,Y}(v_e)$, where $v_e$ is the head of $e$ such that if $e_1,e_2$ are distinct edges of $A$ between $X_A$ and $Y_A$, then $P_{e_1.R}$ and $P_{e_2,R}$ are edge-disjoint.

For each loose $(A,B) \in \Se_j$ and each $e$ of $A$ between $X_A$ and $Y_A$, we define the following:
	\begin{itemize}
		\item Define $P_e$ to be the directed path obtained by concatenating the following three directed paths:
			\begin{itemize}
				\item $P_{e,L}$,
				\item the directed path obtained from the semifree-substitution of $\pi_{A,M}(e)$ by deleting $v_{\iota_L(A),L}, v_{\iota_L(A),M}, v_{\iota_L(A),R}, v_{\iota_R(A),L}, v_{\iota_R(A),M},v_{\iota_R(A),R}$, and
				\item $P_{e,R}$.
			\end{itemize}
		\item If $e$ is not incident with $v_{A,0}$ nor $v_{A,1}$, then define $\eta^*(e)=P_e$.
		\item If $e$ is from $v_{A,0}$ to $v_{A,1}$, then define $\eta^*(e)$ to be the directed path obtained by concatenating the following three directed paths:
			\begin{itemize}
				\item the directed path obtained from the free-substitution of $\pi_{A,L}(e)$ by deleting $v_{\iota_L(A),L}$,
				\item $P_e$, and
				\item the directed path obtained from the free-substitution of $\pi_{A,R}(e)$ by deleting $v_{\iota_R(A),R}$.
			\end{itemize}
	\end{itemize}
Note that it extends the domain of $\eta^*$ to be a set containing $E(D_j)$. 
It is clear that $\eta^*$ is a strong immersion embedding from $(D_j,r_j)$ to $(D_{j'},r_{j'})$.
This completes the proof.
\end{pf}

\begin{theorem} \label{wqo no loop}
For every positive integer $t$, $\F_t^*$ is well-behaved.
\end{theorem}

\begin{pf}
By Lemma \ref{2-conn wqo}, $\F_t'$ is well-behaved.
By Lemma \ref{2-conn is suff}, $\F_{t,t+1}$ is well-behaved.
By Lemma \ref{root conn suff}, $\F_t^*$ is well-behaved.
\end{pf}

\bigskip

Now we are ready to prove Theorem \ref{main_label}. 
The following is a restatement.

\begin{corollary} \label{wqo_final}
Let $(Q,\leq_Q)$ be a well-quasi-order.
Let $t$ be a positive integer.
For each $i \in {\mathbb N}$, let $D_i$ be a digraph with loops allowed such that there exists no $t$-alternating path in $D_i$, and let $\phi_i: V(D_i) \rightarrow Q$ be a function.
Then there exist $1 \leq j < j'$ and a strong immersion embedding $\eta$ from $D_j$ to $D_{j'}$ such that $\phi_j(v) \leq_Q \phi_{j'}(\eta(v))$ for every $v \in V(D_j)$.
\end{corollary}

\begin{pf}
Let $(Q',\preceq)$ be the well-quasi-order obtained by the Cartesian product of $(Q,\leq_Q)$ and $({\mathbb N} \cup \{0\}, \leq)$.
For each $i \in {\mathbb N}$, let $D_i'$ be the digraph obtained from $D_i$ by deleting all loops, and let $\phi'_i:V(D_i') \rightarrow Q'$ be the function such that for every $v \in V(D'_i)$, $\phi'_i(v)=(\phi_i(v),\ell_v)$, where $\ell_v$ is the number of loops of $D_i$ incident with $v$.
By Theorem \ref{wqo no loop}, there exist $1 \leq j <j'$ and a strong immersion embedding $\eta$ from $D_j'$ to $D'_{j'}$ such that $\phi_j'(v) \preceq \phi'_{j'}(\eta(v))$ for every $v \in V(D_j')$.
Hence for every $v \in V(D_j)$, the number of loops of $D_j$ incident with $v$ is at most the number of loops of $D_{j'}$ incident with $\eta(v)$.
Therefore, we can extend $\eta$ to be a strong immersion embedding from $D_j$ to $D_{j'}$ such that $\phi_j(v) \leq_Q \phi_{j'}(\eta(v))$ for every $v \in V(D_j)$.
\end{pf}

\bigskip

\bigskip

\noindent{\bf Acknowledgement:}
The authors thank Paul Wollan for some discussions about this work. 
This work was initiated when the first author visited the National Center for Theoretical Sciences (NCTS) in Taiwan at 2017.
The first author thanks NCTS for the hospitality.

\end{document}